\documentclass[11pt,oneside, reqno]{amsart}
\usepackage{epsfig,psfrag,subfigure}
\usepackage{amsmath,amssymb,amsfonts,latexsym, mathrsfs,euscript}
\usepackage{graphicx,graphics,cite}
\usepackage{wrapfig}
\usepackage{slashbox}
\usepackage[usenames]{color}
\usepackage{geometry,pbox}                
\usepackage{algorithm,algorithmicx,algpseudocode}

\usepackage[table]{xcolor}
\usepackage{epstopdf}

\newtheorem*{conjecture}{Conjecture}

\newcommand{\ignore}[1]{}

\newcommand{\beq}{\begin{equation}}
\newcommand{\eeq}{\end{equation}}

\theoremstyle{remark}
\newtheorem{note}{Note}

\title[Numerical renormalization group algorithms]{Numerical renormalization group algorithms for self-similar solutions of partial differential equations}

\author{Gast\~ao A. Braga}
\address{Department of Mathematics, Federal University of Minas Gerais, Brazil}
\email{gdabraga@gmail.com}

\author{Frederico C. Furtado}
\address{Department of Mathematics and Statistics, University of Wyoming, Laramie, WY 82071-3036}
\email{furtado@uwyo.edu\\ TEL: 307-766-4321}

\author{Vincenzo Isaia}
\address{Department of Mathematics and Computer Science, Indiana State University, Terre Haute, IN 47809-1902}
\email{Vin.Isaia@indstate.edu\\ TEL: 812-237-2157}

\author{Long Lee}
\address{Department of Mathematics and Statistics, University of Wyoming, Laramie, WY 82071-3036}
\email{llee@uwyo.edu\\ TEL: 307-766-4368}

\date{}
\begin{document}
\begin{abstract}
We systematically study a numerical procedure that reveals the asymptotically self-similar dynamics of solutions of partial differential equations (PDEs). This procedure, based on the renormalization group (RG) theory for PDEs, appeared initially in a conference proceeding by Braga et al. \cite{BFI04}. This numerical version of RG method, dubbed as the numerical RG (nRG) algorithm, numerically rescales the temporal and spatial variables in each iteration and drives the solutions to a fixed point exponentially fast, which corresponds to the self-similar dynamics of the equations. In this paper, we carefully examine and validate this class of algorithms by comparing the numerical solutions with either the exact or the asymptotic solutions of the model equations in literature. The other contribution of the current paper is that we present several examples to demonstrate that this class of nRG algorithms can be applied to a wide range of PDEs to shed lights on longtime self-similar dynamics of certain physical models that are difficult to analyze, both numerically and analytically.



\end{abstract}
\maketitle

\begin{description}
\item[{\footnotesize\bf keywords}]
{\footnotesize  Renormalization group, asymptotic limit, self-similar dynamics, critical exponents}
\end{description}

\section{Introduction}
Solutions of certain partial differential equations (PDEs) exhibit self-similar dynamics in the asymptotic limit. A self-similar solution in the region of large-time asymptotic limit normally is in the form of 
\begin{equation}\label{eq:structure}
u(x, t) \sim \frac{A}{t^{\alpha}}\phi\left(B\frac{x}{t^{\beta}}\right), \quad\text{as}\,\,t\rightarrow\infty.
\end{equation}
This particular form of solution indicates that in the large-time asymptotic region the dynamics of the solution are controlled by two factors, the decay in the magnitude of $u$ and the spread of the spatial distribution of $u$, while the profile of the self-similar solution is described by the function $\phi$. The constant $A$ and $B$ in Eq. (\ref{eq:structure}) are usually related to some conservation laws of the PDE under study. In the literature, physicists were able to predict or determine the scaling exponents $\alpha$ and $\beta$ for critical phenomena by using the renormalization group (RG) approach for a variety of physical models in equilibrium statistical mechanics or quantum field theory. \cite{GL54,G92,SO98, W71-I, W71-II}. In the early 1990s, Goldenfeld et al. developed a perturbative renormalization group method for PDEs and applied it to the study of a number of large-time asymptotic problems \cite{CGO91, GMOL90, GMO92, G92}. A slight twist of the original method was reported in \cite{C97, C01, MC03}. Almost at the same period, a nonperturbative RG approach was introduced by Bricmont et al. \cite{BK95, BKL94}, and was applied to the study of nonlinear dispersive and dissipative wave equations and thermal-diffusive combustion system \cite{BPW94, BKX96, P02}.

Also in the 1990s, Chen and Goldenfeld proposed a numerical RG (nRG) calculation for similarity solutions of porous medium (Barenblatt) equation and traveling waves \cite{CG95}. Their numerical procedure inspired mesh renormalization schemes for studying focusing problems arising in porous medium flow \cite{AABL01, BAA00}. Inspired by the numerical approach of Chen and Goldenfeld and the nonperturbative RG approach of Bricmont et al., Braga et al., introduced a class of  nRG algorithms in a short conference  paper \cite{BFI04} that allow them to systematically search for the critical exponents and the hidden decay in asymptotically self-similar dynamics through repeated scalings in time and space. In this paper, we carefully examine and validate the nRG algorithms of Braga et al. by comparing the numerical solutions of the nRG algorithms with either the exact or the asymptotic solutions of the model equations in literature. We show that the self-similar dynamics captured by the nRG algorithms agree with the theoretical results for the scalar equations, as well as the system of equations. 
Furthermore, a novel contribution of this paper is that we demonstrate that this nRG procedure could shed light on the behavior of self-similar solutions of certain physical models, such as the nonlinear diffusion-absorption model, that are difficulty to analyze, both numerically and analytically.

It is worth noting that, a procedure introduced by Braga et al. for studying a nonlinear diffusion equation with periodic coefficients in \cite{BFMR03} shares the similar spirit as an nRG algorithm studied in this paper.
Numerical procedures based on rescaling the solutions and the time and spacial variables have previously been developed and were used to study the solutions of PDE that blow up in finite time \cite{BK88, FW03, LPSS86, RW00}. Such procedures exploit the known self-similar structure of the solutions under study to determine the appropriate rescalings. The nRG algorithms in this paper, however,  are unique in exploiting fixed points by generating successive iterations of a discrete RG transformation in space and time that drive the system towards a fixed point, and the current implementation of the algorithms in this paper is not suitable for studying blow-up problems. A numerical procedures that mimics the renormalization group theory to compute the spatial profile and blow-up time for self-similar behavior was proposed by Isaia \cite{bib:blow_up}. This version of nRG algorithm, however, uses the Berger-Kohn's time step \cite{BK88} that assumes the time decay exponent is known. Modification of the nRG algorithms presented in this paper for studying the blow-up problems is currently under our investigation and will be reported in a separate paper. 

Finally, to aid the reader, we now outline the contents of the remainder of this paper. In Section \ref{sec:Burgers}, using the Burgers equation, we explain the fundamental idea of the nRG algorithm and outline the procedure of the algorithm. We show that the algorithm allows us for detecting and further finding the critical exponents in the self-similar solutions. In particular, we demonstrate the ability of the algorithm for determining the critical scaling exponents in time and space that render explicitly the distinct physical effects of the solutions of the Burgers equation, depending on the initial conditions. In Section \ref{sec:KdV}, we use the phenomena of dispersive shock waves of the Korteweg-de Vries equation to show that the algorithm can be used as a time integrator for investigating intermediate asymptotic behavior of solutions. In Section \ref{sec:diffusion_absorption}, we study a class of nonlinear diffusion-absorption models. A conjecture on the existence of a critical exponent of the nonlinear absorption term is proposed for problems with discontinuous diffusivities. We also present a marginal case, for which the phenomenon of anomalous decay is observed, as a motivation for the next section. Finally, we present a modified nRG algorithm and illustrate the ability of the modified algorithm for detecting and capturing the hidden logarithmic decay through a nonlinear system of cubic autocatalytic chemical reaction equations in Section \ref{sec:cubic_autocatalytic}.


\section{Scaling Transformation for Burgers Equation}\label{sec:Burgers}

The nRG method studied in this paper is simply the integration of the PDE over a finite time-interval with fixed length followed by a rescaling. To explain this idea, we use the Burgers equation as an example to illustrate the scaling transformation procedure of nRG algorithms. We compare the asymptotic solutions obtained by the nRG algorithm with the exact solutions in the asymptotic region to demonstrate the robust nature of our algorithm. 

The Burgers equation with initial data at $t=1$ are written as 
\begin{equation}\label{eq:burgers}
\begin{split}
&u_t+uu_x=\nu u_{xx},\,\,t >1,\\
\text{I. C.\,\,:}\quad&u(x,1) = f(x),
\end{split}
\end{equation}
where $\nu > 0$ is the viscosity. Let the time and space variables be scaled by powers of a fixed length  $L>1$,
\begin{equation}\label{eq:scale_tx}
t = L\tilde{t},\quad x=L^{\beta}\tilde{x},
\end{equation}
where $\beta> 0$, $\tilde{t}$ and $\tilde{x}$ are new variables.
Suppose the solution of the initial value problem (IVP) (\ref{eq:burgers}), $u(x, t)$, is scaled by
\begin{equation}\label{eq:scale_u}
u_L(\tilde{x},\tilde{t}) = L^{\alpha}\,u(x,t) = L^{\alpha}\,u(L^{\beta}\tilde{x},L\tilde{t}),
\end{equation}
where $\alpha > 0$. This implies that  
\begin{equation}\label{eq:scale_u2}
u(x, t) =L^{-\alpha}u_{L}(\tilde{x},\tilde{t}) = L^{-\alpha}u_{L}(L^{-\beta}x, L^{-1}t). 
\end{equation}
With the above scalings, each term in the Burgers equation is scaled as follows:
\begin{equation}\label{eq:u_t}
u_t =  L^{-(\alpha+1)}\frac{\partial}{\partial \tilde{t}}\left(u_{L}(\tilde{x}, \tilde{t})\right),\,\, u_x =  L^{-(\alpha+\beta)}\frac{\partial}{\partial \tilde{x}}\left(u_{L}(\tilde{x}, \tilde{t})\right),\,\, u_{xx} = L^{-(\alpha+2\beta)}\frac{\partial^{2}}{\partial \tilde{x}^2}\left(u_{L}(\tilde{x}, \tilde{t})\right).
\end{equation}
Substituting Eq. (\ref{eq:u_t}) into Eq. (\ref{eq:burgers}) yields 
\begin{equation}\label{eq:s_burgers}
\begin{split}
&L^{-(\alpha+1)}(u_L)_{\tilde{t}}+L^{-(2\alpha+\beta)}u_L(u_L)_{\tilde{x}}=\nu L^{-(\alpha+2\beta)} (u_L)_{\tilde{x}\tilde{x}},\,\,\tilde{t} >1,\\
\text{I. C.\,\,:}\quad&u_L(\tilde{x},1) = \tilde{f}(\tilde{x}),
\end{split}
\end{equation}
We rewrite the above equation as 
\begin{equation}\label{eq:s2_burgers}
\begin{split}
&(u_L)_{\tilde{t}}+L^{-\alpha-\beta+1}u_L(u_L)_{\tilde{x}}=\nu L^{-2\beta+1} (u_L)_{\tilde{x}\tilde{x}},\,\,\tilde{t} >1,\\
\text{I. C.\,\,:}\quad&u_L(\tilde{x},1) = \tilde{f}(\tilde{x}).
\end{split}
\end{equation}
The integration length for time is from $\tilde{t}=1$ to $\tilde{t}=L$, while the transformed initial condition is $\tilde{f}(\tilde{x})=L^{\alpha}u(L^{\beta}\tilde{x}, L)$.

\subsection{Sequence of Scaling Transformations}\label{sec:seq}
If we perform a sequence of scalings (iterations), then with a fixed $L>1$, and sequences of scaling exponents $\{\alpha_n\}$ and $\{\beta_n\}$, we can define a sequence of rescaled functions $\{u_n\}$ by rewriting Eq. (\ref{eq:scale_u}) (dropping $\tilde{}$ in $\tilde{x}$ and $\tilde{t}$ ) as
\begin{equation}\label{eq:scale_un}
u_n(x,t)  = L^{\alpha_n}\,u_{n-1}(L^{\beta_n}x,Lt),
\end{equation} 
with $u_0=u$ of the original IVP,  Eq. (\ref{eq:burgers}). 
A simple calculation reveals that 
\begin{equation}\label{eq:u_seq}
\begin{split}
u_n(x,t)&=L^{(\alpha_1+\alpha_2+\alpha_3+\cdots+\alpha_n)}u(L^{(\beta_1+\beta_2+\beta_3+\cdots+\beta_n)}x, L^nt)\\
             &=L^{n\bar{\alpha}_n}u(L^{n\bar{\beta}_n}x, L^nt),
\end{split}
\end{equation}
where $\bar{\alpha}_n = \frac{1}{n}(\alpha_1+\alpha_2+\cdots+\alpha_n)$ and $\bar{\beta}_n = \frac{1}{n}(\beta_1+\beta_2+\cdots+\beta_n)$.   Eq. (\ref{eq:u_seq}) shows how $u_n$ in the time interval $t \in [1, L]$ is related to $u$ in the time interval in $t\in [L^{n}, L^{n+1}]$.
Since at each iteration, the scaling of the PDE shown in Eq.(\ref{eq:s2_burgers}) is applied to the previous scaled equation,  the solution of the $n^{th}$ iteration, $u_n(x,t)$, is the solution of the following scaled initial value problem 
\begin{equation}\label{eq:s_burgers_un}
\begin{split}
&(u_n)_{t}+L^{n(-\bar{\alpha}_n-\bar{\beta}_n+1)}u_n(u_n)_{x}=\nu L^{n(-2\bar{\beta}_n+1)} (u_n)_{xx},\,\,t >1,\\
\text{I. C.\,\,:}\quad&u_n(x,1) = f_n(x),
\end{split}
\end{equation}
where $f_n(x)=L^{\alpha_n}u_{n-1}(L^{\beta_n}x, L)$, with $f_0=f$, the initial condition of  the IVP (\ref{eq:burgers}).

A simple observation for the scaled IVP  (\ref{eq:s_burgers_un}) is that if $\bar{\alpha}_n\rightarrow1/2$ and $\bar{\beta}_n\rightarrow1/2$ as $n\rightarrow\infty$, the diffusion term has unscaled viscosity and the pre-factor of the quasilinear term in the left-hand-side is of order 1 for sufficiently large $n$. In this case, for small viscosity, the long-time behavior of the solution is dominated by the quasilinear term. On the other hand, if $\bar{\alpha}_n\rightarrow 1$ and  $\bar{\beta}_n\rightarrow 1/2$ as $n\rightarrow\infty$, the  viscosity in the diffusion term  remains unscaled, but the advection term on the left-hand-side has $L$ to the power of negative $n/2$. With $L > 1$ and $n\rightarrow\infty$, the advection term eventually drops out the equation, and the diffusion will be the dominant term. The numerical renormalization group calculations based on the algorithm in the next section confirm that for positive mass initial data, it is the former case, while for the zero mass initial data, it is the later one. There is a third possibility, $\bar{\alpha}_n\rightarrow 0$ and $\bar{\beta}_n\rightarrow 1$ as $n\rightarrow\infty$, corresponding to the traveling wave solutions of the Burgers equation, that we do not discuss in this paper.

\subsection{The nRG procedure}\label{sec:nrg}


We describe the nRG procedure in Algorithm \ref{alg1:nrg}.
\begin{algorithm}[t]
\begin{algorithmic}
\For {$n=0,1,2,\ldots,$ until convergence}
\begin{enumerate}
\item[1.] Start with the IVP (\ref{eq:burgers}) for $n=0$. Evolve $u_n$ from $t=1$ to $t=L$, using the IVP (\ref{eq:s_burgers_un})  for $n\ge 1$.
\item[2.] Compute $\alpha_{n+1}$ by
\begin{equation*}
L^{\alpha_{n+1}} = \frac{||u_n(\cdot,1)||_{\infty}}{||u_n(\cdot,L)||_{\infty}} =\frac{||u(L^{n\bar{\beta}_n}\cdot,L^n)||_{\infty}}{||u(L^{n\bar{\beta}_n}\cdot,L^{n+1})||_{\infty}}.
\end{equation*}
\item[3.] Compute $\beta_{n+1}$ from an appropriate scaling relation between $\alpha_{n+1}$ and $\beta_{n+1}$. Normally, $\beta_{n+1}= g(\alpha_{n+1})$, for some function $g$.
\item[4.] Compute $A_n=L^{n(\alpha_n-\bar{\alpha}_n)}$,  $B_n=L^{n(\beta_n-\bar{\beta}_{n})}$, and $f_{n+1}(x)=L^{\alpha_{n+1}}u_n\left(L^{\beta_{n+1}}x,L\right)$, where $\bar{\alpha}_n = \frac{1}{n}(\alpha_1+\alpha_2+\cdots+\alpha_n)$ and $\bar{\beta}_{n} = \frac{1}{n}(\beta_1+\beta_2+\cdots+\beta_n)$.
\end{enumerate}
\EndFor
\end{algorithmic}
\caption{The nRG procedure for the Burgers equation}
\label{alg1:nrg}
\end{algorithm}
Note that in Step 4. of Algorithm \ref{alg1:nrg}, the variable $A_n$ is defined, as now we explain, based on the assumption that for the solution of Eq. (\ref{eq:burgers}), denoted $u(x, t)$,  there exists a self-similar profile function $\phi$ such that 
\begin{equation}\label{eq:rg_similarity}
u(x,t)\sim \frac{A}{t^{\alpha}}\phi\left(B\frac{x}{t^{\beta}}\right),\quad {\text{as}}\,\,\,t\rightarrow\infty,
\end{equation} 
where $A$ and $B$ are constants and $\alpha$ and $\beta$ are the powers of decay and spreading with respect to time, respectively.  After $n$ iterations,  $t=L^{n}$, and hence Eq. (\ref{eq:rg_similarity}) becomes 
\begin{equation}\label{eq:rg_similarity_Ln}
u(x,L^{n})\sim \frac{A}{L^{n\alpha}}\phi\left(B\frac{x}{L^{n\beta}}\right),
\end{equation}
which implies that for some large enough $n$
\begin{equation}\label{eq:rg_similarity_Ln_inv}
L^{n\alpha}u(L^{n\beta}x,L^{n})\sim A\phi\left(Bx\right).
\end{equation}
Suppose that $\alpha_n\rightarrow\alpha$ and $\beta_n\rightarrow\beta$, as $n\rightarrow\infty$, Eq. (\ref{eq:rg_similarity_Ln_inv}) is equivalent to 
\begin{equation}\label{eq:rg_similarity_Ln_inv_n}
L^{n\alpha_n}u(L^{n\beta_n}x,L^{n})\sim A\phi\left(Bx\right),\quad n\rightarrow\infty.
\end{equation}
From Eq. (\ref{eq:u_seq}),  we have 
\beq\label{eq:u_seq_inv}
L^{-n\bar{\alpha}_n}u_n(L^{-n\bar{\beta}_n}x, t) = u(x, L^{n}t).
\eeq
Letting $t=1$ in Eq. (\ref{eq:u_seq_inv}), Eqs. (\ref{eq:rg_similarity_Ln_inv_n}) and (\ref{eq:u_seq_inv}) together imply that 
\beq\label{eq:u_n_phi_1}
L^{n(\alpha_n - \bar{\alpha}_n)}u_n(L^{n(\beta_n-\bar{\beta}_n)}x, 1) \sim A\phi(Bx).
\eeq
If we define $A_n=L^{n(\alpha_n-\bar{\alpha}_n)}$,  $B_n=L^{n(\beta_n-\bar{\beta}_n)}$, we expect that 
\begin{equation}\label{eq:An_Bn}
\lim_{n\rightarrow\infty}A_n=\lim_{n\rightarrow\infty}L^{n(\alpha_n-\bar{\alpha}_n)}\rightarrow A,\quad \lim_{n\rightarrow\infty}B_n=\lim_{n\rightarrow\infty}L^{n(\beta_n-\bar{\beta}_n)} \rightarrow B,
\end{equation}
provided 
\beq\label{eq:u_n_phi_2}
u_n(L^{n(\beta_n-\bar{\beta}_n)}x, 1)\rightarrow \phi(Bx),\,\,\text{as}\,\,\,n\rightarrow\infty.
\eeq

\subsection{Numerical experiments}\label{sec:implement}

We now describe the simple explicit algorithm used for solving the scaled PDE (\ref{eq:s_burgers_un}). Suppose that within each iteration of the nRG algorithm, we discretize the spatial derivatives by the centered difference scheme and use the first-order Euler method for our temporal discretization. If we denote
\begin{equation}\label{eq:k_n}
\kappa_n = L^{n(-\bar{\alpha}_n-\bar{\beta}_n+1)},\quad \nu_n=\nu L^{n(-2\bar{\beta}_n+1)},\quad v=u_n,
\end{equation}
at time $t^{j+1}=(j+1)\Delta t$, the fully discretized scaled equation of Eq. (\ref{eq:s_burgers_un}) at the $i^{th}$ grid point is
\begin{equation}\label{eq:Euler-center}
v_i^{j+1} = v_i^j -\Delta t\frac{\kappa_n}{2}\left(\frac{(v_{i+1}^j)^2-(v_{i-1}^j)^2}{2\Delta x_n}\right) + \Delta t \nu_n\left (\frac{v_{i-1}^j-2v_i^j+v_{i+1}^j}{\Delta x_n^2}\right ),
\end{equation}
where $\Delta x_n$ is the spatial grid size for the $n^{th}$ iteration. From Eq. (\ref{eq:scale_un}), suppose that 
we denote the $i^{th}$ spatial node on a uniform mesh at the $n^{th}$ iteration as $(x_i)_n=i\Delta x_n$, where $i=0,\pm1,\pm2\cdots$, we have 
\begin{equation}\label{eq:scaled_xi}
u_n(i\Delta x_n,t)=L^{\alpha_n}\,u_{n-1}\left(L^{\beta_n}i\Delta x_{n},Lt\right)=L^{\alpha_n}\,u_{n-1}\left(i\Delta x_{n-1},Lt\right).
\end{equation}
This implies that 
\begin{equation}\label{eq:scaled_dx}
\Delta x_n = L^{-\beta_n}\Delta x_{n-1}. 
\end{equation}
Since $L>1$, if $\beta_n > 0$ for all $n$ then $\Delta x_n < \Delta x_{n-1} < \Delta x_{n-2}<\cdots<\Delta x_0$, and for sufficiently large $n$, $\Delta x_n<<\Delta x_0$. If a uniform $\Delta t$ is used for all iterations, this could lead to  numerical instability at later iterations for an explicit time integrator, such as the Euler method adopted in this paper.  
Conversely, if $\Delta t$ decreases accordingly to maintain the stability requirement, eventually the time integration becomes too costly for the nRG algorithm. A remedy for this situation is that instead of scaling the mesh size, we keep the same $\Delta x$ at all time, i.e. $\Delta x_n=\Delta x_{n-1}= \cdots = \Delta x_0 = \Delta x$, through interpolation. 

To explain this idea, we see that at the end of the first iteration ($n=0$), we are supposed to set the initial data at the $j^{th}$ node for the next iteration to be  $f_1(x_j) = f_1(L^{\beta_1}j\Delta x_1)=L^{\alpha_1}u_0(j\Delta x_0, L)$, or equivalently  $f_1(x_j)= f_1(j\Delta x_1)=L^{\alpha_1}u_0(L^{-\beta_1}j\Delta x_0, L)$. If we choose $x_j$ on the new grid to be the same as that of the old grid, i.e. $j\Delta x_1 = j\Delta x_0$, we need the value $u_0(L^{-\beta_1}j\Delta x_0, L)$.  Since $u_0(j\Delta x_0, L) $ has been computed for each $j$, to obtain $u_0(L^{-\beta_1}j\Delta x_0, L)$, we can simply linearly interpolate $u_0(k\Delta x_0, L)$ and $u_0((k+1)\Delta x_0, L)$, where  $k\Delta x_0< L^{-\beta_1}j\Delta x_0< (k+1)\Delta x_0$ for some $k$. By repeating this process, we have $\Delta x=\Delta x_0=\Delta x_1=\Delta x_2=\cdots = \Delta x_n$. 
The above interpolation principle can be extended to quadratic or cubic interpolation, and spline functions as well. This interpolation-resampling strategy was previously proposed in \cite{SO98} as a means to capture the consequences of space-time translational symmetry on a discrete lattice.
\begin{note}
If an interpolation scheme is used,  the last step in Algorithm \ref{alg1:nrg} is modified by $L^{\alpha_{n+1}}u_{n}\left(L^{\beta_{n+1}}x,L\right) = \displaystyle\frac{1}{\max|\bar{u}_n\left(L^{\beta_{n+1}}x,L\right)|}\bar{u}_{n}\left(L^{\beta_{n+1}}x,L\right)$, where $\bar{u}_{n}\left(L^{\beta_{n+1}}x,L\right)$ is the interpolation of $u_n\left(L^{\beta_{n+1}}x,L\right)$. This normalization step is necessary to avoid diminishing of the amplitude due to the interpolation.
\end{note}

\subsubsection{Positive initial mass}\label{sec:positive}
The first initial condition we consider for our numerical experiments is the characteristic function  
\begin{equation}\label{eq:init}
u(x,0)=\chi_{[-\ell,\ell]}(x) =\begin{cases}
1,  & -\ell \le x \le \ell, \\
0,  &\text{else}.
\end{cases}
\end{equation}
It is known that a conserved quantity of the Burgers equation is the total mass defined by 
\begin{equation}\label{eq:tm}
M=\int_{-\infty}^{\infty}u(x,t)dx. 
\end{equation}
For the characteristic function, the total mass is trivial to compute initially.  

Whitham \cite{whitham74} showed that a special asymptotic self-similar solution of single hump for the Burgers equation with initial data possessing positive total mass has the following explicit formula (page 104, Eq. (4.32)):  
 \begin{equation}\label{eq:whitham} 
 u(x,t) = \sqrt{\frac{2M}{t}} g(z,R),
\end{equation}
where  $M>0$ is the total mass of initial data, $z= x/\sqrt{2Mt}$ is the similarity variable, and $R=M/(2\nu)$ is the Reynolds number,  where $\nu$ is the viscosity. The function $g(z,R)$ is 
\begin{equation}\label{eq:g}
g(z,R) = \frac{(e^R-1)}{2\sqrt{R}}\frac{e^{-z^2R}}{\sqrt{\pi}+(e^R-1)\displaystyle\int_{z\sqrt{R}}^{\infty}e^{-\xi^2}d\xi}.
\end{equation}
Based on dimensional arguments, Whitham indicated that the similarity form of the above solution is 
\begin{equation}\label{eq:similarity}
u(x,t) = \displaystyle\sqrt{\frac{\nu}{t}}\phi\left(\displaystyle\frac{x}{\sqrt{\nu t}};\displaystyle\frac{M}{\nu}\right).
\end{equation}
Comparing Eq. (\ref{eq:rg_similarity}) and Eq. (\ref{eq:similarity}), we expect that the sequences of exponents $\{\alpha_n\}$ and $\{\beta_n\}$ in the nRG calculation converge to $\alpha=1/2$ and $\beta=1/2$, respectively, for sufficiently large iteration numbers. Hence our nRG procedure starts with letting $\beta_n=1/2$ for all $n$. This means that the spatial variable is always scaled by $L^{-1/2}$ for the next iteration. As we mentioned earlier, this assumption also corresponds to ensuring that the viscosity remains the same at all time, that is $\nu_n=\nu$. 

\begin{figure}[bhtp]
\centering
(a)\includegraphics[width=2.8in]{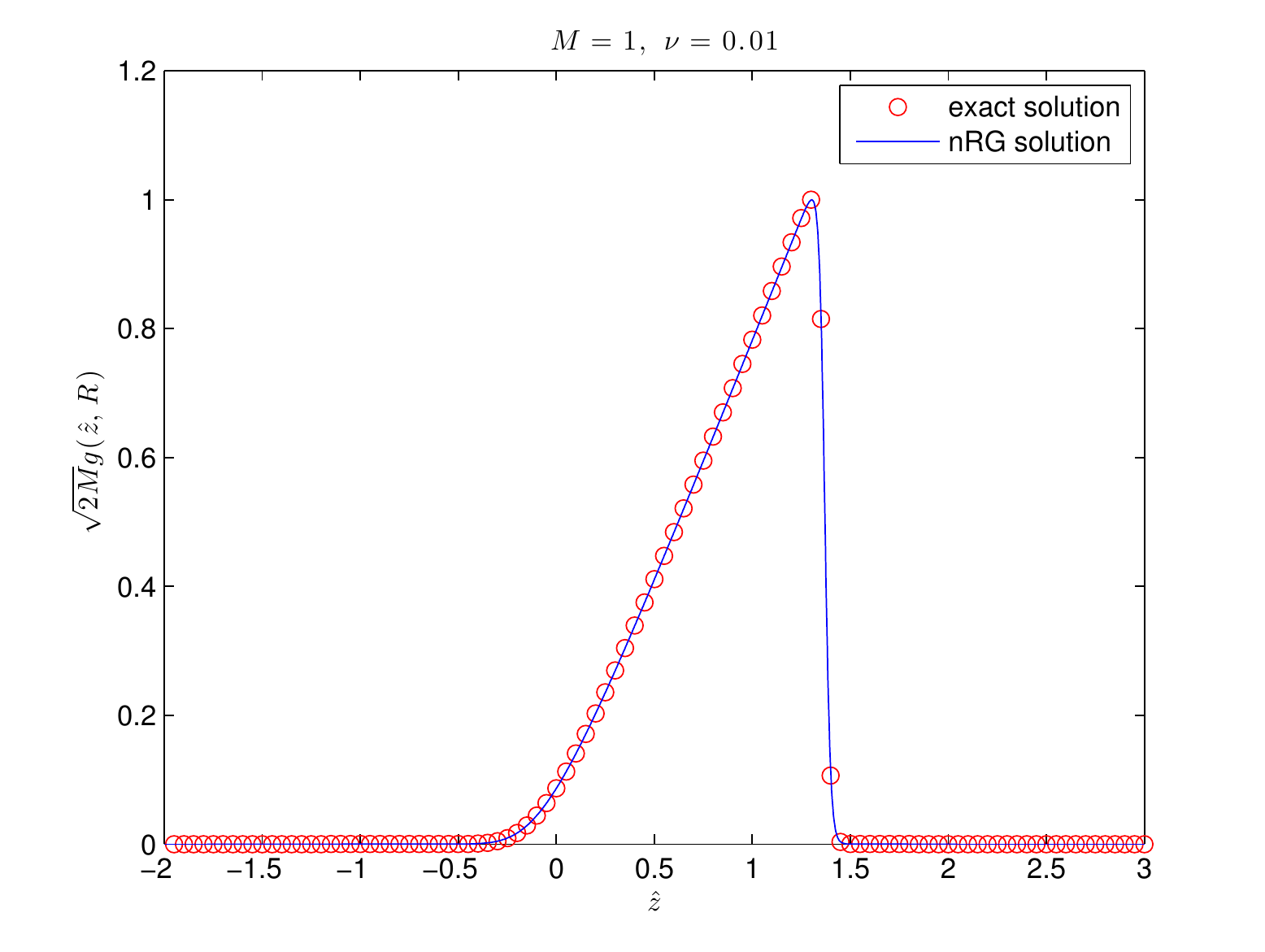}
\includegraphics[width=2.8in]{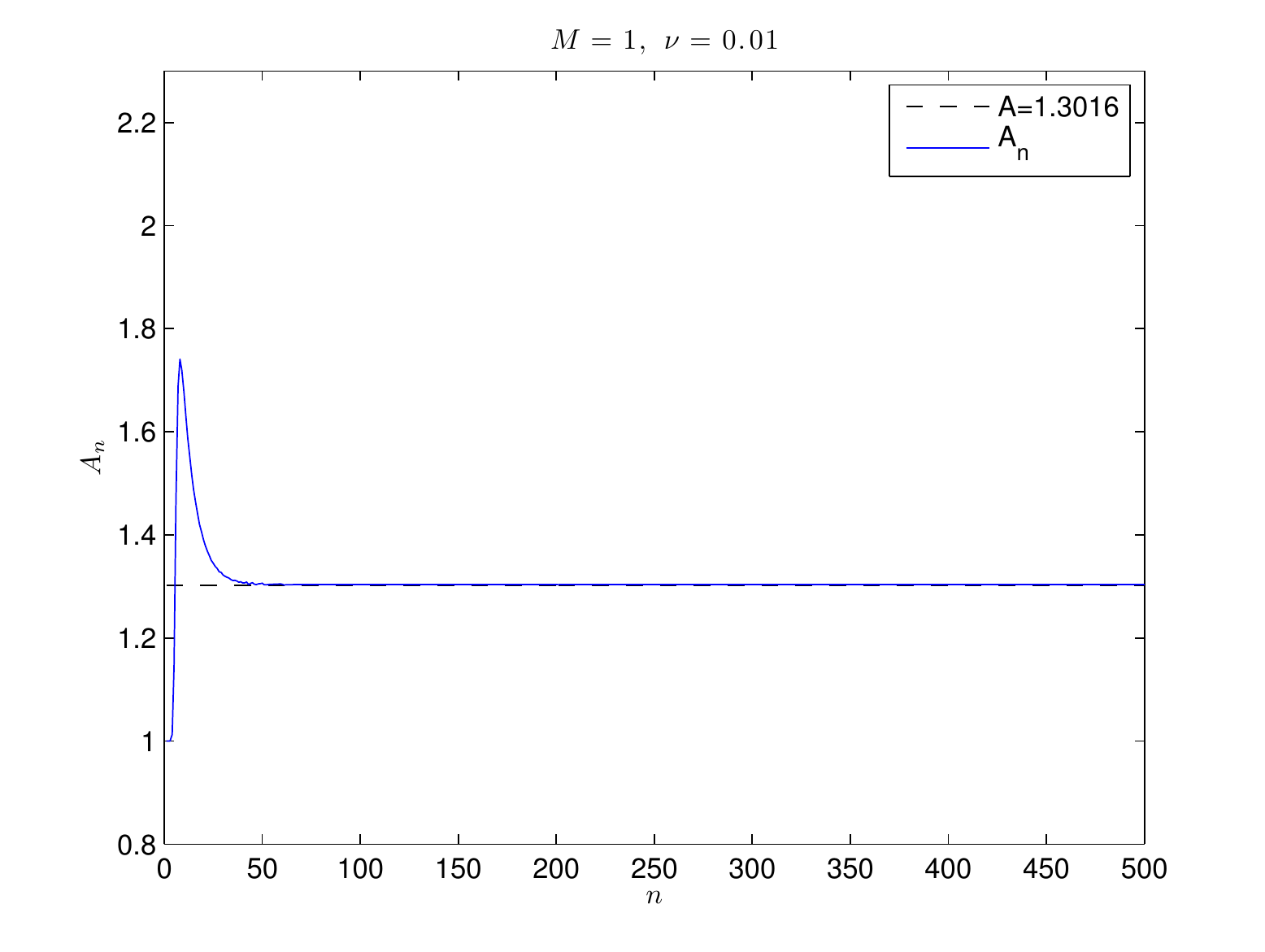}\\
(b)\includegraphics[width=2.8in]{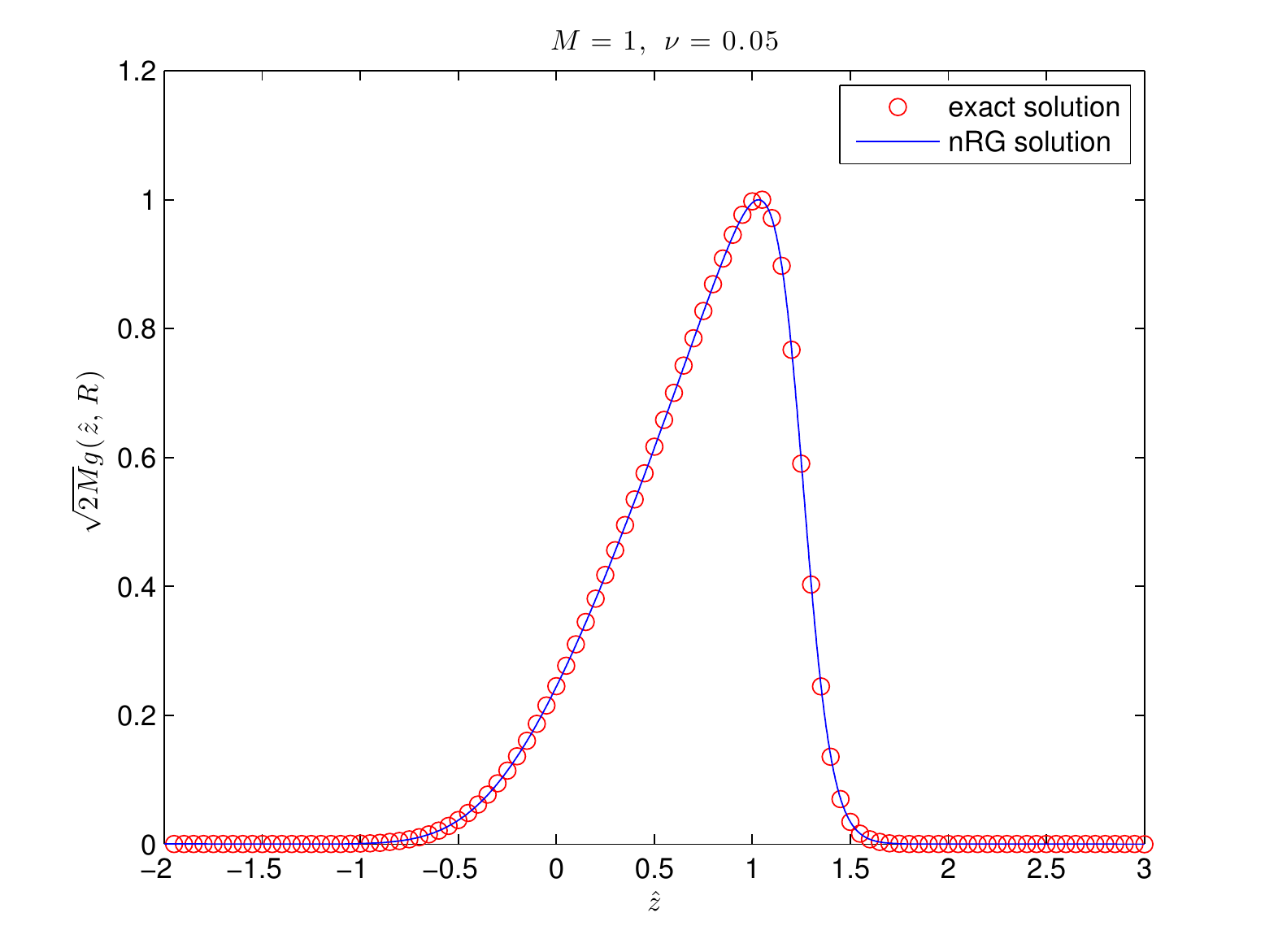}
\includegraphics[width=2.8in]{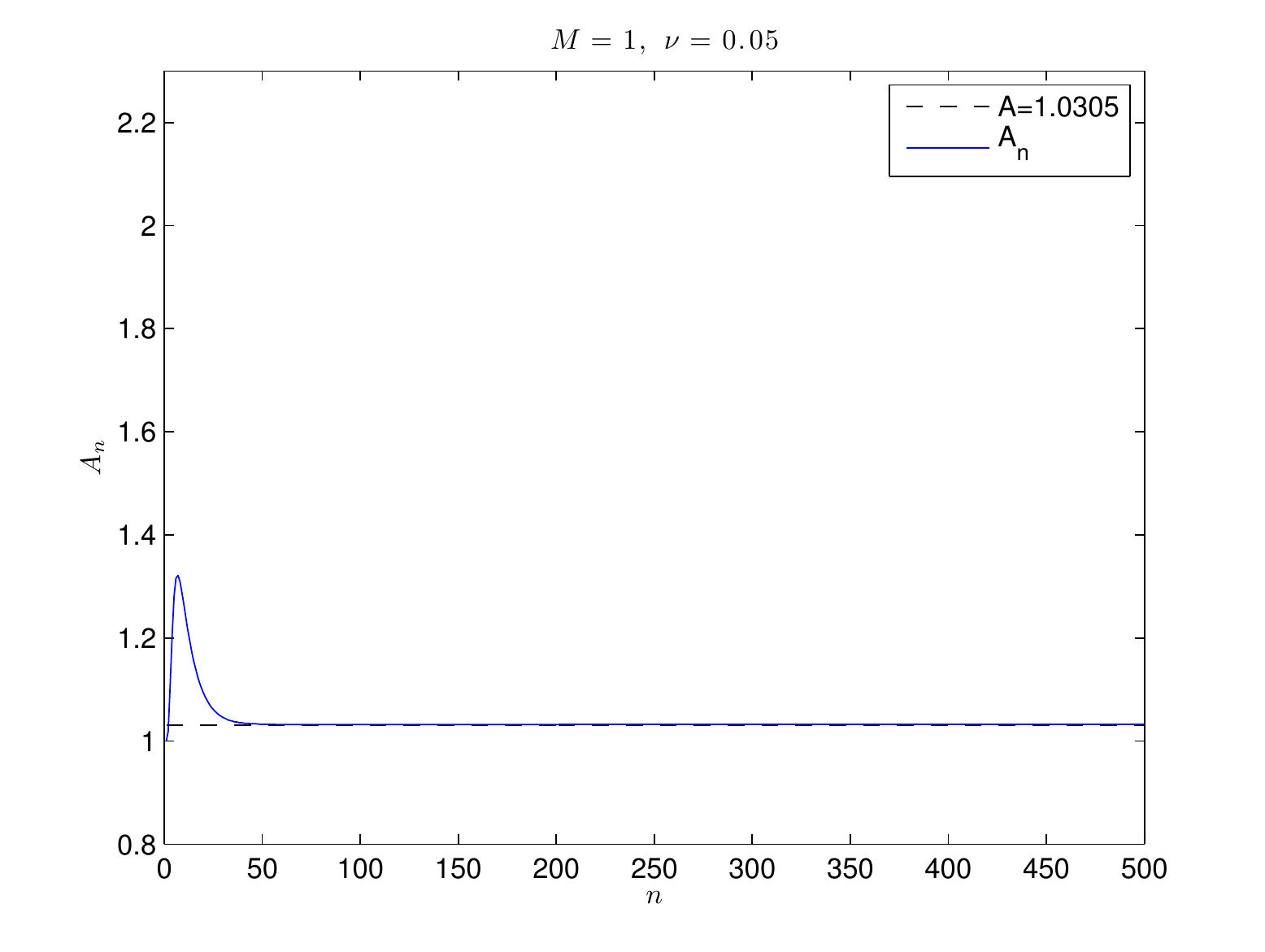}\\
(c)\includegraphics[width=2.8in]{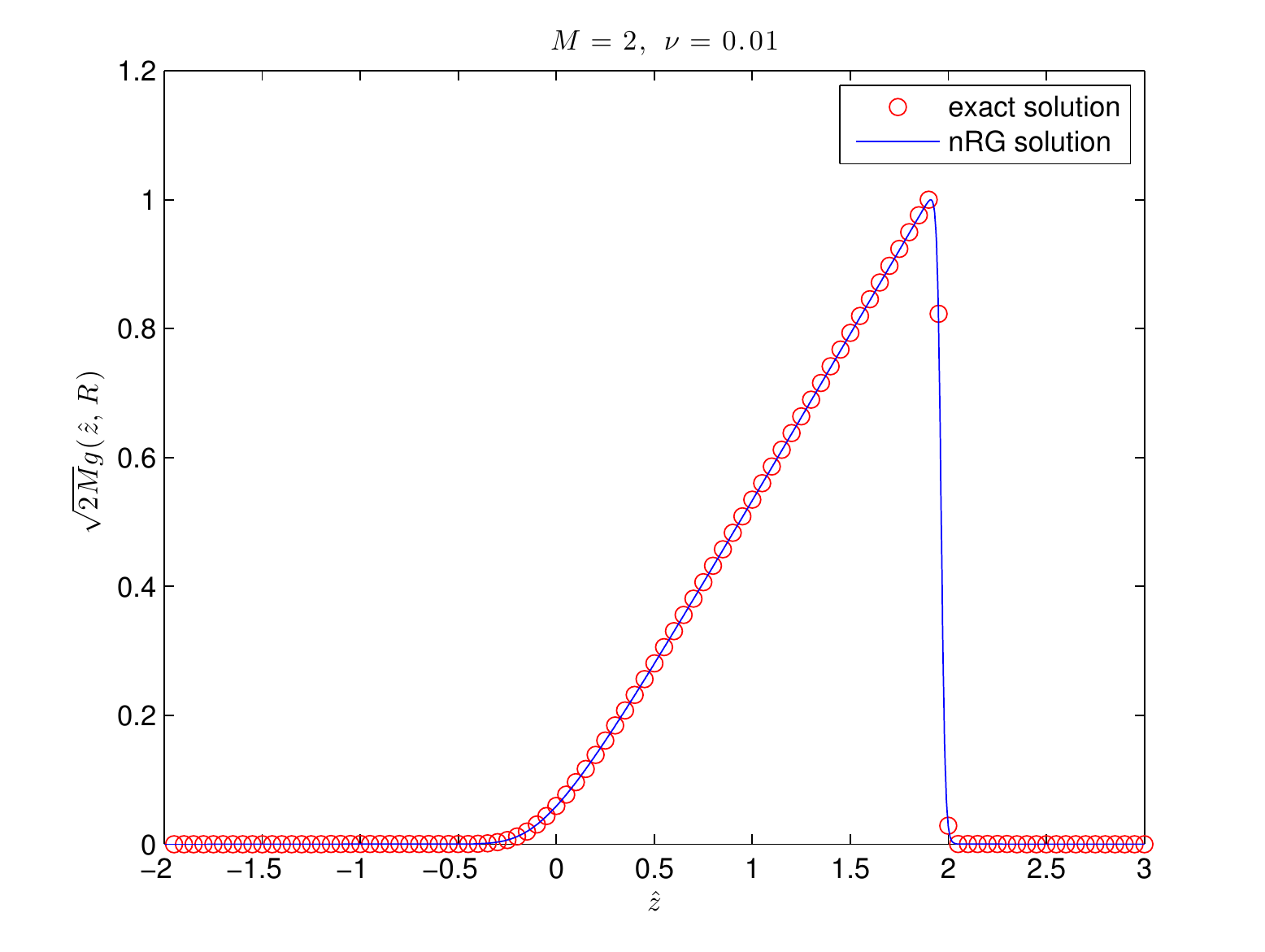}
\includegraphics[width=2.8in]{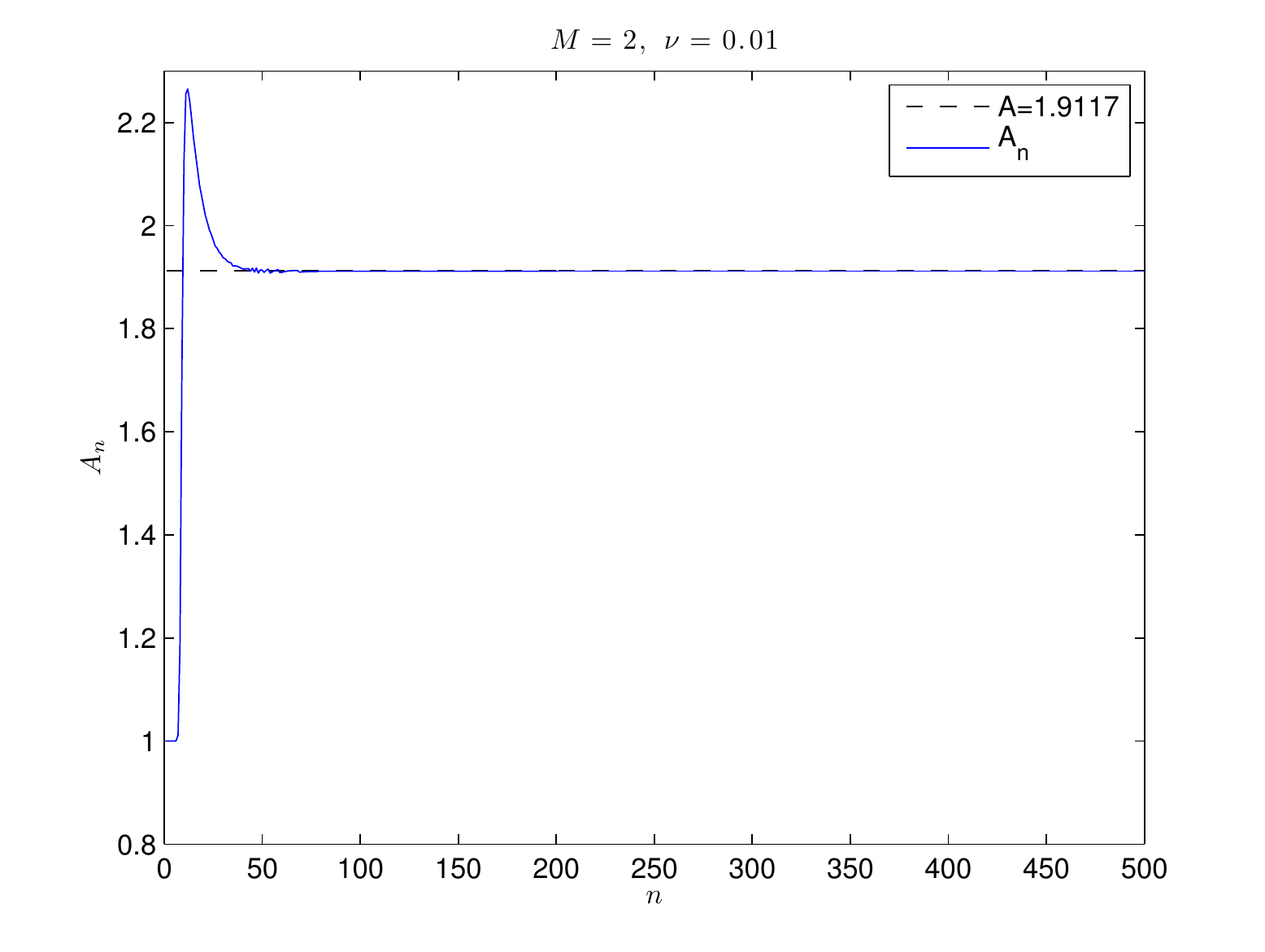}
\caption{Comparisons between the nRG calculations and the asymptotic self-similar solutions shown in \cite{whitham74}. The left panel shows the waveforms of the RG solutions at the final iteration $n=500$ and the theoretic predictions for various $M$ and $\nu$, while the right panel shows the convergence of $A_n$  to its theoretical value $A$. The initial mass and viscosity used for the comparisons are (a) $M=1$, $\nu=0.01$, (b)  $M=1$, $\nu=0.05$, and (c)  $M=2$, $\nu=0.01$. All numerical calculations use $\Delta x=0.0032$ and $\Delta t=0.0001$.}
\label{fig:positive_comp1}
\end{figure}

\begin{figure}[bhtp]
\centering
(a)\includegraphics[width=2.8in]{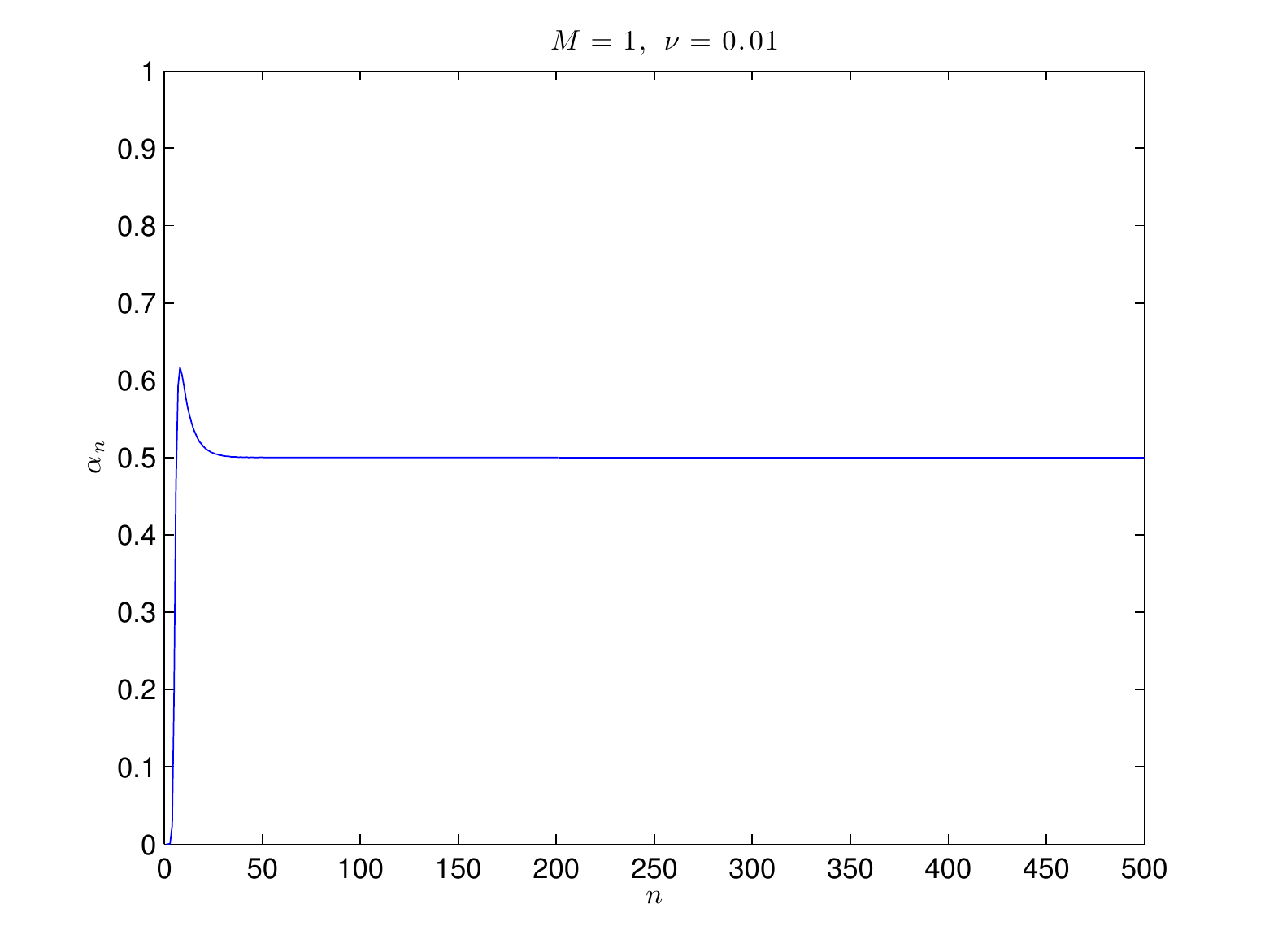}
\includegraphics[width=2.8in]{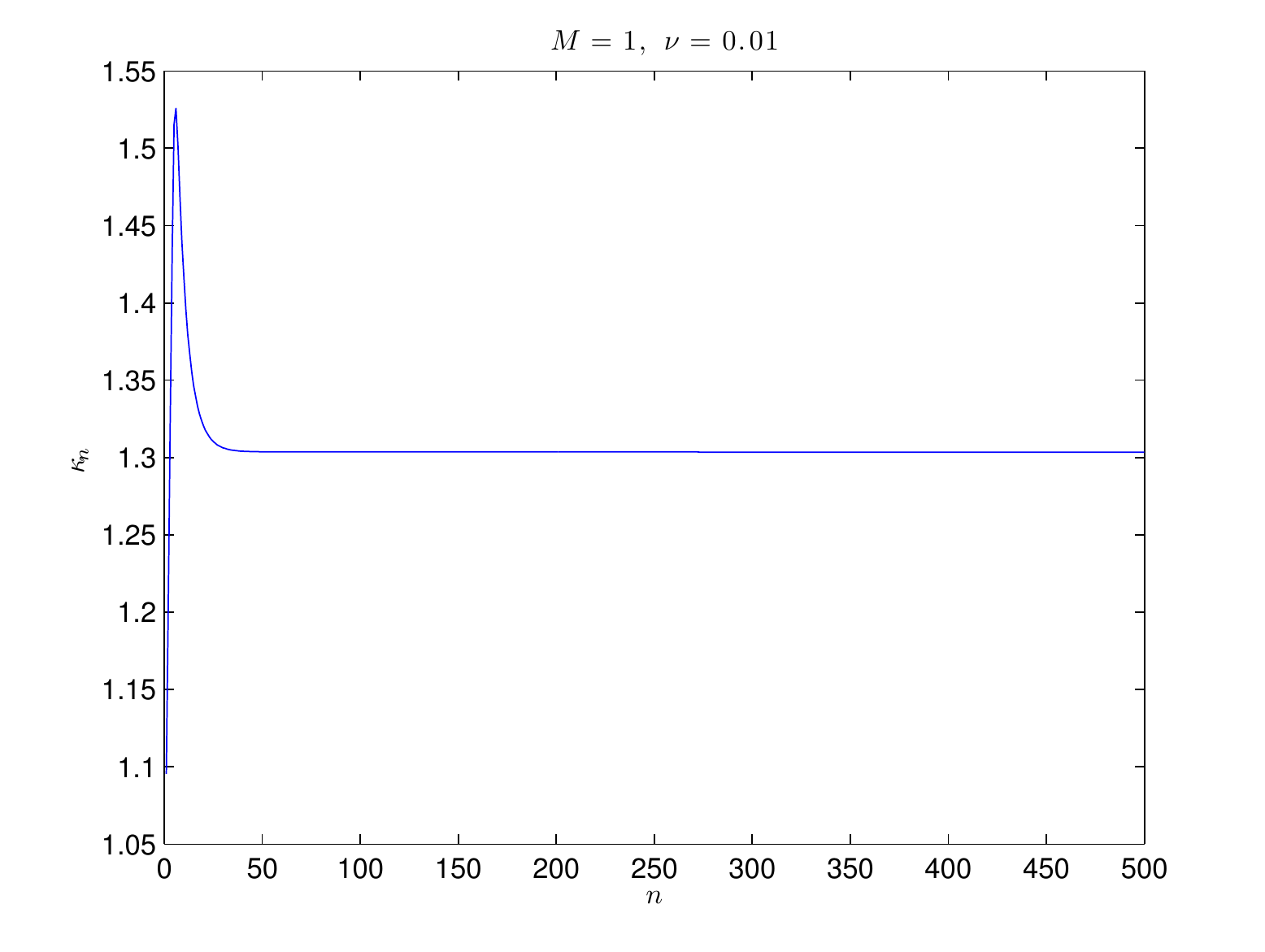}\\
(b)\includegraphics[width=2.8in]{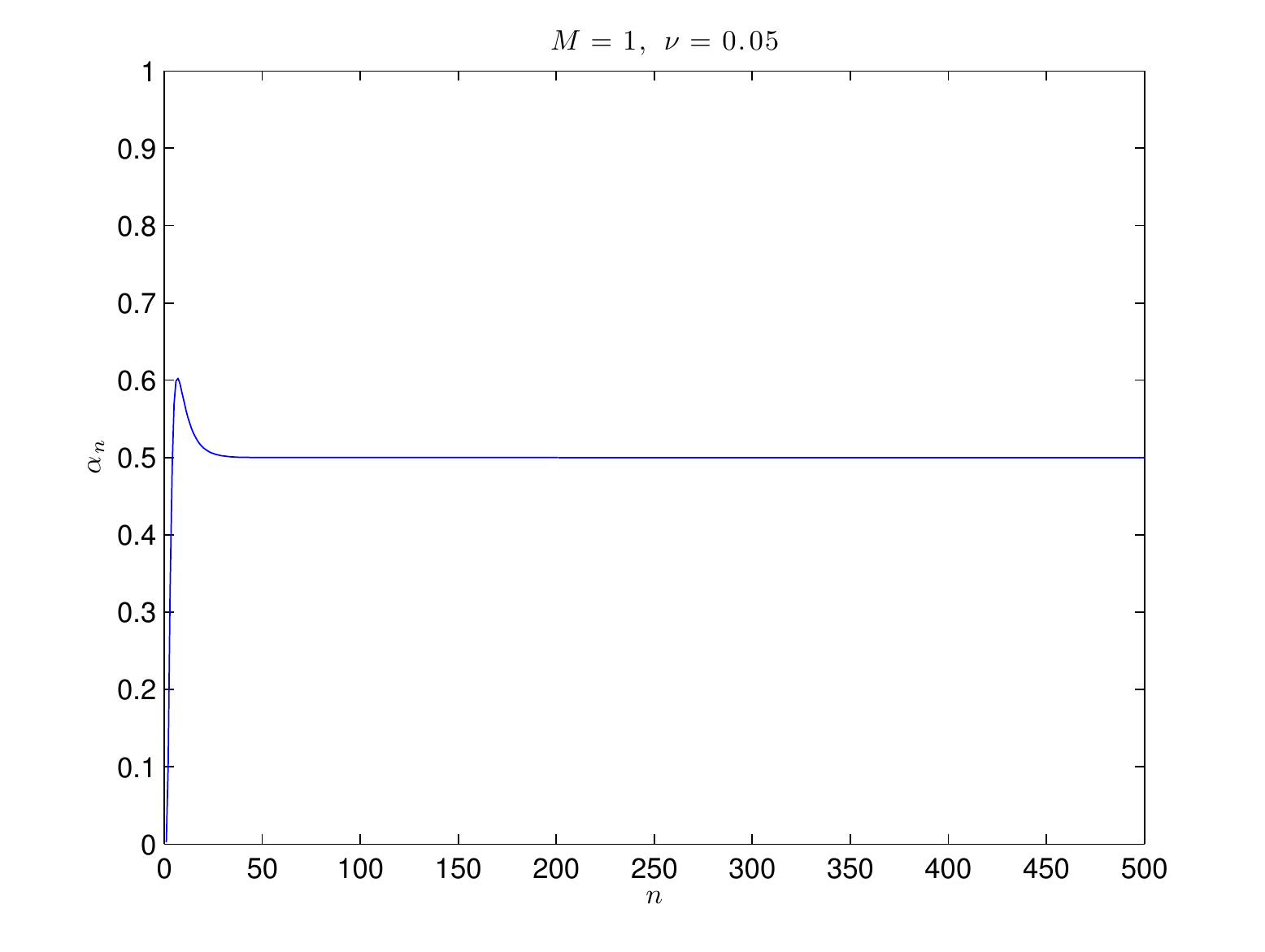}
\includegraphics[width=2.8in]{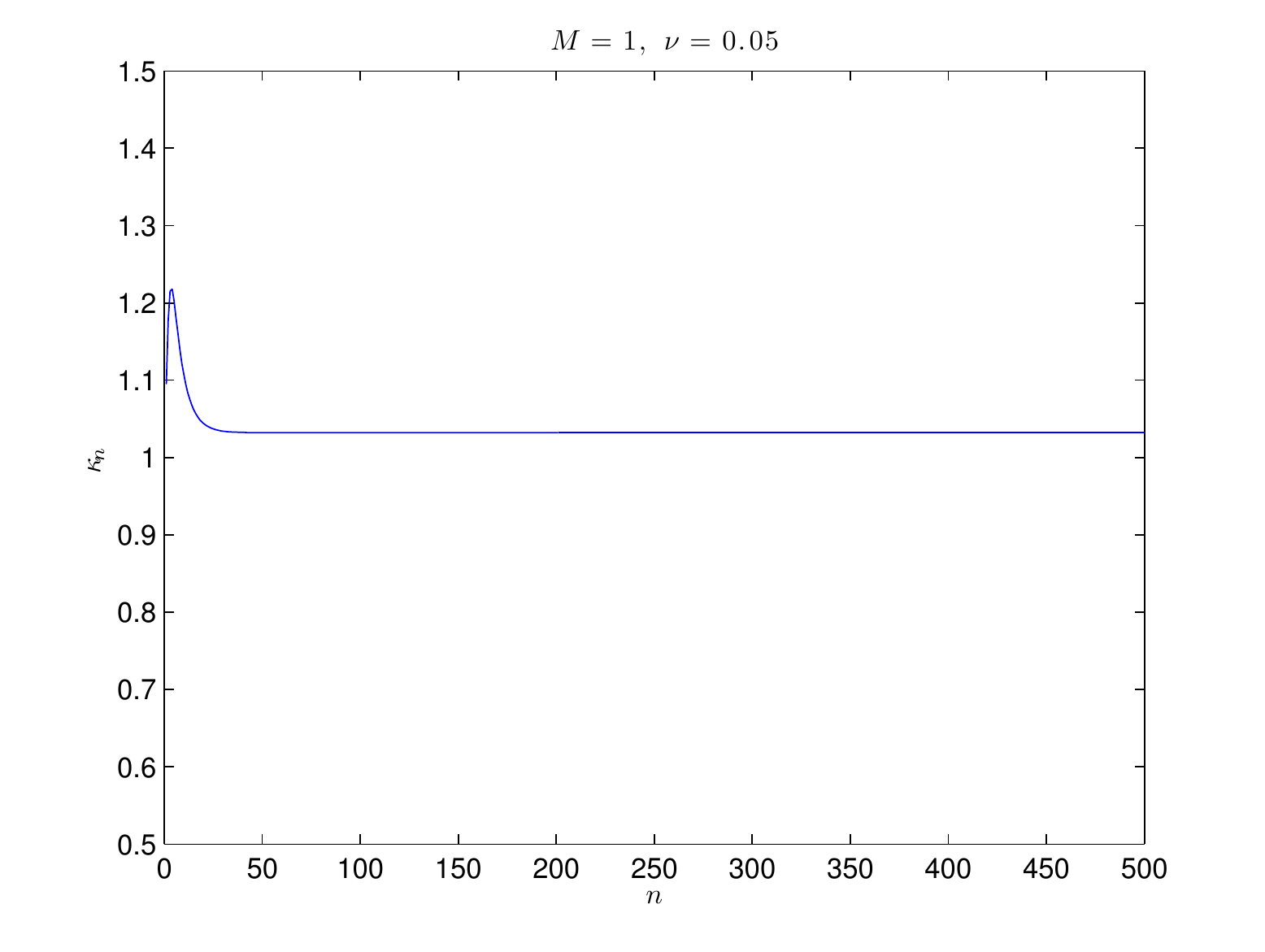}\\
(c)\includegraphics[width=2.8in]{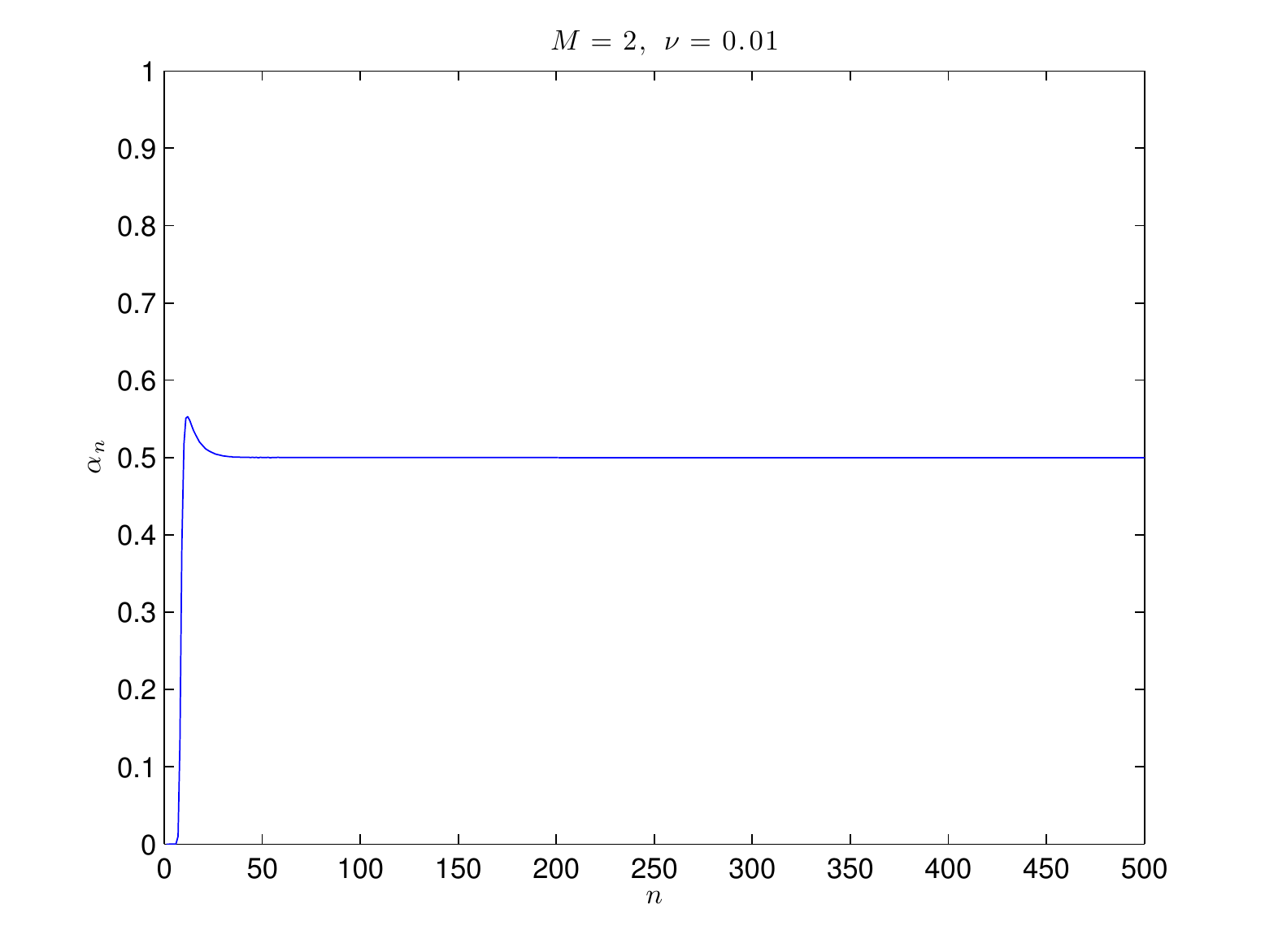}
\includegraphics[width=2.8in]{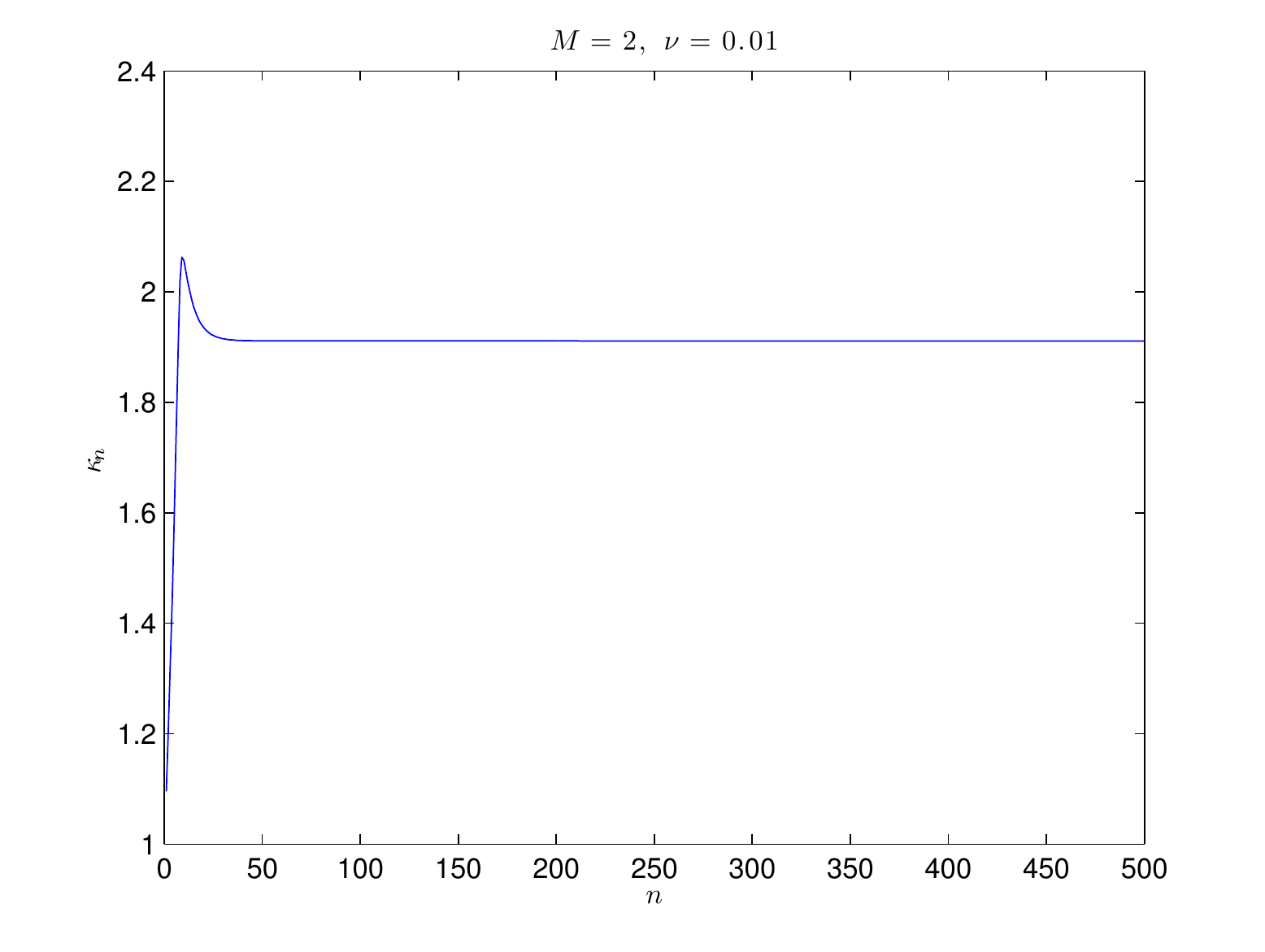}\\
\caption{$\alpha_n$ (left) and $\kappa_n$ (right) versus $n$ (number of iteration). The initial mass and viscosity are the same as that in Figure \ref{fig:positive_comp1} (a), (b), and (c). The figure indicates that $\alpha_n\rightarrow\frac{1}{2}$ and $\kappa_n \ne 0$ and is of $O(1)$, as expected.}
\label{fig:positive_comp2}
\end{figure}

Comparing  Eqs. (\ref{eq:u_seq}) and (\ref{eq:whitham}), we see that the term $\sqrt{t}\,u(\hat{z},t)$ is equivalent to $L^{n\bar{\alpha}_n}u(\hat{z}, L^n)$, where  $\hat{z}=\sqrt{2 M} z = \frac{1}{\sqrt{t}} x = L^{-n/2} x$ with $t=L^n$. With the assumption $\beta_n=1/2$ for all $n$, we have $\bar{\beta}_n=1/2$, and $L^{-n\bar{\beta}_n} = \frac{1}{\sqrt{t}}$. Moreover, if $\{\alpha_n\}$ approaches to $1/2$ (later we will show in our numerical experiments that this is in fact the case), then $\bar{\alpha}_n \approx 1/2$ for sufficiently large $n$. Hence $L^{n\bar{\alpha}_n} \approx \sqrt{t}$. This implies
\begin{equation}\label{eq:2mg}
u_n(L^{-n/2}x,L^n) = L^{n/2} u(x, t) = \sqrt{t}\, u(x,t)=\sqrt{2M}\,g(z,R).
\end{equation}
Now plugging $z=\displaystyle\frac{\hat{z}}{\sqrt{2M}}$ into $g(z,R)$ to obtain $g(\hat{z},R)$, we have
\begin{equation}\label{eq:2mg_hat}
u_n(L^{-n/2}x,L^n) = \sqrt{2M}\,g(\hat{z},R).
\end{equation}
Using Eq. (\ref{eq:2mg_hat}), if we plot the normalized function of $\sqrt{2M}\,g(\hat{z}, R)$ versus $\hat{z}$ (so that the amplitude is one), this should be equivalent to plotting $u_n$ against $L^{-n/2}x$ at the time evolution $t=L^{n}$.  We are now ready to compare the theoretical similarity profile derived by Whitham \cite{whitham74} (page 106, Figure 4.1) and the one obtained from our nRG procedure. Figure \ref{fig:positive_comp1} is the comparison of nRG calculations and the theoretical asymptotic solutions. For each case (a), (b), and (c) of Figure \ref{fig:positive_comp1}, the initial mass is  $M=1,\,1$, and $2$, respectively ($\ell =1/2, 1/2$ and 1, respectively in Eq. (\ref{eq:init})), while the diffusivity constant is $\nu=0.01,\, 0.05$, and $0.01$, respectively. The figure shows that the nRG calculations remarkably capture the theoretical predictions for various initial data and diffusivity constant. In particular, the left panel plots the final waveform of the solution between the predicted theoretical profiles and that of our RG calculations, while the right panel indicates the convergence of $A_n$  to its theoretical value $A$. Note that, integrating Eq. (\ref{eq:u_n_phi_1}), and using Eqs. (\ref{eq:An_Bn}) and (\ref{eq:u_n_phi_2}), the theoretical value $A$ is given by
\beq
A = \displaystyle\frac{M}{\int_{\mathbb{R}}\phi(x)dx},
\eeq
where $\phi(x)$ is the computed RG profile. The total number of iterations for all three cases is $n=500$. Within each iteration, the calculation was carried out in a periodic domain, $-8\le x\le 8$, with $\Delta x= 16/5000$, while the time integration is from $t=1$ to $t=L=1.2$ with $\Delta t=0.2/2000$. 
Figure \ref{fig:positive_comp2} shows the calculations of $\alpha_n$ and $\kappa_n$. The figure indicates that $\alpha_n\rightarrow\frac{1}{2}$ and $\kappa_n \ne 0$ and is of $O(1)$. This suggests that in the asymptotical region the diffusion constant is unscaled. Since the diffusion constant is small and the coefficient of the advection is of order one in all three cases, the equation is advection dominant in the asymptotical region.

\subsubsection{Zero initial mass}

In this example, we consider the following initial condition for the Burgers equation
\begin{equation}\label{eq:init_zero}
u_0(x)=-\chi_{[-\ell, 0]}+\chi_{[\ell, 0]}=\begin{cases}
-1,  & -\ell \le x \le 0, \\
1,  & 0 < x \le \ell, \\
0,  &\text{else}.
\end{cases}
\end{equation}
The total mass of the above function is zero. For this initial condition, Whitham showed that for small diffusivity constants, the inviscid theory is adequate to explain the solutions of the Burgers equation for most of the range, except in the final decay period. The solutions are of typical $N$ wave structures before the final decay. However, as $t\rightarrow\infty$, in the final decay, for any fixed diffusivity constant, no matter how small, the solution of the Burgers equation is 
\begin{equation}\label{eq:diploe}
u(x,t)\sim \frac{x}{t}\sqrt{\frac{a}{t}}e^{-x^2/(4\nu t)},
\end{equation}
for some fixed $a$.
\begin{figure}[hbtp]
\centering
(a)\includegraphics[width=2.8in]{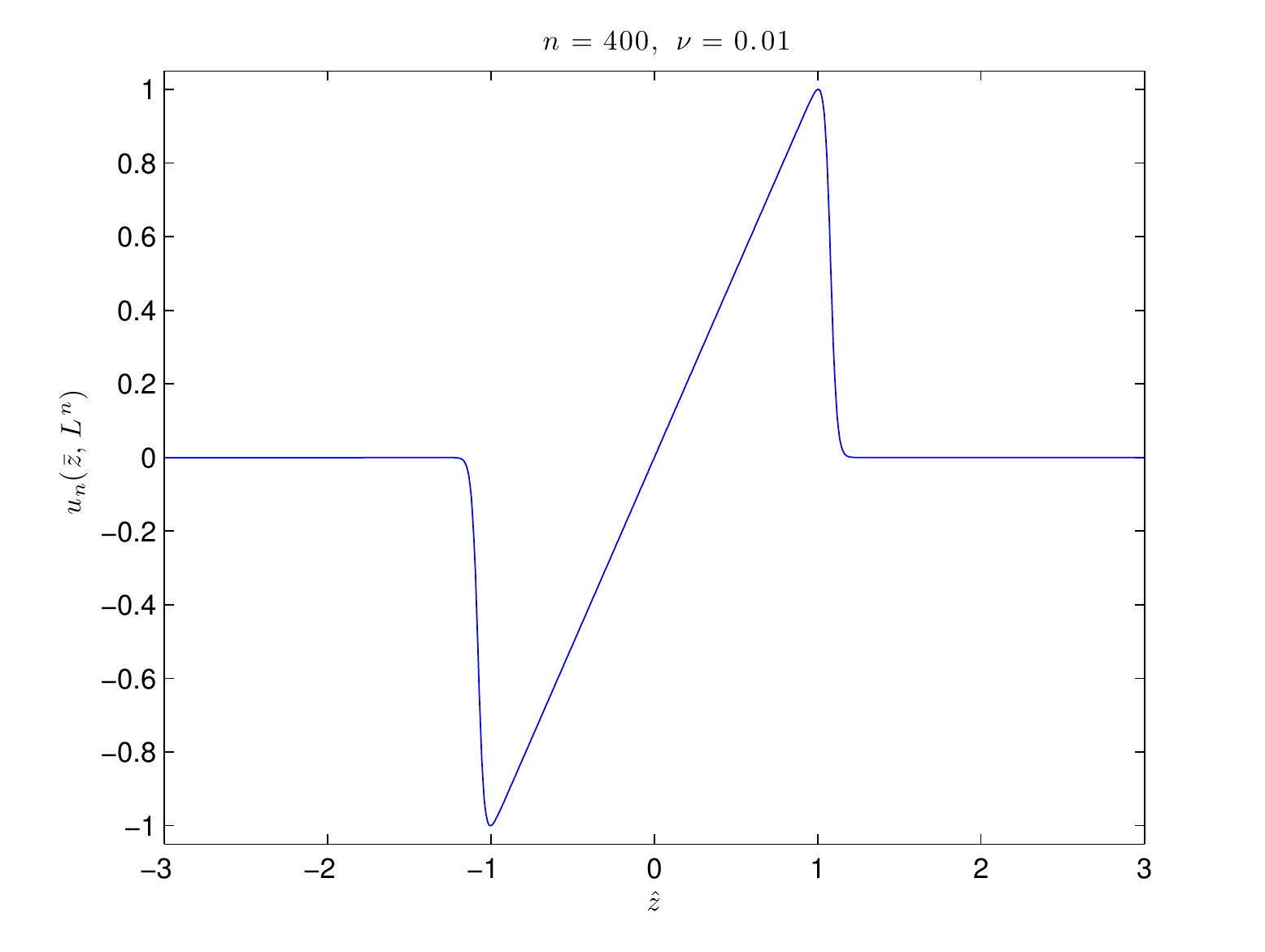} 
(b)\includegraphics[width=2.8in]{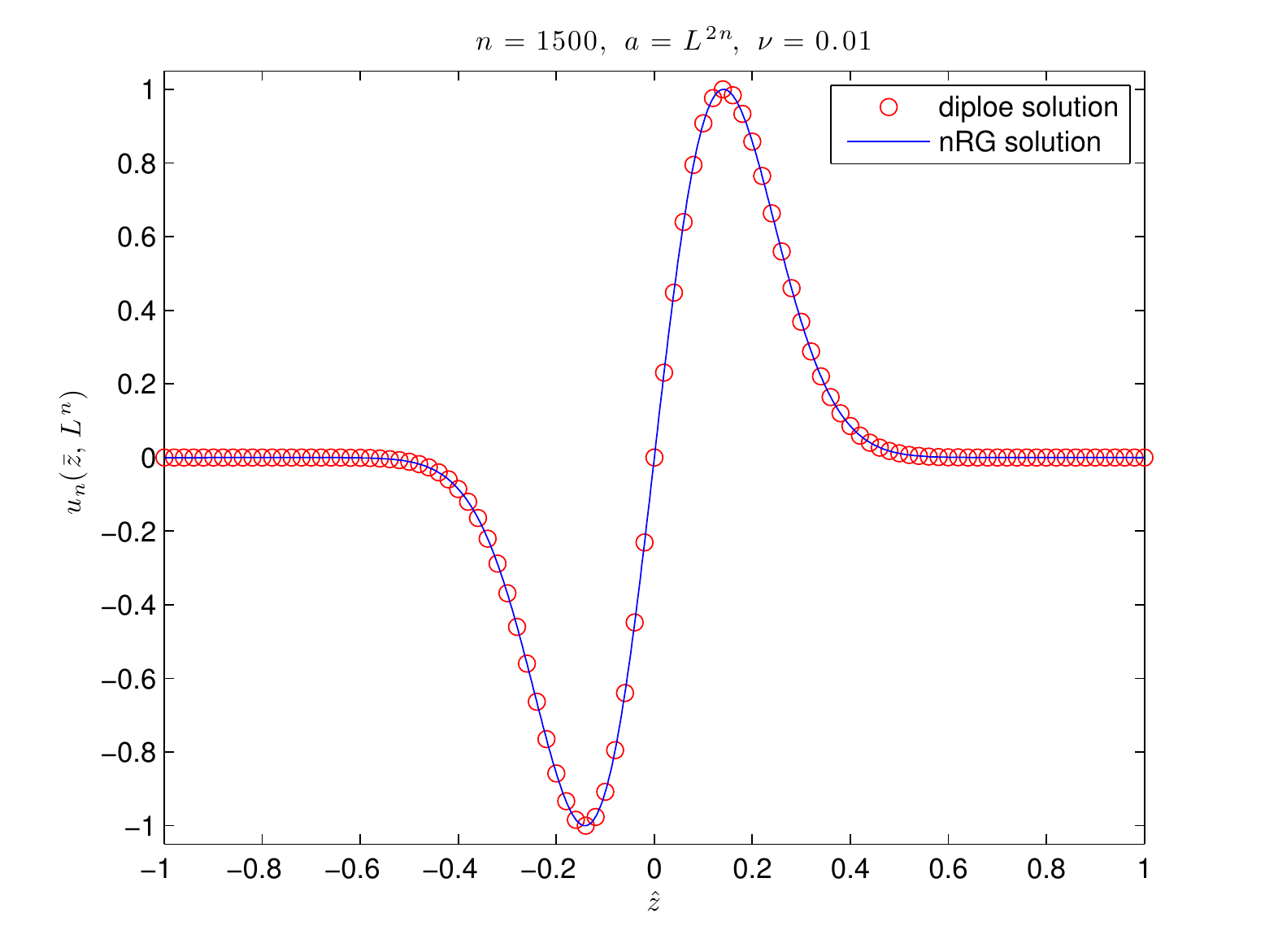}
\caption{Numerical experiment with zero initial mass ($M=0$). The viscosity is $\nu=0.01$. (a) Waveform exhibiting $N$-wave structure computed by the nRG algorithm at the $400^{th}$ iteration. (b) Comparison between the asymptotic dipole-like solution and the nRG calculation at the $1500^{th}$ iteration.}
\label{fig:zero_comp}
\end{figure} 
Eq. (\ref{eq:diploe}) is the dipole solution of the heat equation, which means  that the diffusion dominates the nonlinear term in the final decay, regardless the magnitude of the diffusivity constant. To study the final decay of solution by the nRG algorithm for this class of initial data, we let $\beta$ in the spacial scaling be fixed and equal to $1/2$. This ensures that the coefficient in front of the diffusion term remains unscaled at all time. We set the diffusivity constant $\nu=0.01$ and the parameter $L=1.2$. The total number of iteration is $n=1500$, and $\ell=1$ in the initial condition. Figure \ref{fig:zero_comp}(a) shows the snap shot of the nRG solution at the $400^{th}$ iteration. It clearly shows that the solution is of the $N$ wave structure at this stage. In Figure \ref{fig:zero_comp}(b), the final profile of the nRG calculation at $n=1500$ is compared with the dipole solution of heat equation. The fixed number $a$ in  Eq. (\ref{eq:diploe}) is chosen to be $L^{2n}$. Recall that the scaled spatial variable in the nRG calculation is  $\hat{z}=L^{-n(1/2)}x=x/\sqrt{t}$. Substituting $a=L^{2n}$ and $\hat{z}=L^{-n(1/2)}x=x/\sqrt{t}$ into Eq. (\ref{eq:diploe}) and knowing that $L^{2n}=t^2$ yields the dipole solution as a function of the similarity variable $\hat{z}$
\begin{equation}\label{eq:dipole_scaled}
u(\hat{z})=\hat{z}e^{-\hat{z}^{2}/(4\nu)}.
\end{equation}
The circles in Figure \ref{fig:zero_comp}(b) are the above dipole solution, with the amplitude being normalized to one, plotted against the similarity variable. The result of the nRG calculation at $n=1500$ (the solid line in  Figure \ref{fig:zero_comp}(b)) correctly predicts the final decay.
\begin{figure}[hbtp]
\centering
(a)\includegraphics[width=2.8in]{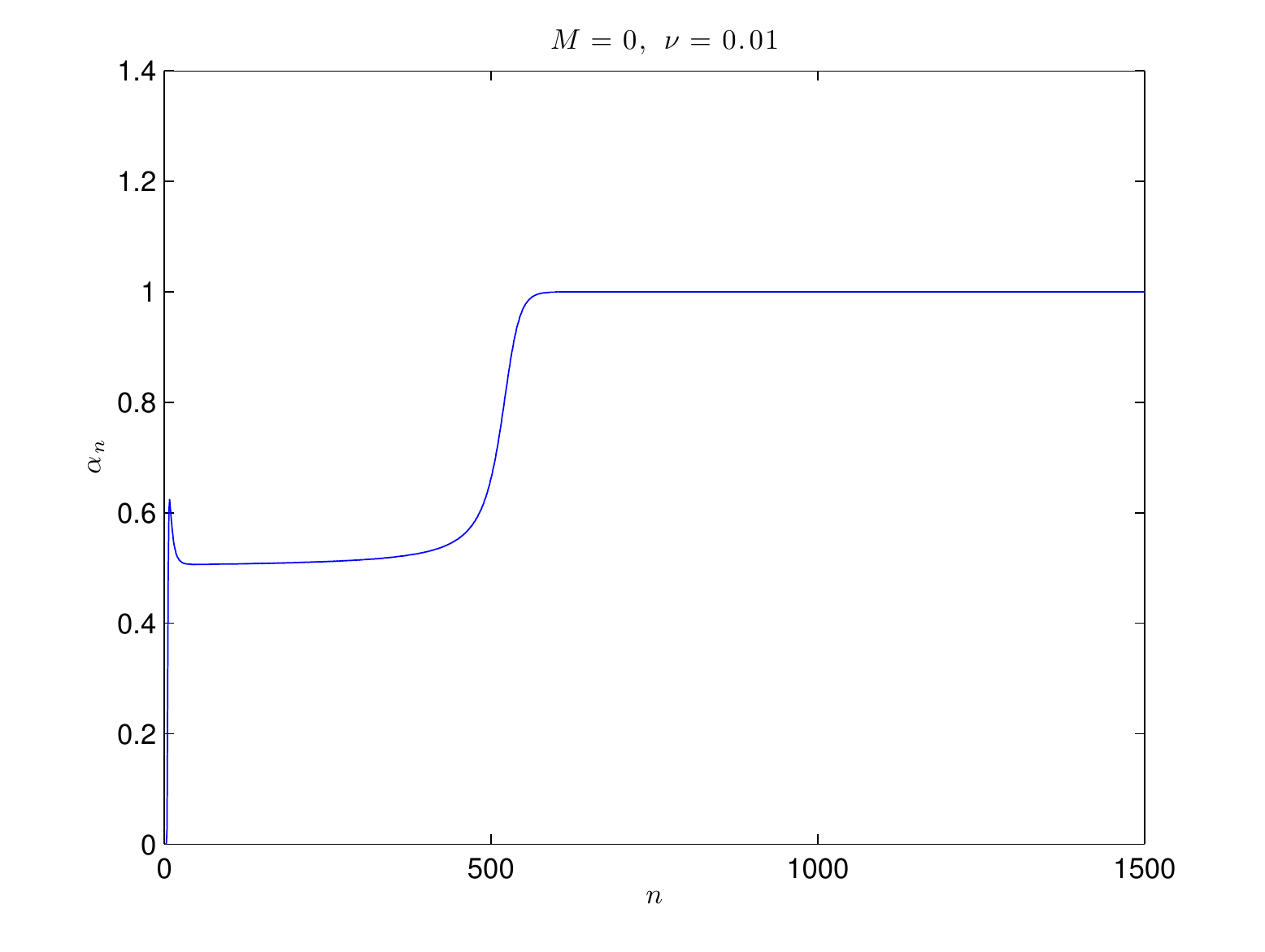} 
(b)\includegraphics[width=2.8in]{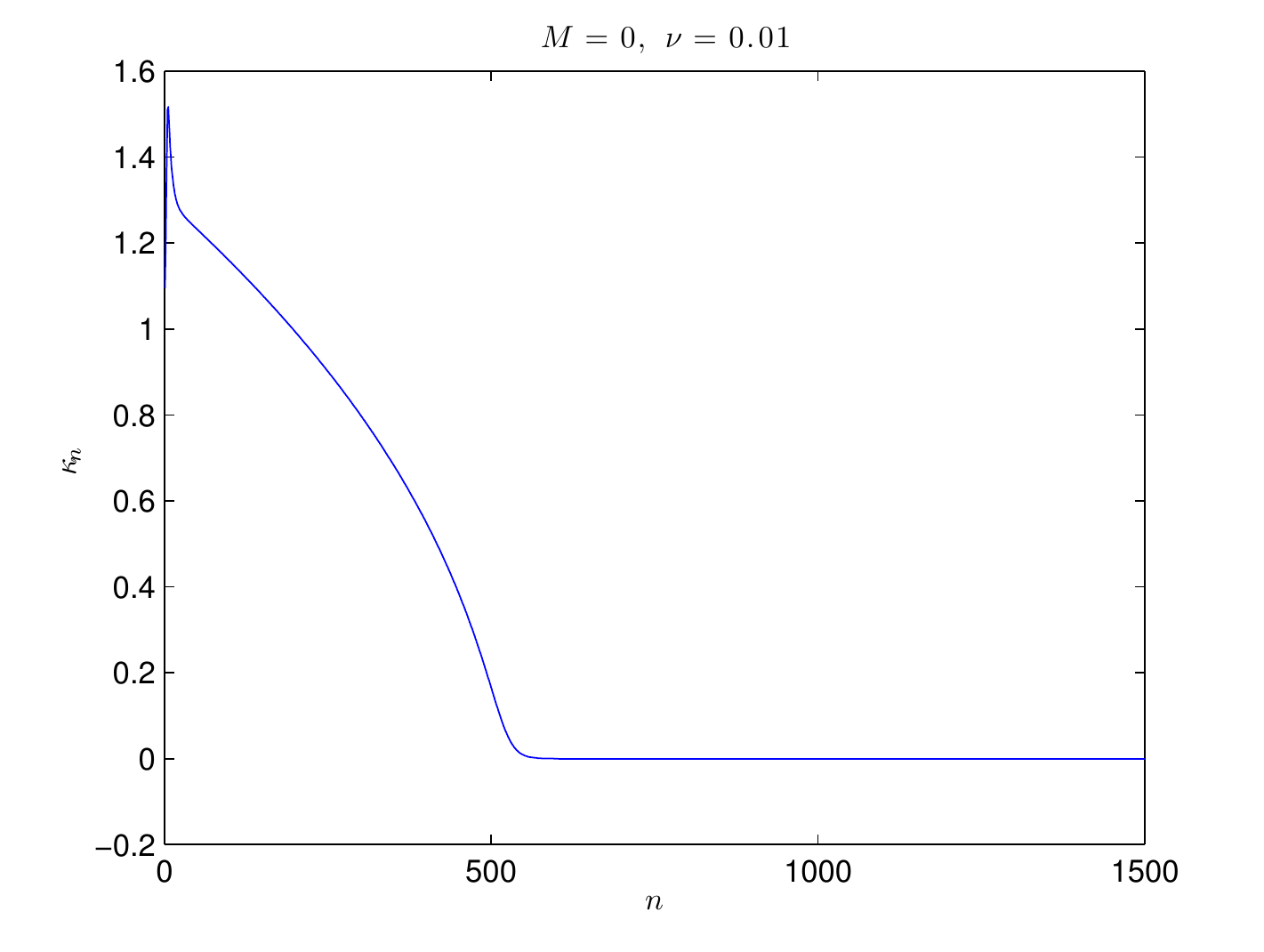}
\caption{Scaled coefficients in the previous calculation (Figure \ref{fig:zero_comp}). (a) $\alpha_n$ and (b) $\kappa_n$  versus $n$ (number of iteration). }
\label{fig:zero_mass_param}
\end{figure} 
The convergences of the scaling parameters $\alpha_n$ and $\kappa_n$ are shown in Figure \ref{fig:zero_mass_param}(a) and (b), respectively. They indicate that $\alpha_n\rightarrow 1$, which results in $\kappa_n\rightarrow 0$, where $\kappa_n$ is the coefficient in front of the advection term of the scaled PDE. Such a convergence of $\kappa_n$ suggests that the diffusion term dominates the final decay, regardless the magnitude of the diffusion constant, which is kept unscaled through out the calculation. These figures show the mechanism of the solution transiting from the $N$-wave structure to the dipole solution, as predicted in Whitham's analysis.
The choice of the periodic domain, $\Delta x$, and $\Delta t$ is the same as that used in Figure \ref{fig:positive_comp1}.  

We remark that to prevent the numerical artifacts from distorting the symmetry of the solution during the step of normalizing the amplitude, we apply the negative mirror image of the right hump for the left one at the end of each time evolution.

\section{Korteweg-de Vries equation and dispersive shock waves}\label{sec:KdV}

In this section, we illustrate that the nRG procedure is an efficient method for studying asymptotic behavior of dispersive shock waves (DSWs). DSWs appear when dispersion dominates dissipation for step-like data; they have been seen in plasmas, fluids, superfluids and optics \cite{AB13a}. In a sequence of papers, Ablowitz et al. analyze interactions and asymptotics of DSWs for the Korteweg-de Vries (KdV) equation \cite{ABH09, AB13a, AB13b}.  KdV equation is chosen for their study because it is the leading-order asymptotic equation for weakly dispersive and weakly non-linear systems.

Consider the dimensionless form of KdV equation
\begin{equation}\label{eq:KdV}
u_t+uu_x+\epsilon^2u_{xxx} = 0,
\end{equation}
with $u=u(x,t)$ going rapidly to the boundary conditions
\begin{equation}
\lim_{x\rightarrow -\infty} u(x, t)=0,\quad \lim_{x\rightarrow \infty} u(x, t)=-6c^2.
\end{equation}

Single stepwise initial data for the above problem evolves to a wedge-shape envelope combining three basic regions: exponential decay region on the right, the DSW region in the middle, and the region of oscillating tail on the left \cite{AB13a}. All three regions travel to the left with time at a speed $x=-12c^2 t$, while the DSW region is expanding and is of the order $O(t)$ \cite{AB13a, AB13b}. The amplitude of the DSWs saturates at  $6c^2$.

The KdV equation for the described problem can be posed as an initial-boundary-value problem (IBVP) on a truncated domain $-\ell\le x\le\ell$, where the initial and boundary data are prescribed as follows
\begin{equation}\label{eq:KdV_ibvp}
\begin{cases}
u_t+uu_x+\epsilon^2u_{xxx}&= 0,\quad x\in[-\ell,\ell],\,t>0,\\
u(x,0)&=3(1-\tanh((x-x_0)/w)-2),\\
u(-\ell,t)&=3(1-\tanh((-\ell-x_0)/w)-2),\\
u(\ell,t)&=3(1-\tanh((\ell-x_0)/w)-2),\\
u_x(\ell,t)&=0.
\end{cases}
\end{equation}
For example, if $w=1$ and $\ell-x_0=20$, we have $u(-\ell,t)\approx 0$ and $u(\ell,t)\approx-6$. The existence and uniqueness of solution of the above IBVP is discussed in \cite{BSZ03}. 

The numerical scheme we adopt for solving Eq. (\ref{eq:KdV_ibvp}) is a non-oscillatory explicit finite-difference method \cite{wwh08}.  The spatial and temporal discretization for the algorithm is 
\beq\label{eq:kdv_FD_algorithm}
\begin{split}
\frac{1}{2}\left(\frac{u^{n+1}_{j-1} - u^{n}_{j-1}}{\Delta t} +\frac{u^{n}_{j+1} - u^{n-1}_{j+1}}{\Delta t} \right)  = & - \left(\frac{u^n_{j+1}+u^n_{j}+u^n_{j-1}}{3}\right)\frac{u^n_{j+1}-u^n_{j-1}}{\Delta x}\\
&-\frac{\epsilon^2}{2\Delta x^3}\left(u^n_{j+2}-2u^n_{j+1}+2u^n_{j-1}-u^n_{j-2}\right).
\end{split}
\eeq
Applying the von Neumann analysis for the above scheme yields the stable condition for the scheme
\beq\label{eq:stable-cond}
\frac{\Delta t}{\Delta x} < \frac{2}{\underset{x, t}{\max}|u| + 4 \epsilon^2 / \Delta x}.
\eeq
Using the same scaling variables as Eq. (\ref{eq:scale_tx}) and Eq. (\ref{eq:scale_u}), the scaled KdV equation is 
\begin{equation}
u_t + L^{-\alpha-\beta+1}uu_x +\epsilon^2 L^{-3\beta +1}u_{xxx} = 0.
\end{equation}
Here we drop the subscript $L$ for $u$ and $\,\tilde{}$ for $t$ and $x$.
Since the amplitude of DSWs saturates at $6c^2$, it is not necessary to scale the amplitude and hence we set $\alpha=0$ for our nRG calculations. Also, similar to the Burgers equation, we choose to retain the coefficient of the dispersion term unscaled. This results in $\beta=1/3$ at all time for our nRG calculations, which suggests that the DSW region expands in the order $O(t^{1/3})$ for the RG calculations. We choose $c^2=1$. 

\begin{figure}[bhtp]
\centering
(a)\includegraphics[width=2.8in]{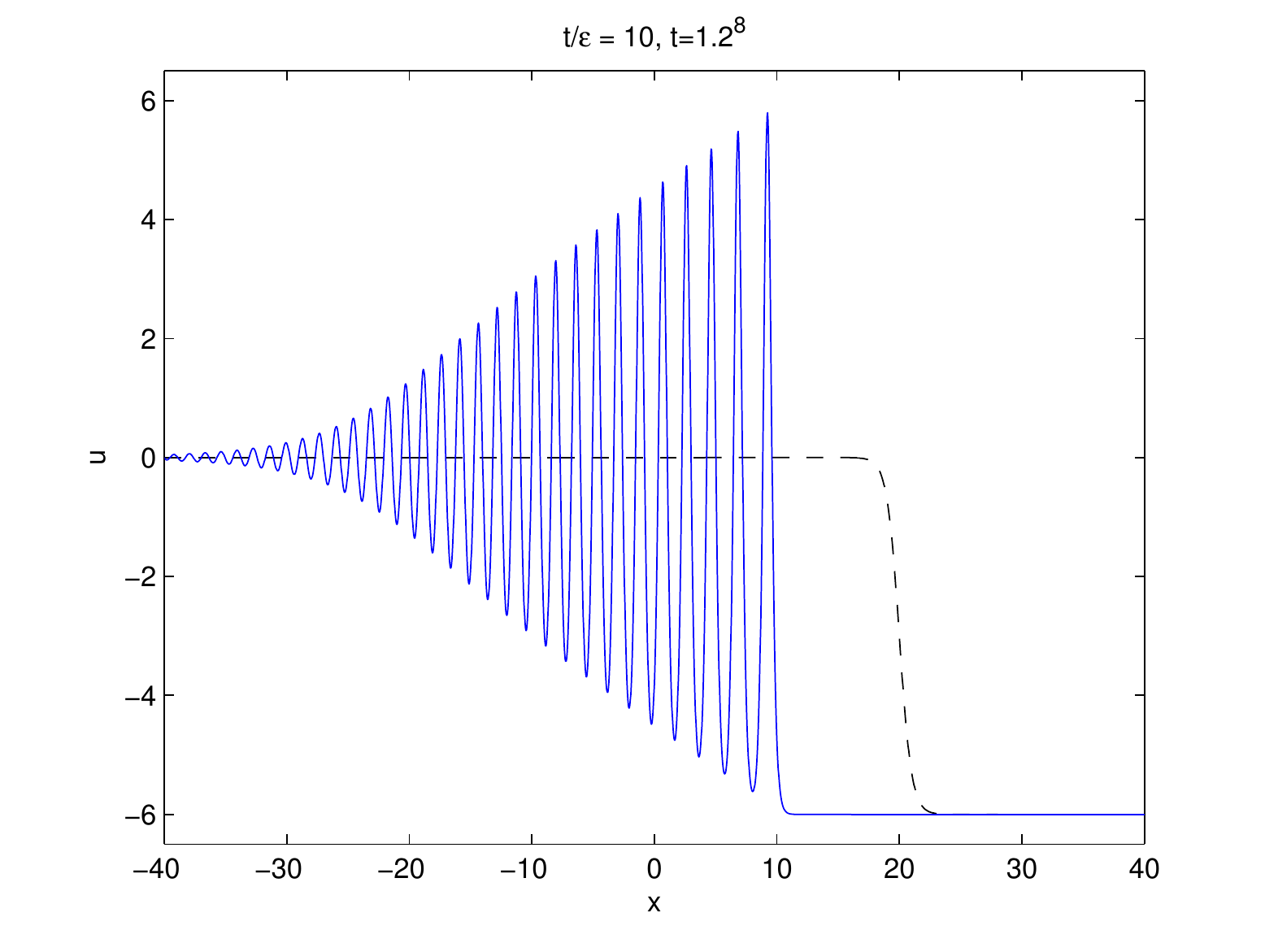} 
(b)\includegraphics[width=2.8in]{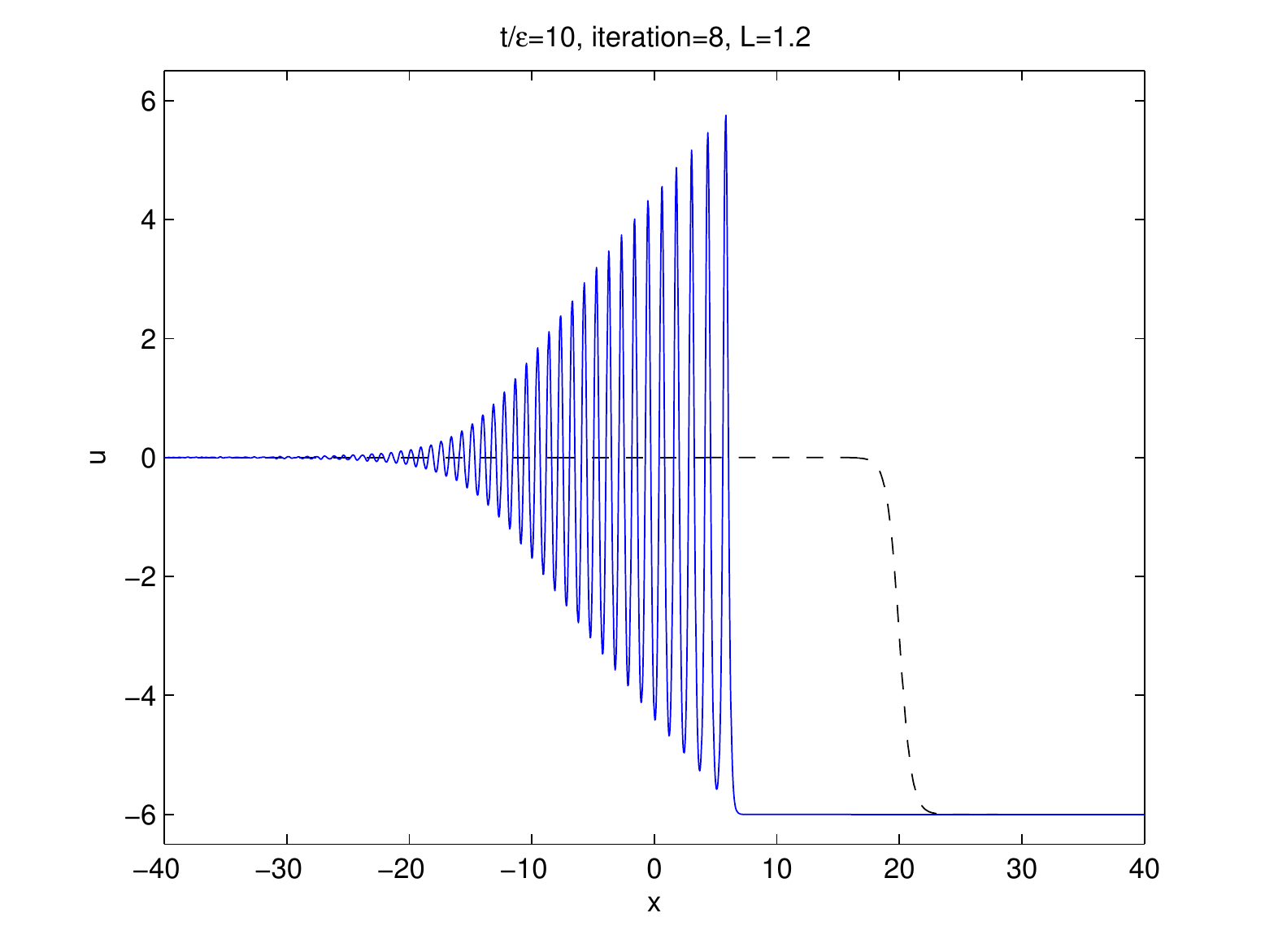} 
\caption{(a) Direct numerical simulation. $t/ \epsilon=10$, where $t=1.2^8$. The spacial and temporal step sizes are $\Delta x=40/2000$ and $\Delta t=1.2^{8}/80000$, respectively. $x_0=20$ and $w=1$. (b) Numerical renormalization group calculation. $\beta=1/3$, $\alpha=0$, $L=1.2$, and $t/ \epsilon =10$. Eight iterations is performed ($n=8$, i.e. $t=1.2^{8}$). For both (a) and (b), dashed line is the initial condition.}
\label{fig:KdV_step}
\end{figure}

\begin{figure}[bhtp]
\centering
(a)\includegraphics[width=2.8in]{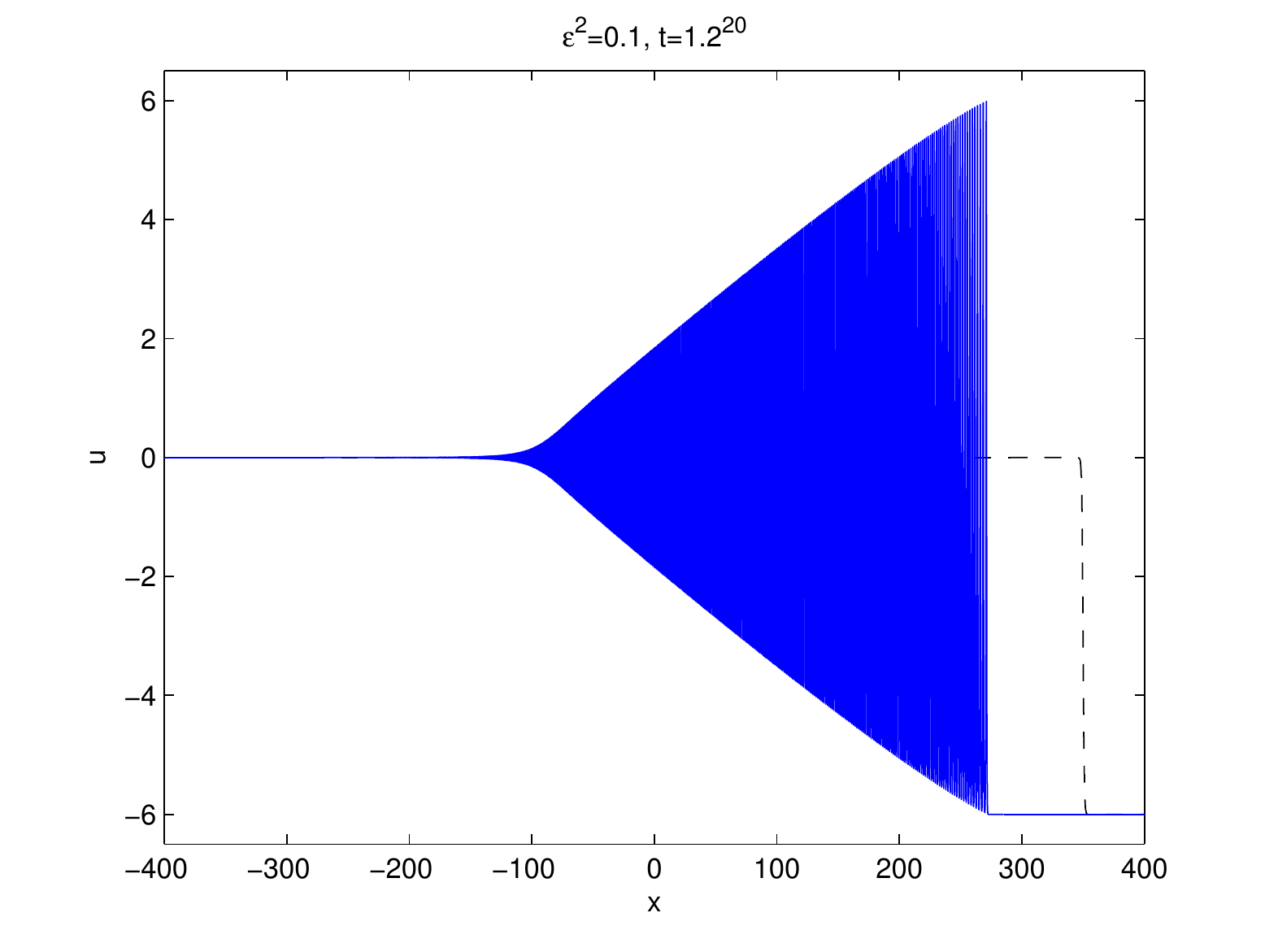} 
(b)\includegraphics[width=2.8in]{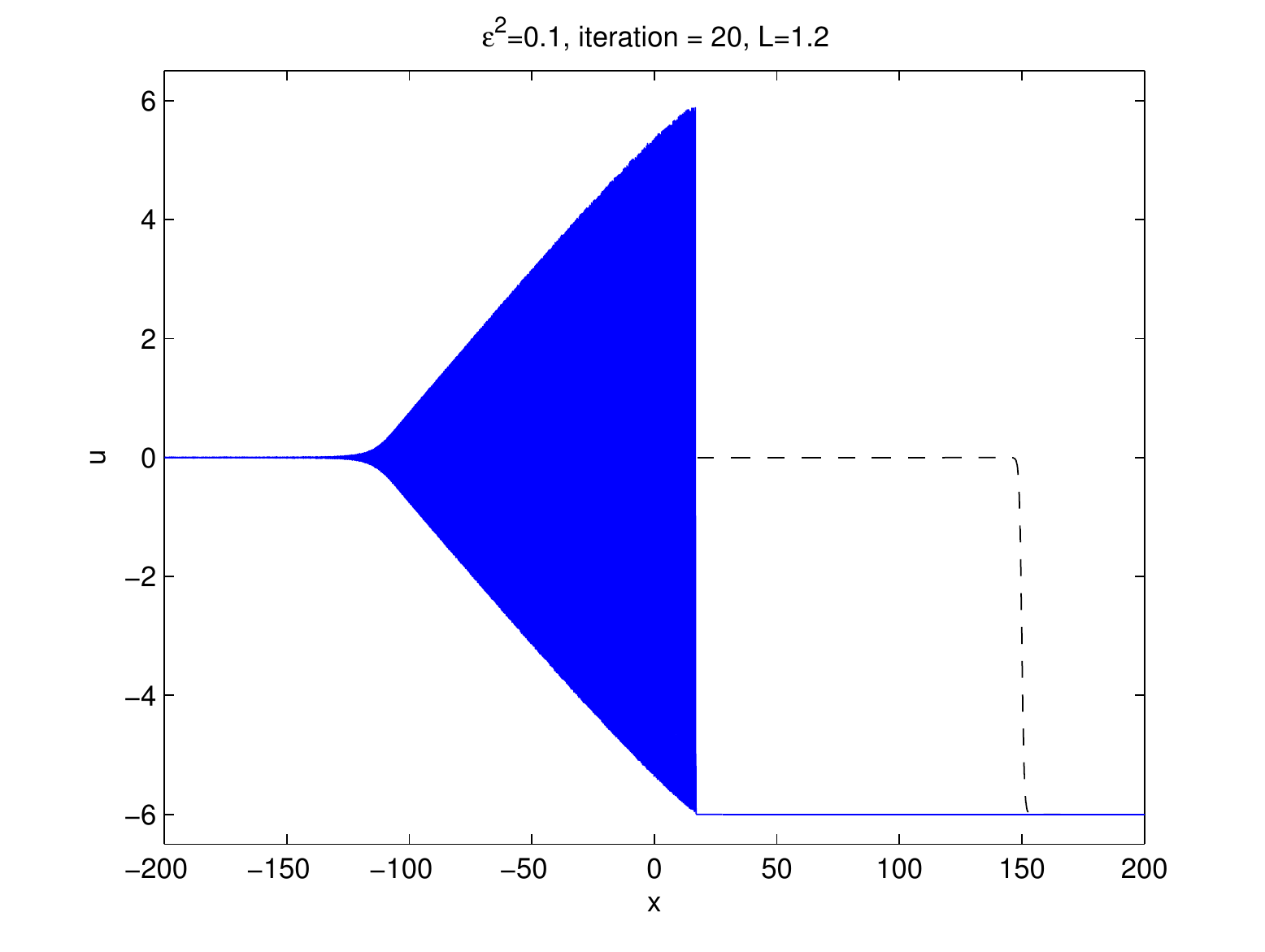} 
\caption{(a) Direct numerical simulation. $x_0=350$, $w=1$, $\Delta x=1/80$, $\Delta t=(1.2)^{20}/32000000\approx 1.19805\times 10^{-6}$. $\epsilon^2=0.1$. (b) Numerical renormalization group calculation. $\beta=1/3$, $\alpha=0$, $L=1.2$, and $\epsilon^2=0.1$. Twenty iterations is performed ($n=20$, i.e. $t=1.2^{20}$), $x_0=150$, $w=1$, $\Delta x=1/80$, and $\Delta t=1\times 10^{-6}$. For both (a) and (b), the dashed line is the initial condition.}
\label{fig:KdV_step2}
\end{figure}

With this set of parameters, Figure \ref{fig:KdV_step} and Figure \ref{fig:KdV_step2} are the comparisons between direct numerical simulations and nRG calculations. The initial condition is a single-step tangent profile. The graphs show that the nRG procedure accurately capture the DSWs in a confined domain and the results are consistent with that in the literature \cite{AB13a, AB13b}. Moreover a simple calculation below will illustrate that the nRG procedure is more efficient than direct numerical calculation. 

Suppose that the final time for a direct simulation is $t=L^n$ and the number of solitons in the region of DSWs at the final time is $N_{s}$. Assume that the spacial grid-size required to resolve these solitons is $\Delta x$. For this $\Delta x$, the required temporal step-size is $\Delta t =O((\Delta x)^2)$, by the stable condition (\ref{eq:stable-cond}). Therefore the total number of time steps required for the simulation is 
$L^{n}/\Delta t \sim L^{n}/(\Delta x)^2$. Now for the nRG procedure, $t=L^n$ means the number of iterations is $n$. Since  the dispersion coefficient is kept the same, the number of solitons after $n$ iterations is also $N_{s}$ for the nRG procedure. However, because of the spacial scaling, the spacial grid-size required for resolving the solitons is now $\Delta\tilde{x}=L^{-n\beta}\Delta x$. Hence the temporal step-size for stability requirement is $\Delta\tilde{t} = O((\Delta\tilde{x})^2)$. The number of time steps for the nRG procedure to the final time is $n(L-1)/\Delta\tilde{t} \sim n(L-1)L^{2n\beta}/(\Delta x)^2$. Since $\beta=1/3$, the numerator of the above equation is $n(L-1)L^{2n/3}$ and this is less than $L^n$ for large $n$. 

Our numerical experiments also confirm that direct numerical simulation is much more time consuming than the nRG calculation for long-time simulations.  We further remark that if an implicit algorithm is used, for which $\Delta t$ could be chosen at the same order of $\Delta x$ \cite{SK09}, then the nRG procedure is even more preferable than direct numerical simulation.

\section{A modified diffusion-absorption model}\label{sec:diffusion_absorption}

In this section we consider a one-dimensional modified diffusion-absorption model 
\begin{equation}\label{eq:diff-absorp}
u_t = D(u_t) \partial_{xx}u^{m+1} -\lambda u^p,
\end{equation}
where $m\ge 0$ and $D(u_t)$ is the transport coefficient. Typically, in the literature the following cases are considered:
\begin{enumerate}
\item $D(u_t)=1$ and $m\ge1$. In this case, if $\lambda >0$ and $p\ge m+1$, this is a model of nonlinear heat equation with absorption \cite{bib:HV87}, whereas if $\lambda > 0$ and $0<p<1$, this is slow diffusion combined with strong absorption \cite{bib:GSV99}. 
\item $D(u_t)$ is a Heaviside function:
\begin{equation}\label{eq:Heaviside}
D(u_t) = \begin{cases}  1+\epsilon & \quad \text{for}\,\, u_t  < 0 \\
1 & \quad \text{for}\,\, u_t  >  0.
\end{cases} 
\end{equation}
For this case, if $\lambda =0$ and $m=0$, this is the so-called Barenblatt's equation. \cite{bib:CW96, GMOL90, Barenblatt87}, whereas if $\lambda =0$ and $m=1$, it is called the modified porous-medium equation \cite{CGO91}.
\end{enumerate}
In this section, our investigation focuses on the comparison of the longtime solution behavior between the discontinuous $D(u_t)$ and the following non-constant continuous alternatives:
\begin{enumerate}
\item $D(u_t)$ is a smooth transitioned profile connecting $1+\epsilon$ and 1, or
\begin{equation}\label{eq:D(u_t)}
D(u_t) = 1+ \epsilon\left(\frac{1}{2}\left(1+\tanh(-u_t/\sigma)\right)\right),
\end{equation}
where $D(u_t)\rightarrow 1+\epsilon$ as $u_t\rightarrow -\infty$  and $D(u_t)\rightarrow 1$ as $u_t\rightarrow\infty$.  $\sigma > 0$ is a parameter to adjust the width of the smooth transition of the hyperbolic tangent profile. 
\item $D(u_t)$ is a piecewise continue function that connects the two constant states $1+\epsilon$ and $1$ by a straight line, or
\begin{equation}\label{eq:D(u_t)2}
D(u_t) = \begin{cases}
1+\epsilon &\quad \text{for}\,\,u_t < -\delta \\
1-\displaystyle\frac{\epsilon}{2\delta}\left(u_t-\delta\right) & \quad \text{for}\,\,-\delta<u_t<\delta\\
1& \quad \text{for}\,\,u_t > \delta,
\end{cases}
\end{equation}
where $\delta$ plays the similar role as  $\sigma$ in Eq. (\ref{eq:D(u_t)}).
\end{enumerate}


\subsection{Validation for the RG algorithm}\label{sec:barenblatt}
The self-similar solution of the Barenblatt's equation studied in the literature \cite{bib:CW96} is an ideal example to validate our nRG algorithm. Consider the Barenblatt's equation
\begin{equation}\label{eq:Baren}
\begin{split}
&u_t = D(u_t) u_{xx},\\ 
&u(t=0,x) = u_0(x),
\end{split}
\end{equation}
where $D(u_t)$ is a Heaviside function defined in Eq. (\ref{eq:Heaviside}). The asymptotic similarity form of the Barenblatt's equation is 
\begin{equation}\label{eq:time_decay}
u(x,t)\simeq \frac{A}{t^\alpha}\phi\Big(\frac{x}{2\sqrt{\kappa_{+} t}}\Big),\quad \kappa_{+} = \,\,\text{diffusivity in the regime where}\,\, \frac{\partial}{\partial t} > 0
\end{equation}
\cite{bib:CW96}, where $A$ is some pre-factor. For our numerical validation $\kappa_{+} =1$ and $\kappa_{-} =1+\epsilon$, as shown in Eq. (\ref{eq:Heaviside}). The parameter $\alpha$ is a function of the diffusivity ratio $\kappa_{-} / \kappa_{+}$. Cole and Wagner \cite{bib:CW96}  showed that the asymptotic expansions of $\alpha$ in $\epsilon$ is 
\begin{equation}
\alpha(\epsilon) = \alpha_0 + \epsilon\alpha_1+\epsilon^2\alpha_2+\cdots,
\end{equation}
where 
\begin{equation}\label{eq:perturbation}
\alpha_0 = \frac{1}{2},\quad \alpha_1= \frac{1}{\sqrt{2\pi e}},\quad \text{and}\,\,\,\,\,\alpha_2 = -0.06354624.
\end{equation}
We remark that since the asymptotic similarity solution of the linear heat equation has the time decay rate $\alpha=\alpha_0=1/2$, the $\epsilon$ and $\epsilon^2$ terms in Eq. (\ref{eq:perturbation}) are sometimes called {\it the anomalous dimension} of the decay.

\begin{figure}[h]
\centering
\includegraphics[width=4.8in]{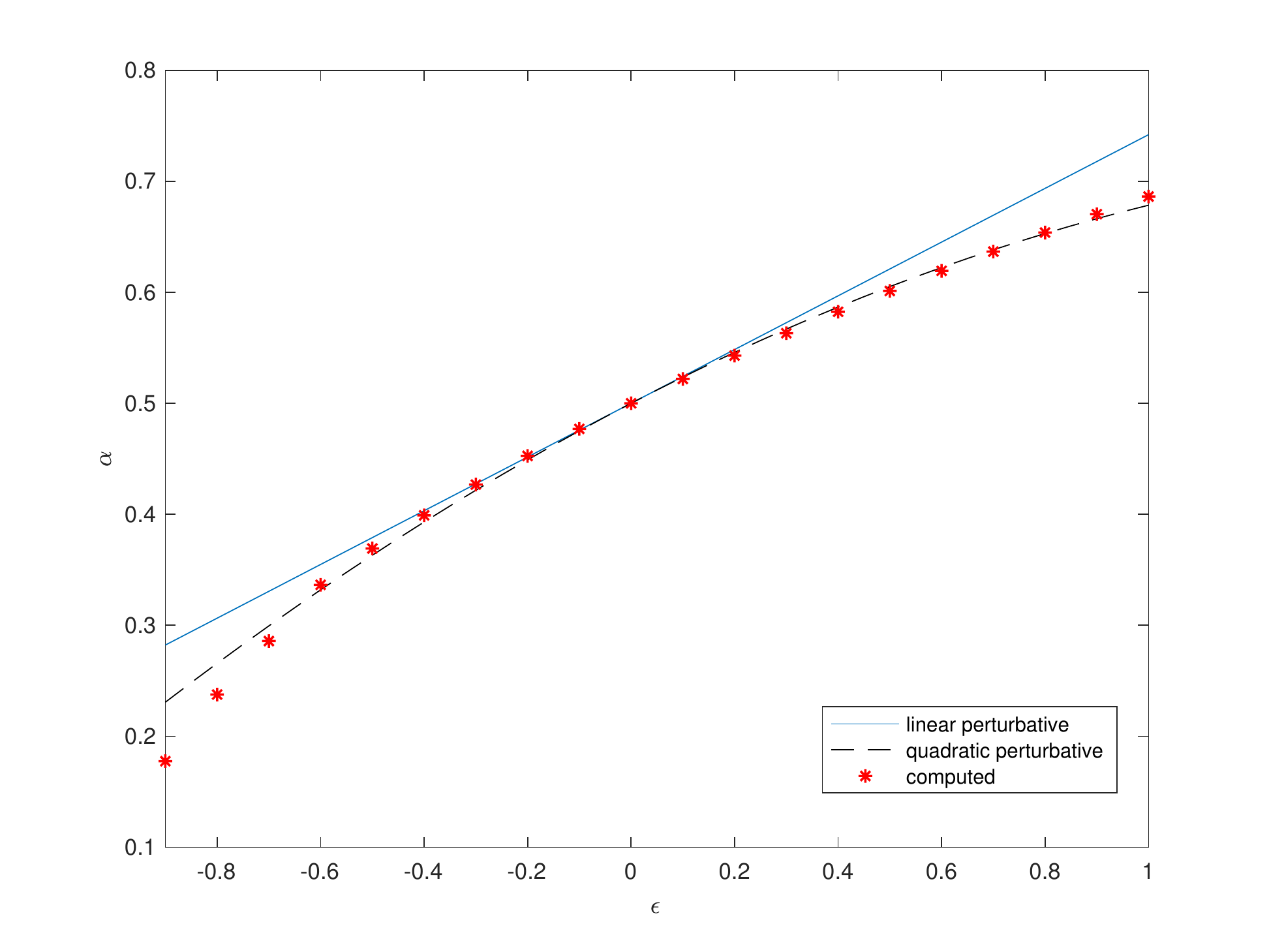} 
\caption{The time decay exponent $\alpha$ in Eq. (\ref{eq:time_decay}) captured  by the nRG algorithm versus the perturbative values predicted by the linear and quadratic formula, Eq. (\ref{eq:perturbation}), given in the literature \cite{bib:CW96}.}
\label{fig:validation}
\end{figure}

Eq. (\ref{eq:Baren}), the Barenblatt's equation, is essentially a nonlinear equation, since the diffusivity is a function of $u_t$.  Suppose that the time and space variables and the amplitude of $u$ are scaled the same way as that in Eqs. (\ref{eq:scale_tx}) and (\ref{eq:scale_u}), the scaled Barenblatt's equation for the $n^{th}$ RG iteration becomes
\begin{equation}\label{eq:scaled_Baren}
\begin{split}
&(u_n)_{t}=D\left(L^{-n(\bar{\alpha}_{n}+1)}(u_{n})_t\right) L^{n(-2\bar{\beta}_{n}+1)} (u_n)_{xx},\,\,t >1,\\
\text{I. C.\,\,:}\quad&u_n(x,1) = f_n(x),
\end{split}
\end{equation}
where $\bar{\alpha}_{n}$, $\bar{\beta}_{n}$ and  $f_n(x)$ have been defined in Section \ref{sec:seq}. We solve the above initial-value-problem by choosing an initial condition
\begin{equation}\label{eq:IC_experiment}
u_0(x, 1) =
\begin{cases}
 \cos x,\quad &-\frac{\pi}{2} \le x\le \frac{\pi}{2},\\
0\quad &x > \frac{\pi}{2}\,\,\text{or}\,\, x < -\frac{\pi}{2}.
\end{cases}
\end{equation} 
and discretizing Eq.(\ref{eq:scaled_Baren}) with the $2^{nd}$-order  Crank-Nicolson scheme. The diffusivity at the $(k+1)^{th}$ time step in the $n^{th}$ RG iteration is linearized by $D(L^{-n(\bar{\alpha}_{n}+1)}(3u_n^{k}-4u_n^{k-1}+u_n^{k-2})/(2\Delta t))$. In our nRG calculation, the spacial scaling parameter $\beta=1/2$ is fixed, so that the magnitude of the diffusivity is unscaled of all time. The time integration span in each nRG iteration is from $L=1$ to $L=2$, while the total iteration number is 100. A periodic domain $x\in [-10, 10]$ with the temporal step $\Delta t =  0.05$ and the spacial step $\Delta x =  \frac{20}{160}= 0.125$ are used for all numerical computations with different $\epsilon$-values. Cubic interpolation scheme is employed for the grid interpolation approach discussed in Section \ref{sec:implement}. Figure \ref{fig:validation} is the comparison between the time decay exponent $\alpha$ in Eq. (\ref{eq:time_decay}) captured  by the nRG algorithm and the linear and quadratic perturbative values predicted by Eq. (\ref{eq:perturbation}) given in the literature \cite{bib:CW96}. The actual values of  $\alpha(\epsilon)$ for the comparison for 20 $\epsilon$-values, $\epsilon = -0.9, -0.8,\cdots, 0.9, 1.0$, are listed in Table \ref{tab:validation} in Appendix \ref{app:data1}. From both the Figure and the Table, we see that when $\epsilon$ is small, the $\alpha$-values given by the nRG algorithm and the perturbative formula are almost coincided, since the formula are derived by assuming small $\epsilon$.

\subsection{Non-constant continuous $D(u_t)$}

\begin{figure}[h]
\centering
(a)\includegraphics[width=2.8in]{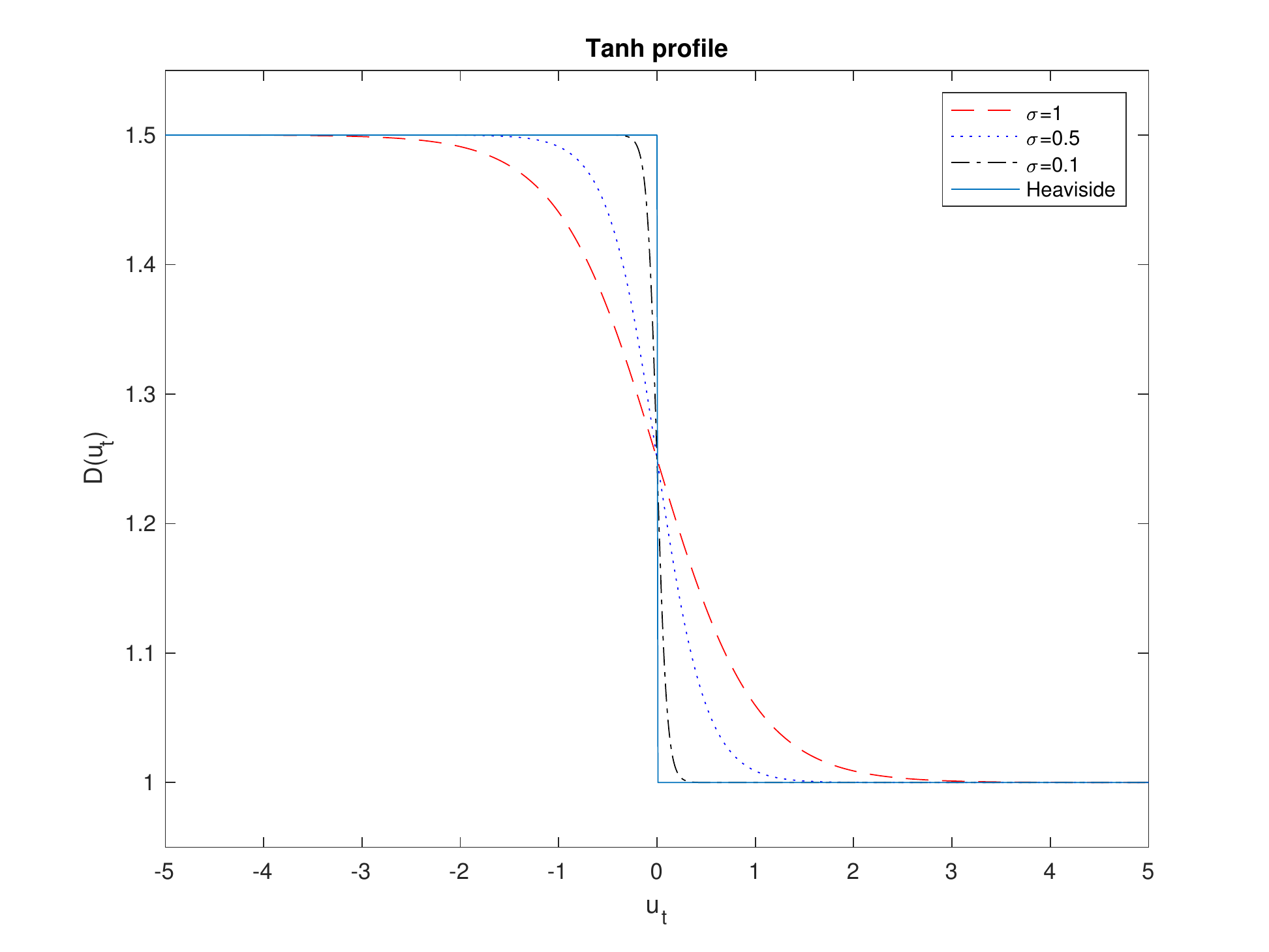} 
(b)\includegraphics[width=2.8in]{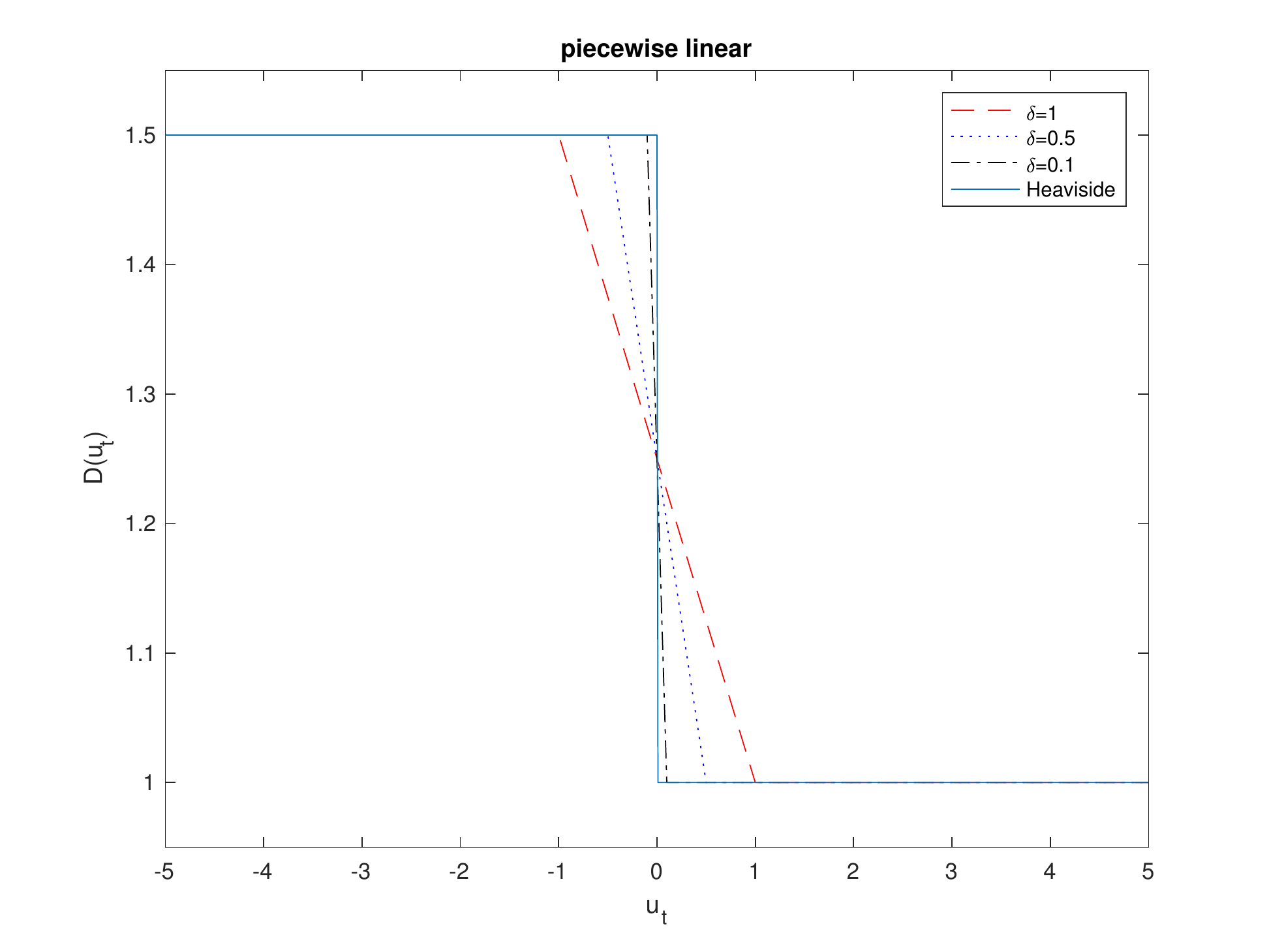} 
\caption{(a)  Smooth transitioned $D(u_t)$ in Eq. (\ref{eq:D(u_t)}), where $\epsilon=0.5$ and $\sigma=1, 0.5$, and 0.1. (b) Piecewise linear $D(u_t)$ in Eq. (\ref{eq:D(u_t)2}), where $\epsilon=0.5$ and $\delta=1, 0.5$, and 0.1.}
\label{fig:D(u_t)}
\end{figure}

Suppose now the Heaviside diffusivity function (\ref{eq:Heaviside}) is replaced by two types of continuous functions: (i) a smooth transitioned function described in Eq. (\ref{eq:D(u_t)}) with $\epsilon=0.5$ and $\sigma=1, 0.5$ and 0.1, as shown in Figure \ref{fig:D(u_t)}(a), and (ii) a piecewise linear function described in Eq. (\ref{eq:D(u_t)2}) with $\epsilon=0.5$ and $\delta=1, 0.5$ and 0.1, as shown in Figure \ref{fig:D(u_t)}(b). We will demonstrate numerically the change of the behavior for the decaying parameter $\alpha$ with such a replacement. Before our numerical experiments, it is worth pointing out that the Heaviside diffusivity function under the time-scaling is
\begin{equation}
D\left(L^{-n(\bar{\alpha}_{n}+1)}((u_{n}))_t\right)=D((u_{n})_t) = \begin{cases}  1+\epsilon & \quad \text{for}\,\, (u_{n})_t  < 0 \\
1 & \quad \text{for}\,\, (u_{n})_t  >  0,
\end{cases} 
\end{equation}
for $n\rightarrow\infty$. i.e. the value of the diffusivity depends only on the sign of $(u_{n})_t$.  The diffusivity functions described in (i) and (ii), however, behave differently under the time-scaling:
\begin{equation}\label{eq:D(Lu_t)}
D\left(L^{-n(\bar{\alpha}_{n}+1)}((u_{n}))_t\right)\longrightarrow D(0)=1+\frac{\epsilon}{2}, \quad \text{as}\,\,\, n\longrightarrow\infty, \quad \text{if}\,\,\, \bar{\alpha}_n > -1\,\,\,\text{and}\,\,\, |(u_{n})_t| < \infty.
\end{equation}
Equivalently to say, the above back-of-the-envelope calculation suggests that if the diffusivity approaches to a constant, the Barenblatt-like's equation approaches to a linear heat equation, and consequently the time decay parameter $\alpha$ approaches to 1/2, i.e. $\bar{\alpha}_n \rightarrow 1/2$.

\begin{figure}[h]
\centering
(a)\includegraphics[width=2.8in]{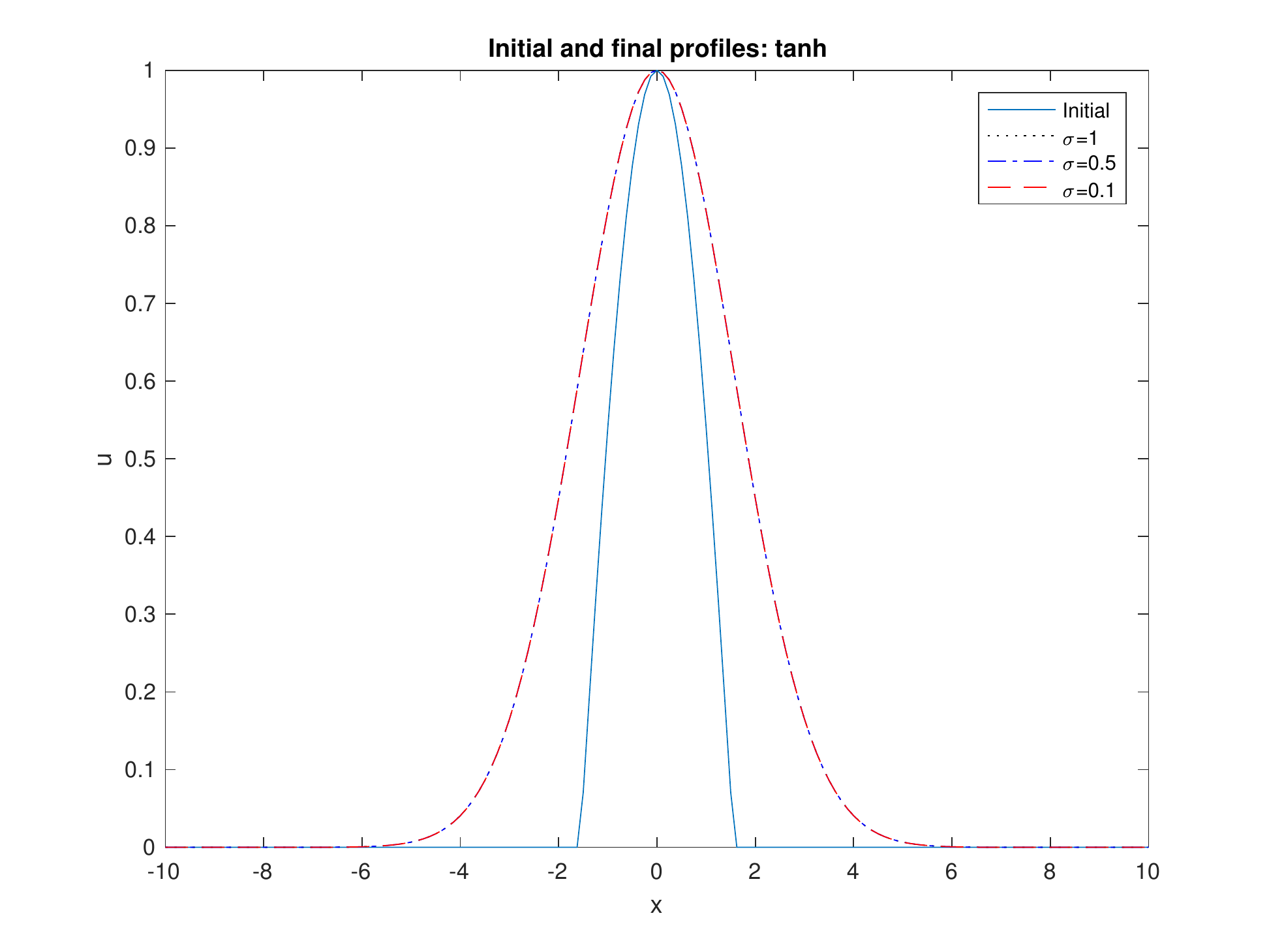} 
(b)\includegraphics[width=2.8in]{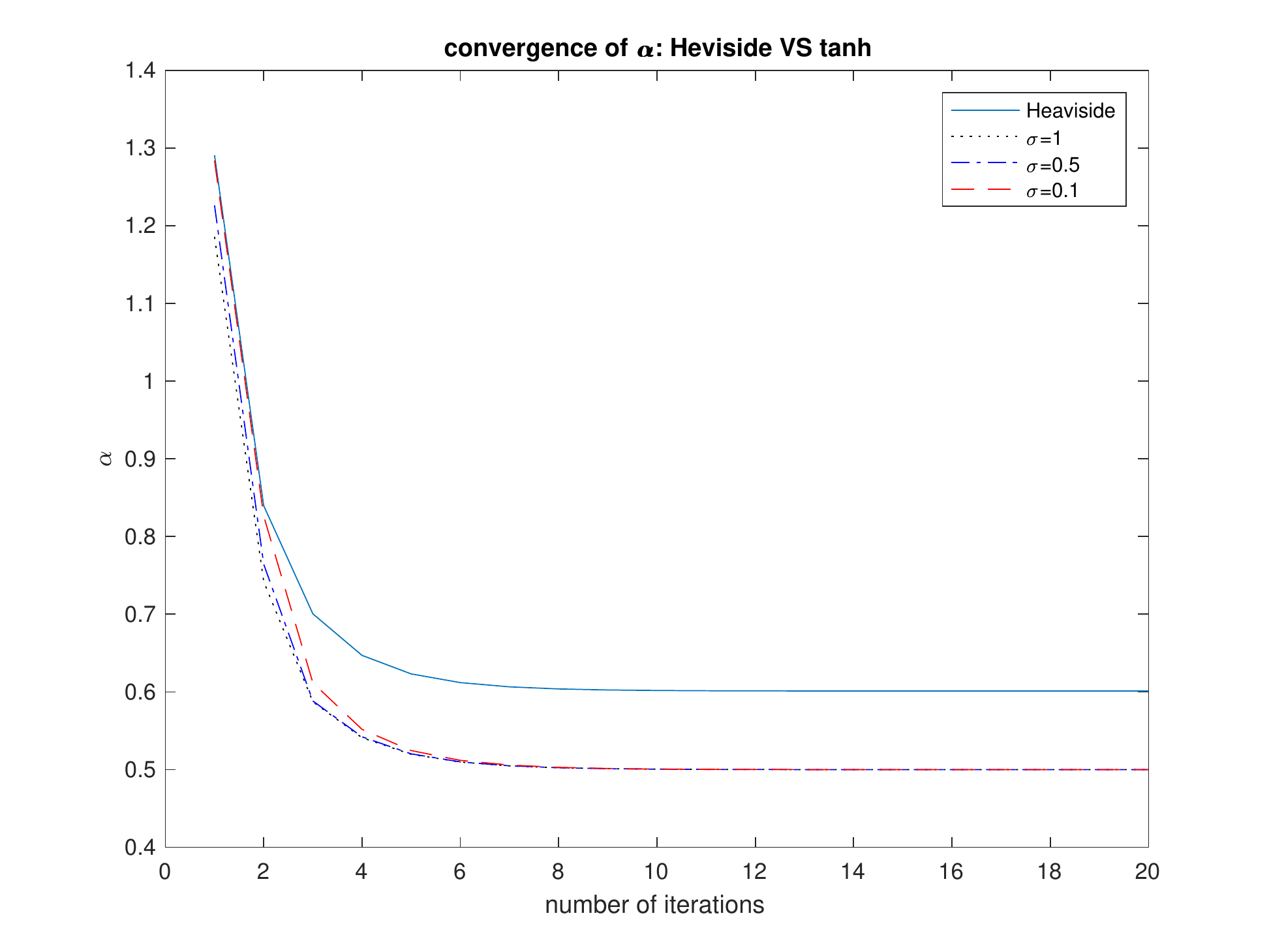} 
\caption{ (a) The initial profile and the computed asymptotic similarity forms for Eq. (\ref{eq:scaled_Baren}). Three smooth transitioned $D(u_t)$ shown in Figure \ref{fig:D(u_t)}(a) are used for the calculations. The final asymptotic similarity forms computed after 200 nRG iterations are visually indistinguishable for the three different bandwidths of hyperbolic tangent profiles. (b) The time decay exponent $\alpha$ captured by the nRG algorithm for the calculations in (a), compared with that of the Barenblatt's equation, for which $D(u_t)$ is a Heaviside function. $\alpha$ approaches to 1/2 for the continuous smooth $D(u_t)$ after only 20 iterations.}
\label{fig:tanh}
\end{figure}

\begin{figure}[h]
\centering
(a)\includegraphics[width=2.8in]{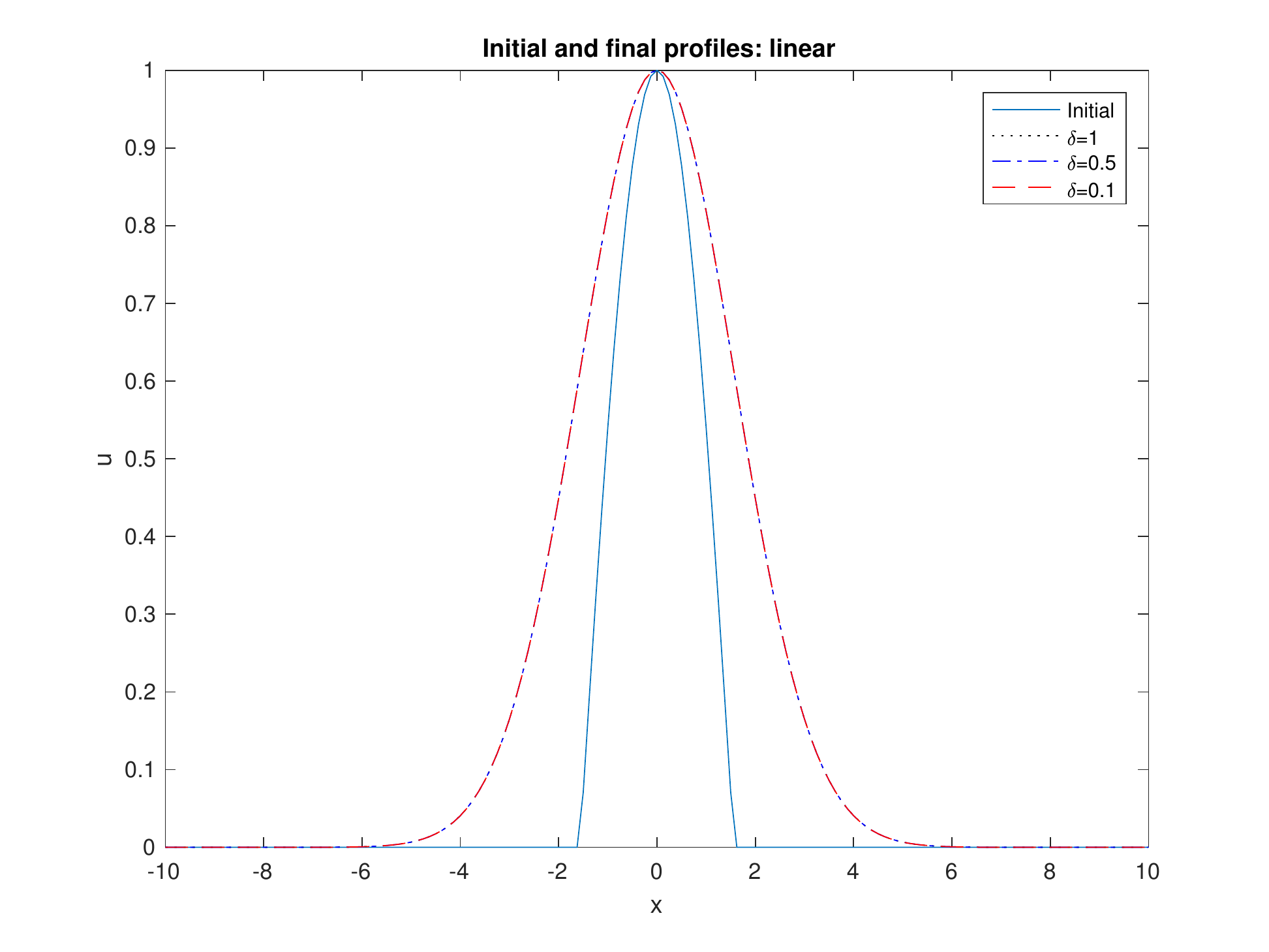} 
(b)\includegraphics[width=2.8in]{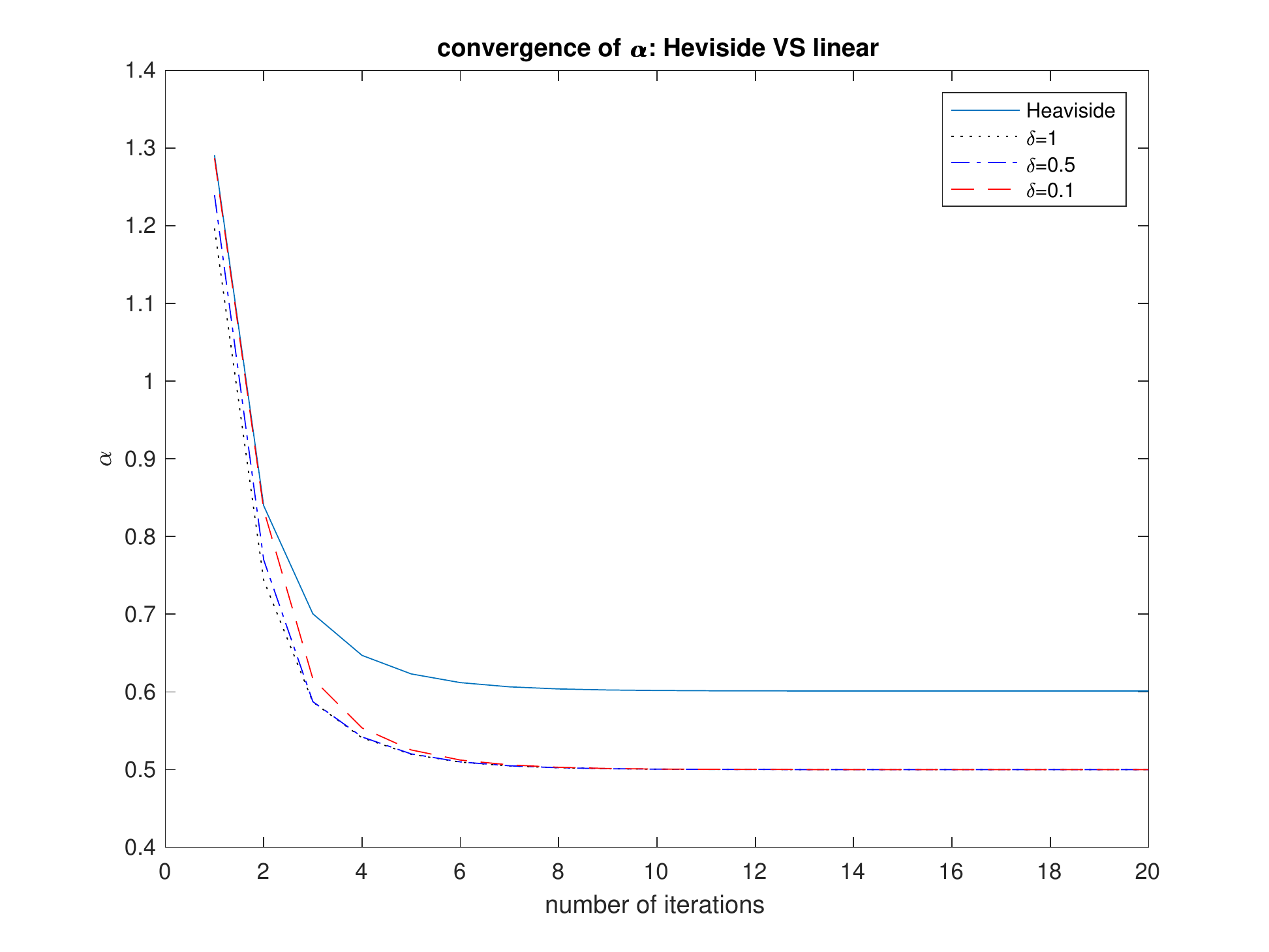} 
\caption{The same calculation as Figure \ref{fig:tanh}, except the diffusivity is a continuous function that connects the two constant states with a straight line.}
\label{fig:linear}
\end{figure}

\begin{figure}[h]
\centering
(a)\includegraphics[width=2.8in]{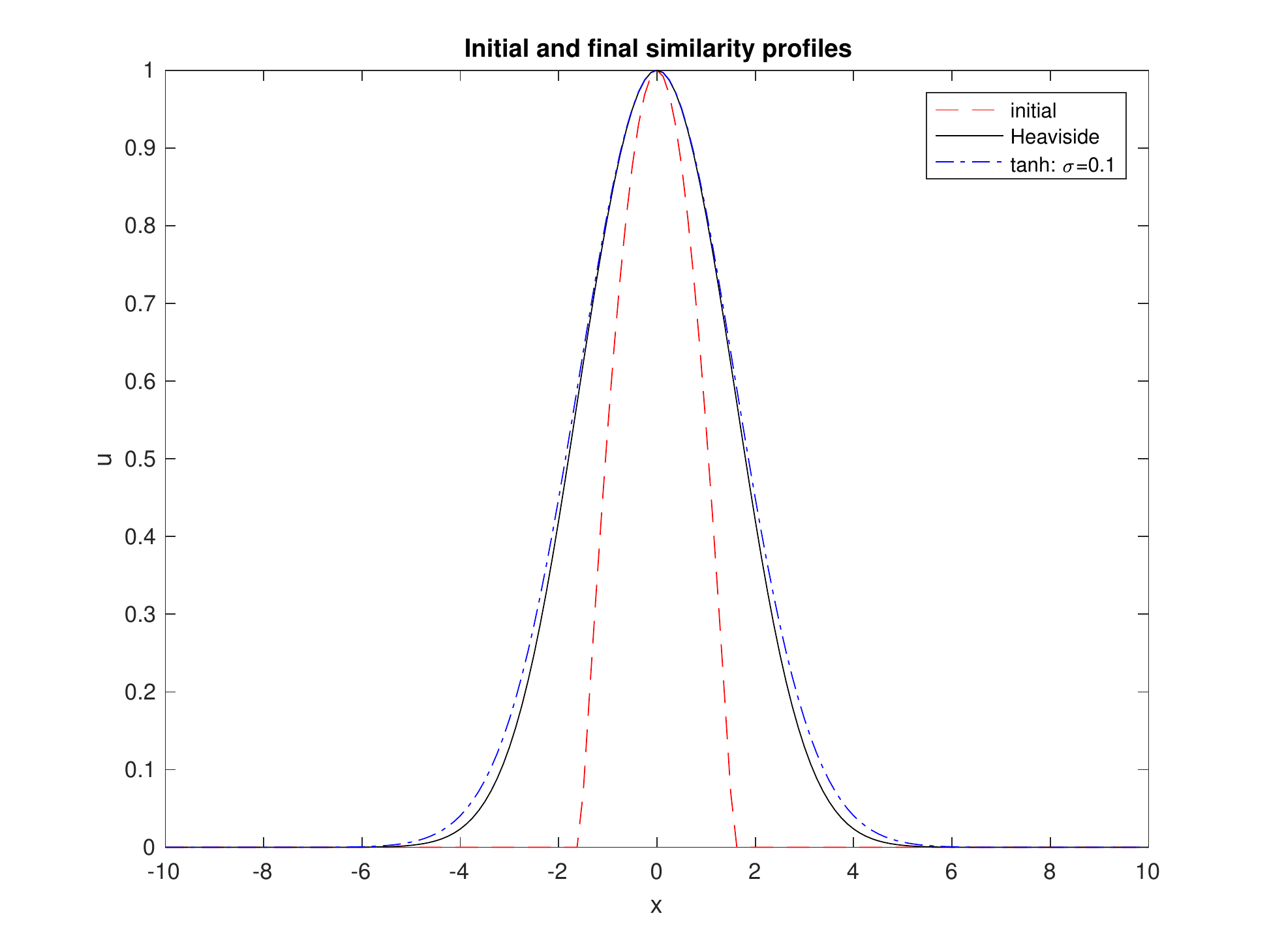} 
(a)\includegraphics[width=2.8in]{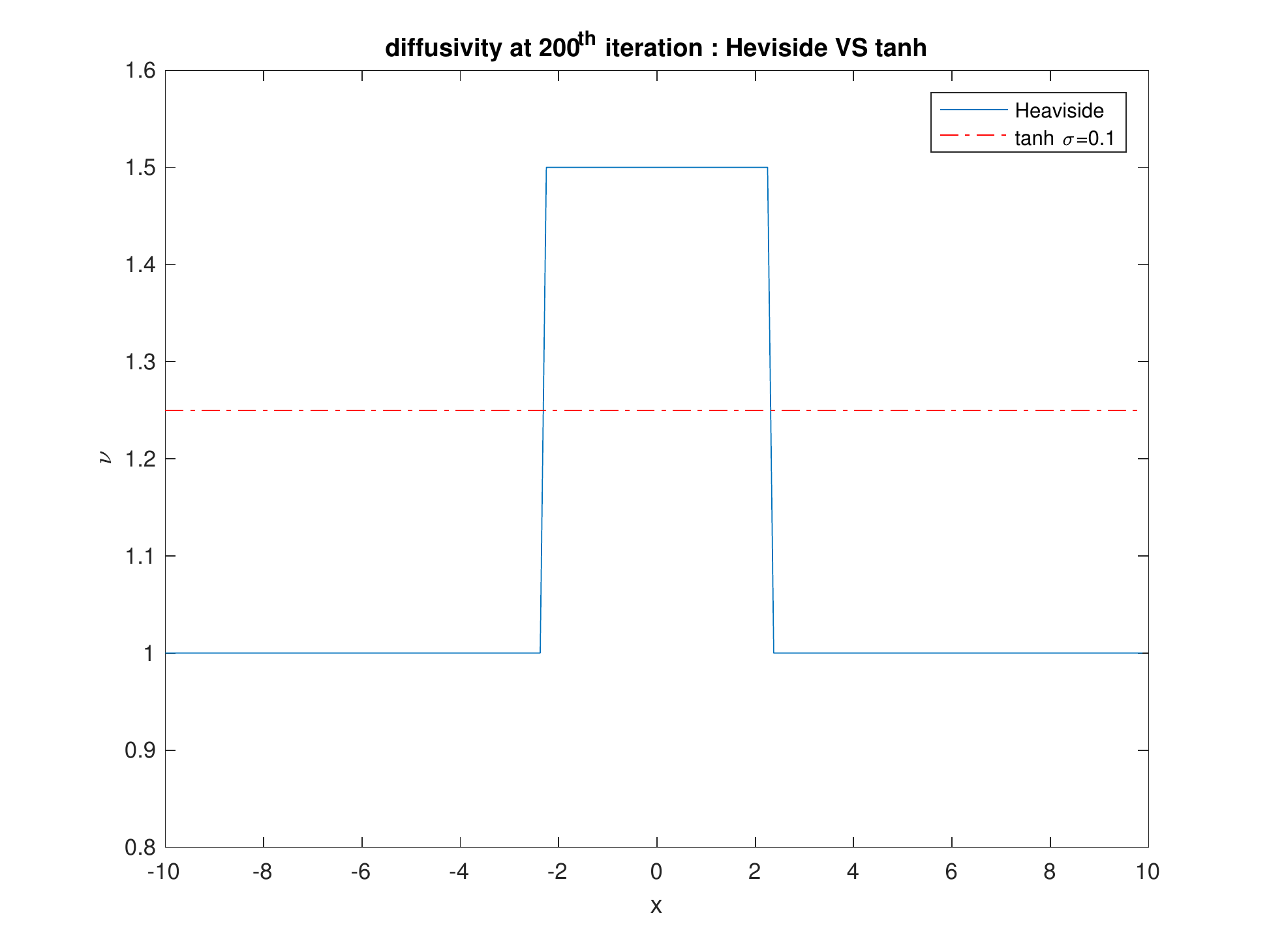} 
\caption{(a) A comparison of the computed asymptotic similarity forms between the discontinuous Heaviside function and the continuous hyperbolic tangent diffusivity function ($\sigma=0.1$). (b) The diffusivity distributions at the end of the $200^{th}$ nRG iterations. The diffusivity distribution (dashed line) is a constant of $1+\frac{\epsilon}{2}$ for the hyperbolic tangent diffusivity function.}
\label{fig:comp_profile)}
\end{figure}

We now repeat the same numerical experiment in Section \ref{sec:barenblatt}, except we replace the Heaviside diffusivity function by the continuous functions in Figure \ref{fig:D(u_t)}. Figure \ref{fig:tanh}(a) shows the initial profile and the computed asymptotic similarity forms. While the three continuous hyperbolic tangent diffusivity functions have different transition bandwidths, the computed final asymptotic similarity forms are visually indistinguishable after 200 nRG iterations. Figure \ref{fig:tanh}(b) shows the time decay parameter $\alpha$ in Eq. (\ref{eq:time_decay}) between the continuous and discontinuous diffusivity function during the nRG calculations. As expected, the time decay parameter approaches to 1/2 for the continuous functions.  Figure \ref{fig:linear} shows the similar calculations to Figure \ref{fig:tanh} for the continuous diffusivity function in Figure \ref{fig:D(u_t)}(b). The results in Figures \ref{fig:tanh} \& \ref{fig:linear} are extremely close, despise the two different types of continuity connecting the jump. Figure \ref{fig:comp_profile)}(a) is the comparison of the computed asymptotic similarity forms between the discontinuous Heaviside and the continuous hyperbolic tangent diffusivities. Figure \ref{fig:comp_profile)}(b) is their diffusivity distributions at the end of the 200 nRG iterations.  As the conjecture in Eq. (\ref{eq:D(Lu_t)}), for the hyperbolic tangent function, the diffusivity distribution approaches to a constant of $1+\frac{\epsilon}{2}$ in the asymptotic regime.

\subsection{The modified diffusion-absorption model ($\lambda > 0$)}

We now consider the full modified diffusion-absorption equation, Eq. (\ref{eq:diff-absorp}). It has been shown that for $D(u_t)\equiv \text{constant},\,\lambda>0,\, p>1+m,\, m\ge 0$ (the absorptive regime), two different regimes of time decay exist, where one is dominated by diffusion and the other by absorption \cite{SGKM95}.  A critical point $p^{*}=p^{*}(m, d)$, where $d$ is the dimension of Eq. (\ref{eq:diff-absorp}), separates the two regimes. In the case, $p > p^{*}$, the time decay parameter $\alpha$ becomes a constant function of $p$, indicating that for these absorptive exponents, the model equation is an irrelevant perturbation to the diffusion equation, which means a regime of diffusion dominant. On the other hand, for $1+m<p<p^{*}$, the absorptive time decay is given by 
\begin{equation}\label{eq:time_decay1}
\alpha = \frac{1}{p-1},
\end{equation}
which indicates that the time decay is strongly influenced by absorption and the longtime behavior of the solution does not correspond to the diffusion equation. Therefore, the model equation is a relevant perturbation to the diffusion equation, and so is in the absorption dominant regime. Moreover, for the marginal case, $p=p^{*}$, the absorptive time decay equals the diffusive time decay, or
\begin{equation}
\frac{1}{p-1} = \frac{d}{md+2},
\end{equation}
and hence the critical value of $p$ is 
\begin{equation}\label{eq:time_decay2}
p^{*}=m+\frac{2}{d} +1.
\end{equation}
Note that for $p>p*$ the time decay $\alpha$ is a constant 
\begin{equation}\label{eq:time_decay3}
\alpha_{c} = \frac{1}{p^{*}-1}.
\end{equation}
Detailed results and their derivations for the model equation in the absorptive regime can be found in \cite{SGKM95}.


For the rest of this section, we will consider the cases $m=0$ and $m=1$ with the normalized coefficient $\lambda =1$. We will investigate the time decay parameter $\alpha$, as a function of the absorptive exponent $p$ for $D(u_t)$ that is beyond a constant. In particular, we will illustrate the different behaviors of $\alpha$ vs $p$ between the discontinuous and continuous diffusivity functions. Similar to Eq. (\ref{eq:scaled_Baren}), the scaled  Eq. (\ref{eq:diff-absorp}) for the $n^{th}$ RG iteration is
\begin{equation}\label{eq:scaled_absorption}
\begin{split}
&(u_n)_{t}=D\left(L^{-n(\bar{\alpha}_{n}+1)}(u_{n})_t\right) L^{-n(\bar{\alpha}_n m +2 \bar{\beta}_n -1)}((u_n)^{m+1})_{xx}-L^{-n(\bar{\alpha}_n(p-1)-1)}(u_n)^{p},\,\,t >1,\\
\text{I. C.\,\,:}\quad&u_n(x,1) = f_n(x),
\end{split}
\end{equation}
We consider two types of diffusivity functions in this section: (i) the discontinuous Heaviside function and (ii) the continuous hyperbolic tangent function. To incorporate the nonlinear diffusion ($\partial_{xx}u^2$ for the case $m=1$) and the absorption term, we will use the second-order explicit centered-difference discretization for the spatial derivative and the Euler's scheme for the time evolution. The constrain of the ratio of temporal and spatial step sizes, $\nu\frac{\Delta t}{(\Delta x)^2} \le 1/2$, will be enforced to ensure the stability, where $\nu$ is the diffusivity. The nonlinear diffusivity $D(u_t)$ is linearized the same way as before. 

\vskip 12pt

\noindent 
{\bf Case I: $m=0$:}  
Using the nRG algorithm, we numerically study $\alpha$ versus $p$ for diffusivities beyond a constant function. In particular, we study a discontinuous $D(u_t)$ that is the Heaviside function described in Eq. (\ref{eq:Heaviside}) with $\epsilon=0.5,\, 0$ and $-0.5$. We also study the continuous hyperbolic tangent profile described in Eq. (\ref{eq:D(u_t)}) with $\sigma=0.1$ and $\epsilon=0.5,\, 0$ and $-0.5$.  We choose $2\le p \le 5$ with the increment $\Delta p=0.1$.  We remark  that in order to study the effect of different diffusivities, it will be a good idea to keep the magnitude of the diffusivity unscaled throughout the RG iterations. From Eq. (\ref{eq:scaled_absorption}), if we choose $\beta_{n}=1/2$ for every $n$, the magnitude of the diffusivity will remain unscaled in each nRG iteration.

\begin{figure}[h]
\centering
(a)\includegraphics[width=2.8in]{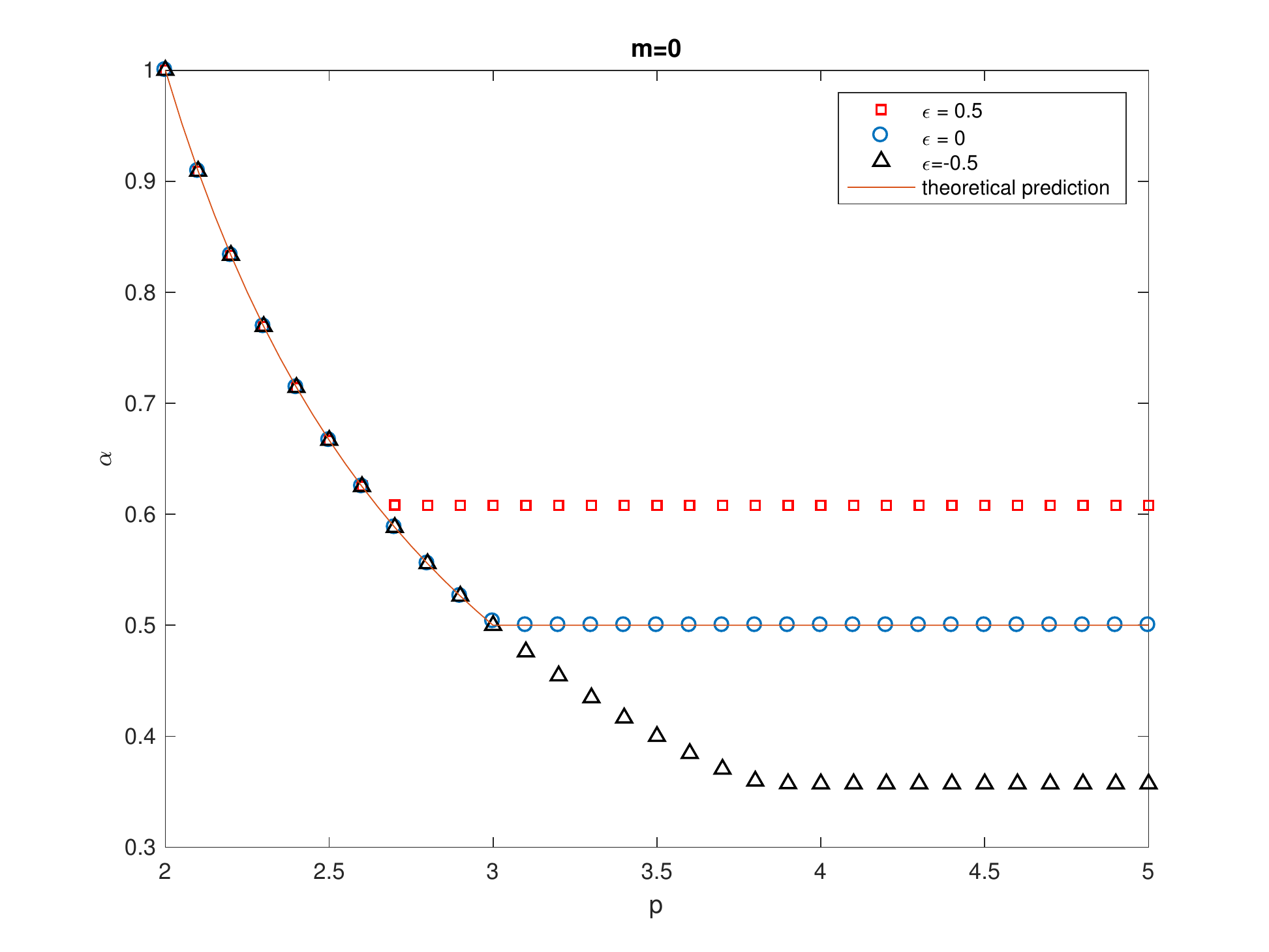} 
(b)\includegraphics[width=2.8in]{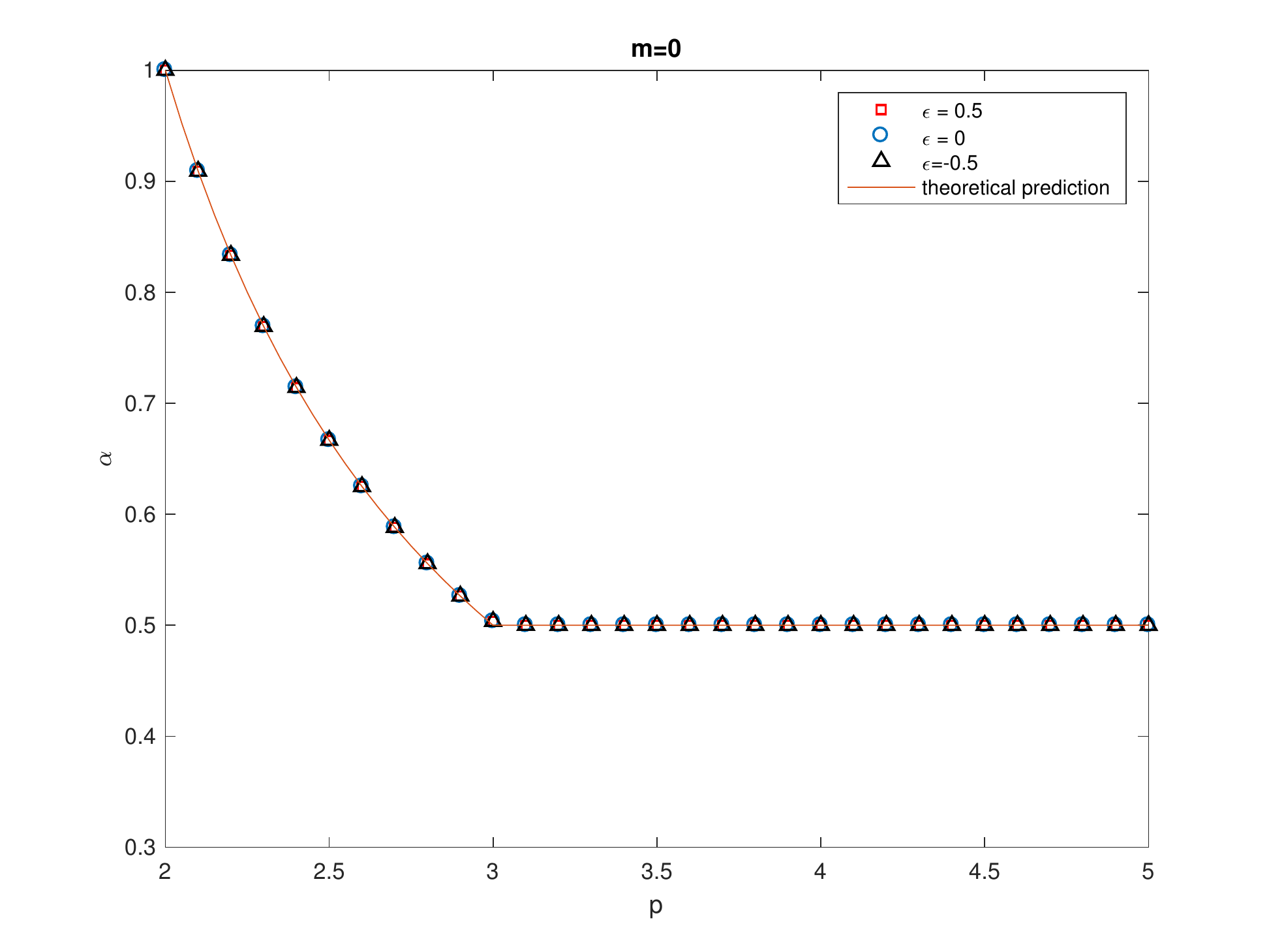}
\caption{$m=0$. (a) $D(u_t)$ is the Heaviside function (\ref{eq:Heaviside}) with $\epsilon=0.5,\, 0$ and $-0.5$. The solid line is the theoretically predicted $\alpha$ values. (b) The same as (a), except $D(u_t)$ is the continuous hyperbolic tangent function (\ref{eq:D(u_t)}) with $\sigma=0.1$.}
\label{fig:p_vs_alpha_m=0}
\end{figure}

Figure \ref{fig:p_vs_alpha_m=0}(a) is the plot for the computed time decay parameter $\alpha$ as a function of $p$ for Heaviside diffusivity with three different jump $\epsilon$-values. The solid line is the theoretically predicted $\alpha$ values, for which, by Eq. (\ref{eq:time_decay2}), the critical value $p^{*}=3$. Hence for $p<3$, $\alpha$ obeys Eq. (\ref{eq:time_decay1}) and for $p\ge 3$, $\alpha = \frac{1}{2}$ by Eq.(\ref{eq:time_decay3}). The numerically computed $\alpha$ values corresponding to the constant diffusivity, $D(u_t)\equiv 1$ ($\epsilon=0$) are marked by the circles. The computed values agree with the theoretical prediction. The computed $\alpha$ values for the diffusivities that are the discontinuous Heaviside function with $\epsilon=0.5$ and $-0.5$ are marked by the squares and the triangles, respectively.  The results suggest that the $\alpha$ values obey Eq. (\ref{eq:time_decay1}) as a function of $p$, until a critical $p$ value is reached. For $p$ that is larger than the critical $p$ value, the time decay parameter $\alpha$ is a constant. Moreover, different $\epsilon$ values give rise to different critical $p$ values. 

Figure \ref{fig:p_vs_alpha_m=0}(b) is the results for diffusivity that is the continuous hyperbolic tangent profile. Unlike the Heaviside diffusivity, the time decay $\alpha$ values for continuous tangent profiles with different jumps $\epsilon$ obey the theoretical prediction for the case of constant diffusivities, or $\alpha=\frac{1}{p-1}$. This behavior is understandable, thanks to the back-of-the-envelope calculation in Eq.(\ref{eq:D(Lu_t)}).

Finally we note that for all our numerical computations, the time is integrated from $L=1$ to $L=2$ with $\Delta t=5\times10^{-4}$, while the computational domain is $-10\le x\le 10$ with $\Delta x=0.125$. The initial condition is described in Eq. (\ref{eq:IC_experiment}). For each simulation, the number of RG iterations is 200.

\vskip 12pt

\begin{figure}[h]
\centering
(a)\includegraphics[width=2.8in]{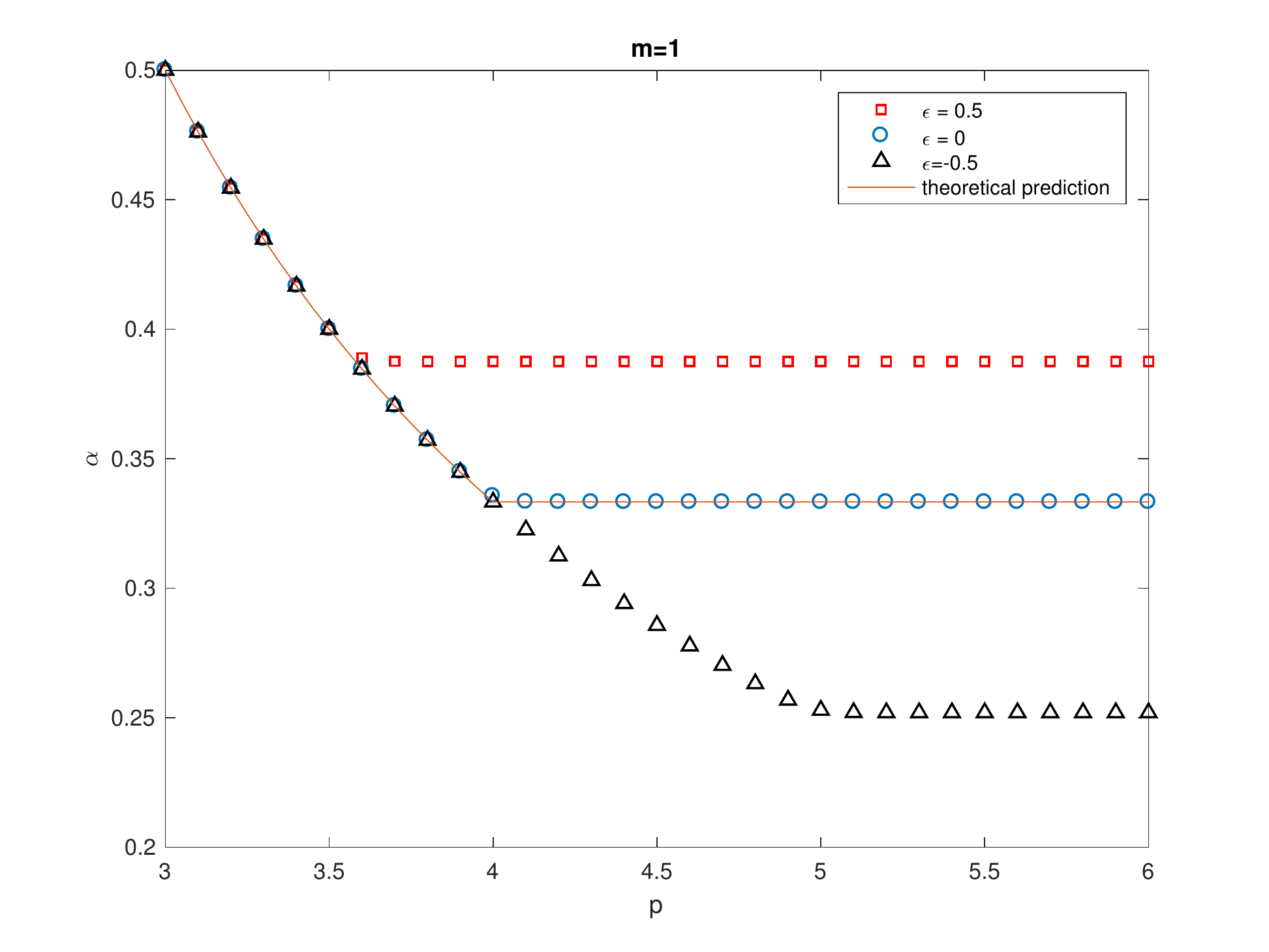} 
(b)\includegraphics[width=2.8in]{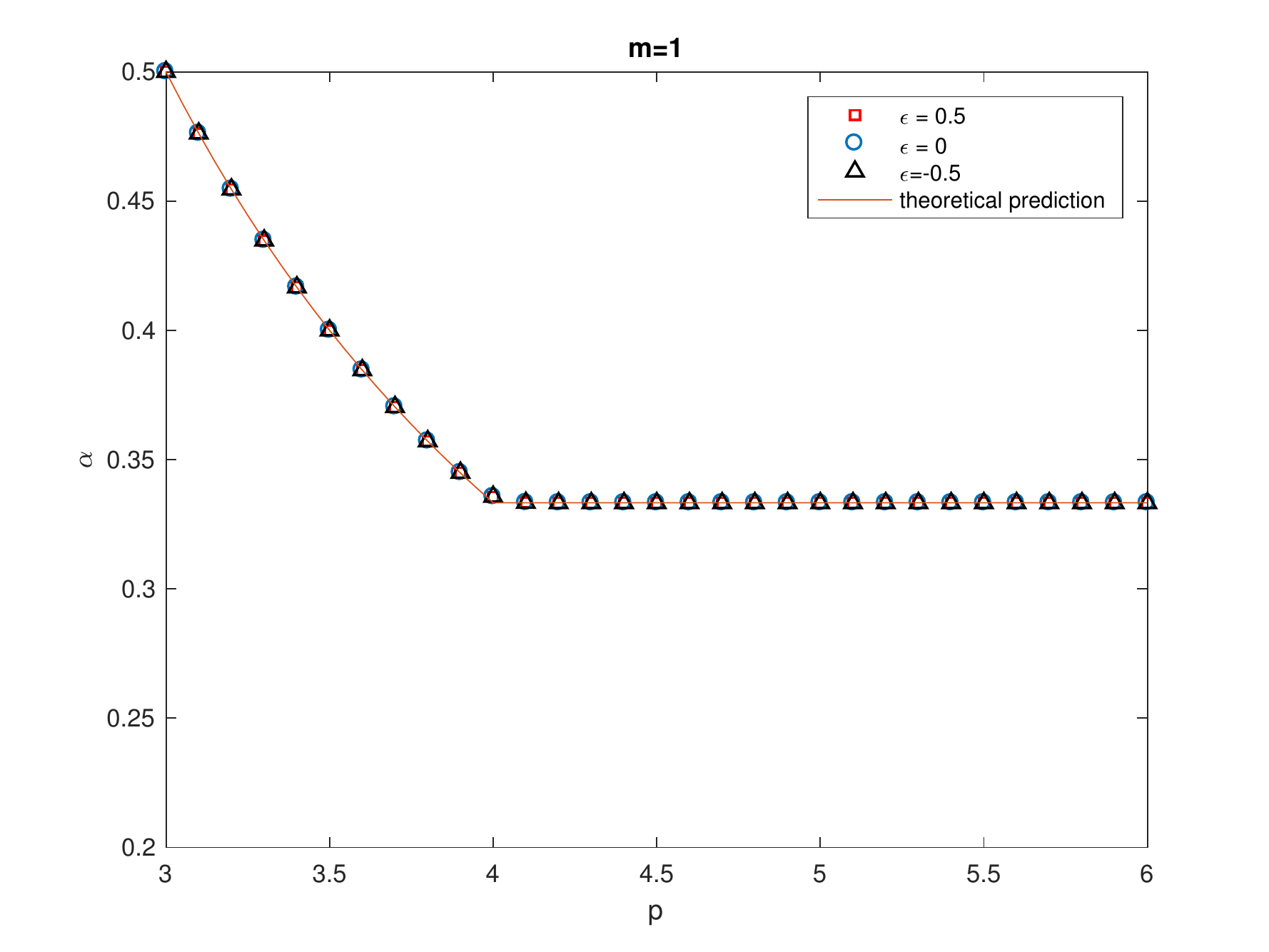}
\caption{$m=1$. (a) $D(u_t)$ is the Heaviside function (\ref{eq:Heaviside}) with $\epsilon=0.5,\, 0$ and $-0.5$. The solid line is the theoretically predicted $\alpha$ values. (b) The same as (a), except $D(u_t)$ is the continuous hyperbolic tangent function (\ref{eq:D(u_t)}) with $\sigma=0.1$.}
\label{fig:p_vs_alpha_m=1}
\end{figure}

\noindent
{\bf Case II: $m=1$:} When $m\ne 0$, based on Eq.(\ref{eq:scaled_absorption}) the scaling factor for the diffusivity depends on both $\alpha_n$ and $\beta_n$. If an unscaled diffusivity is desired, the spacial scaling factor $\beta_n$ should be calculated by $\beta_n = 1-m \alpha_n$, after $\alpha_n$ is computed at the end of each RG iteration. Figure \ref{fig:p_vs_alpha_m=1} shows that the nonlinear diffusion with absorption exhibits the similar behavior as that of the linear case, shown in {\bf Case I}.

The results of the experiments in {\bf Case I \& II} suggest the following conjecture:  
\begin{conjecture}
For the diffusion-absorption model (\ref{eq:diff-absorp}) with $\lambda > 0$, the behavior of the time decay parameter $\alpha$ at the asymptotic regime for the discontinuous Heaviside diffusivity is similar to that for the constant diffusivity, and there exists a critical $p$ value for the Heaviside diffusivity with the jump $\epsilon$, i.e. $p^{*} = 1+\frac{1}{\alpha(\epsilon)}$, where $\alpha(\epsilon)$ depends on $\epsilon$ and is constant for $p\ge p^{*}$.  
\end{conjecture}

\subsection{Marginal case: vanishing pre-factor} 
As we indicated in Section \ref{sec:Burgers} that if the self-similar solution has only power law decay, then we can monitor the sequence of $A_n$, defined in Eq. (\ref{eq:An_Bn}), and expect that $A_n$ will converge to some constant $A\ne 0$.  For the diffusion-absorption model (\ref{eq:diff-absorp}) with a constant diffusivity, $m=0$, and $\lambda > 0$, Bricmont and Kupiainen \cite{BKL94} shows that the longtime self-similar solution for the marginal case, $p=3$, is in the form of
\begin{equation}
u(x, t)\sim \left(\frac{\lambda t \log t}{2\sqrt{3}} \right)^{-\frac{1}{2}}\phi(x t^{-\frac{1}{2}}),
\end{equation}
where $\phi$ is a Gaussian distribution. This suggests that, using Algorithm 1, even though our nRG iterations could successfully capture the decay exponent, $\alpha=\frac{1}{2}$ (as indicated in Figure \ref{fig:p_vs_alpha_m=0} for $\epsilon=0$) and $\beta=1/2$ for keeping the diffusivity unscaled, the sequence of the pre-factor $A_n$ should continue approaching to zero, because
\begin{equation}
\lim_{n\rightarrow\infty} A_n \sim  \left(\frac{\lambda \log t}{2\sqrt{3}} \right)^{-\frac{1}{2}} \rightarrow 0,\quad \text{as}\,\, t \rightarrow\infty.
\end{equation}
Figure \ref{fig:prefactors_porous}(a) shows that for linear diffusion ($m=0$) with constant diffusivity, the computed pre-factor $A_n$ continues to approach to zero after 50,000 nRG iterations for $p=3$, while $A_n$ quickly settles into a nonzero constant for $p=3.05$ and $p=2.95$ that are slightly deviated away from the critical $p$ value, $p=3$. 

\begin{figure}[h]
\centering
(a)\includegraphics[width=2.8in]{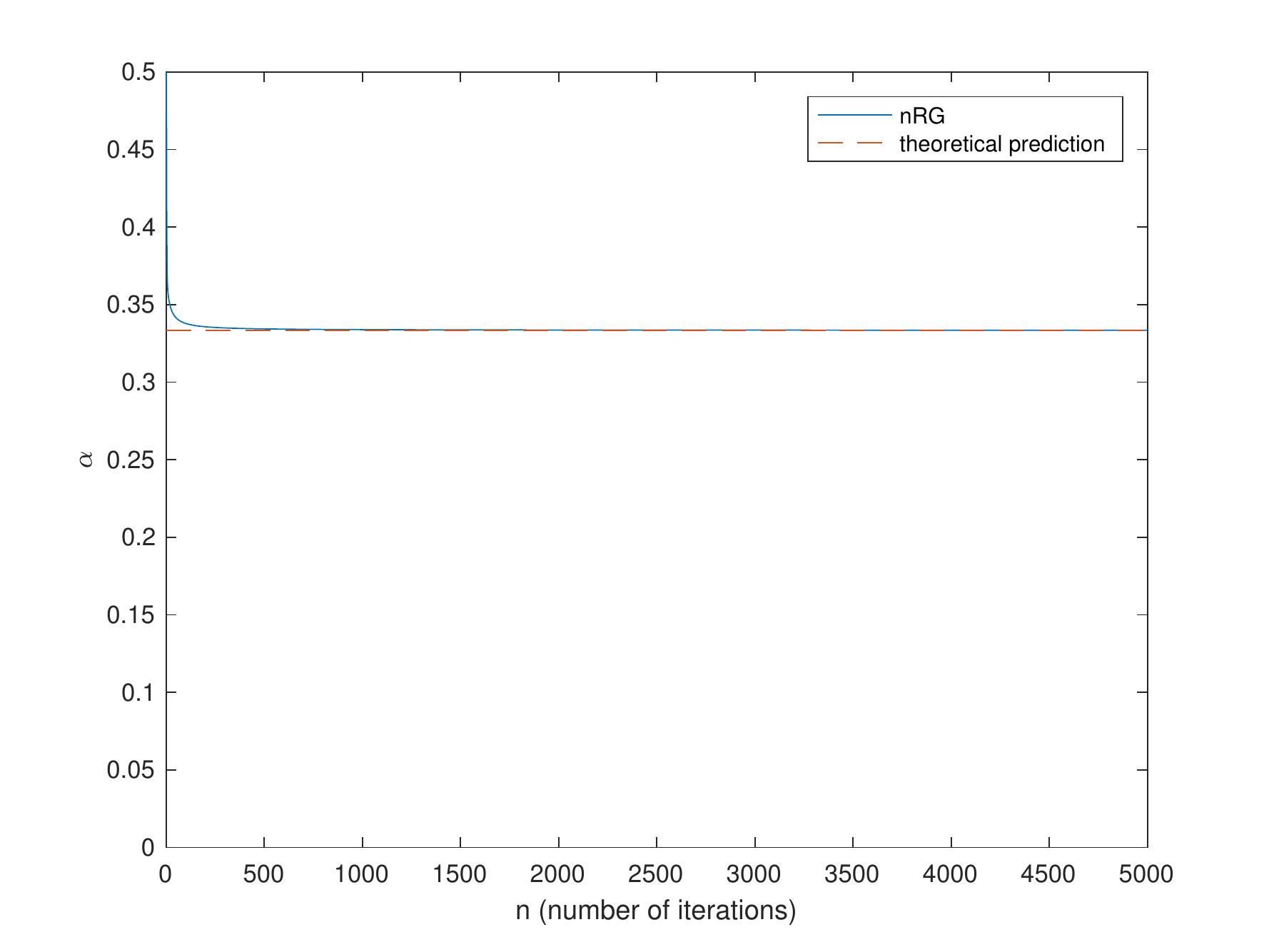} 
(b)\includegraphics[width=2.8in]{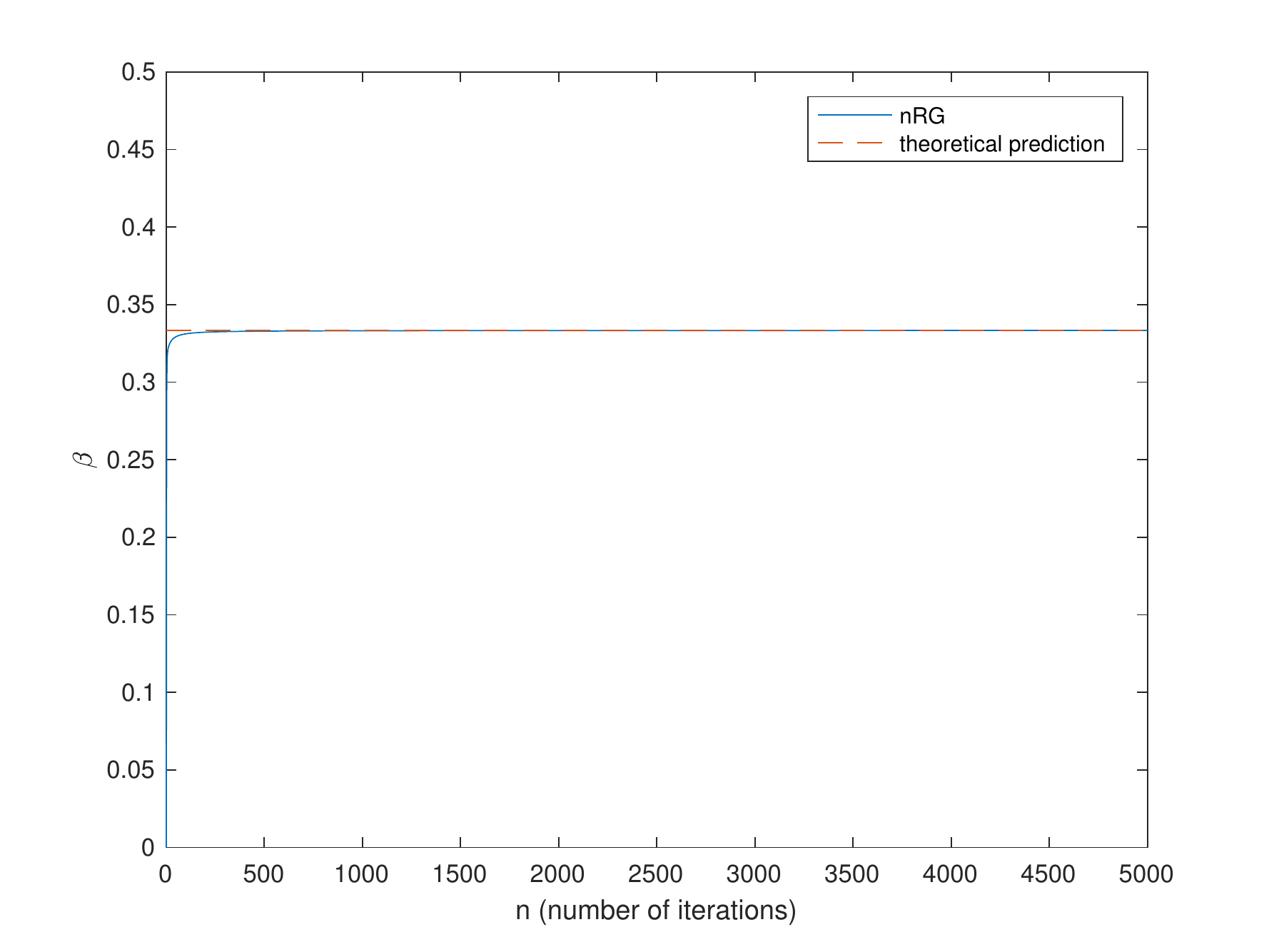}
\caption{Comparison of the time decay factors captured by the nRG and that of the theoretical prediction. (a) $\alpha\rightarrow 1/3$, (b) $\beta\rightarrow 1/3$.}
\label{fig:decayfactors}
\end{figure}

\begin{figure}[h]
\centering
(a)\includegraphics[width=2.8in]{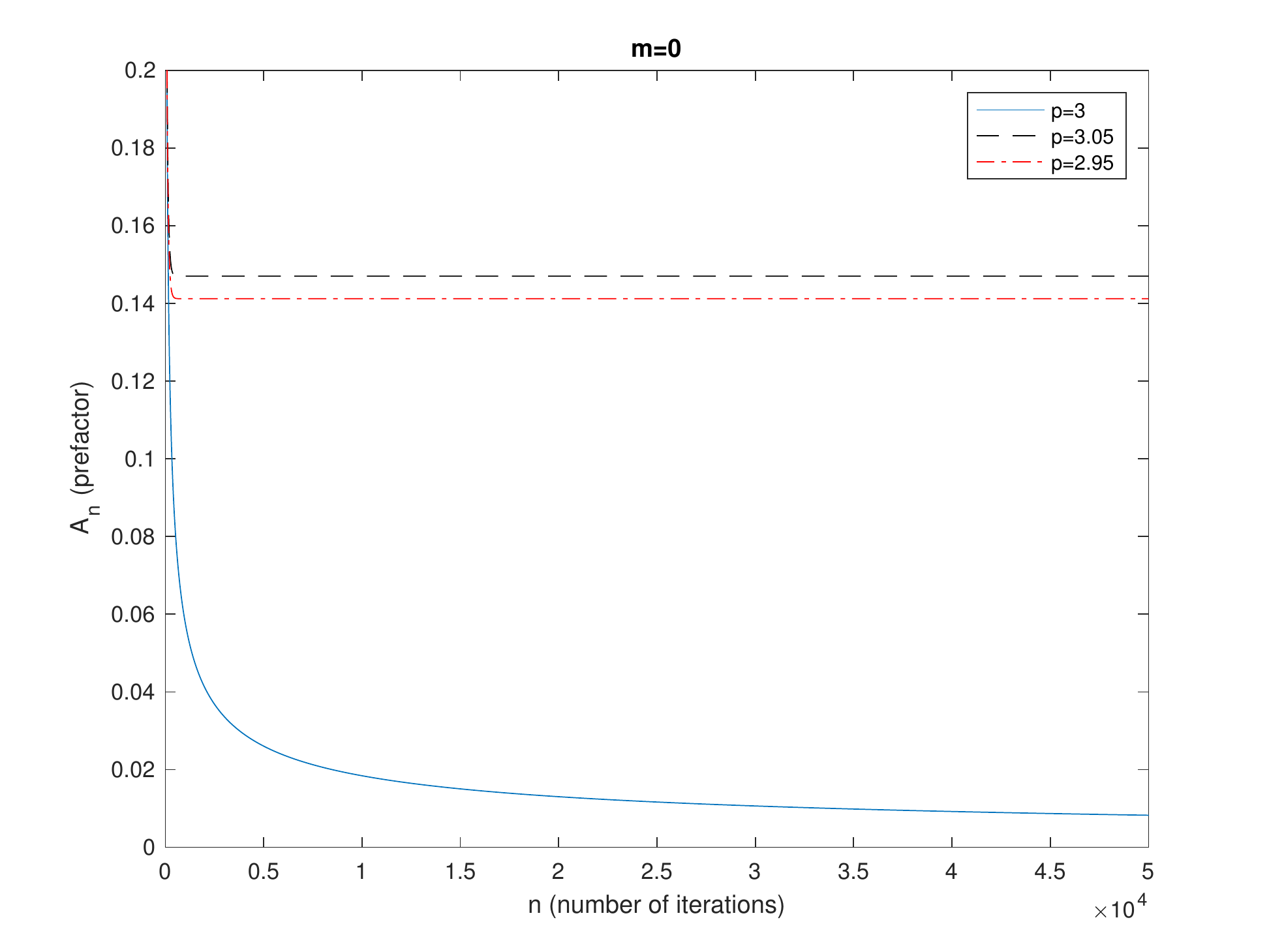} 
(b)\includegraphics[width=2.8in]{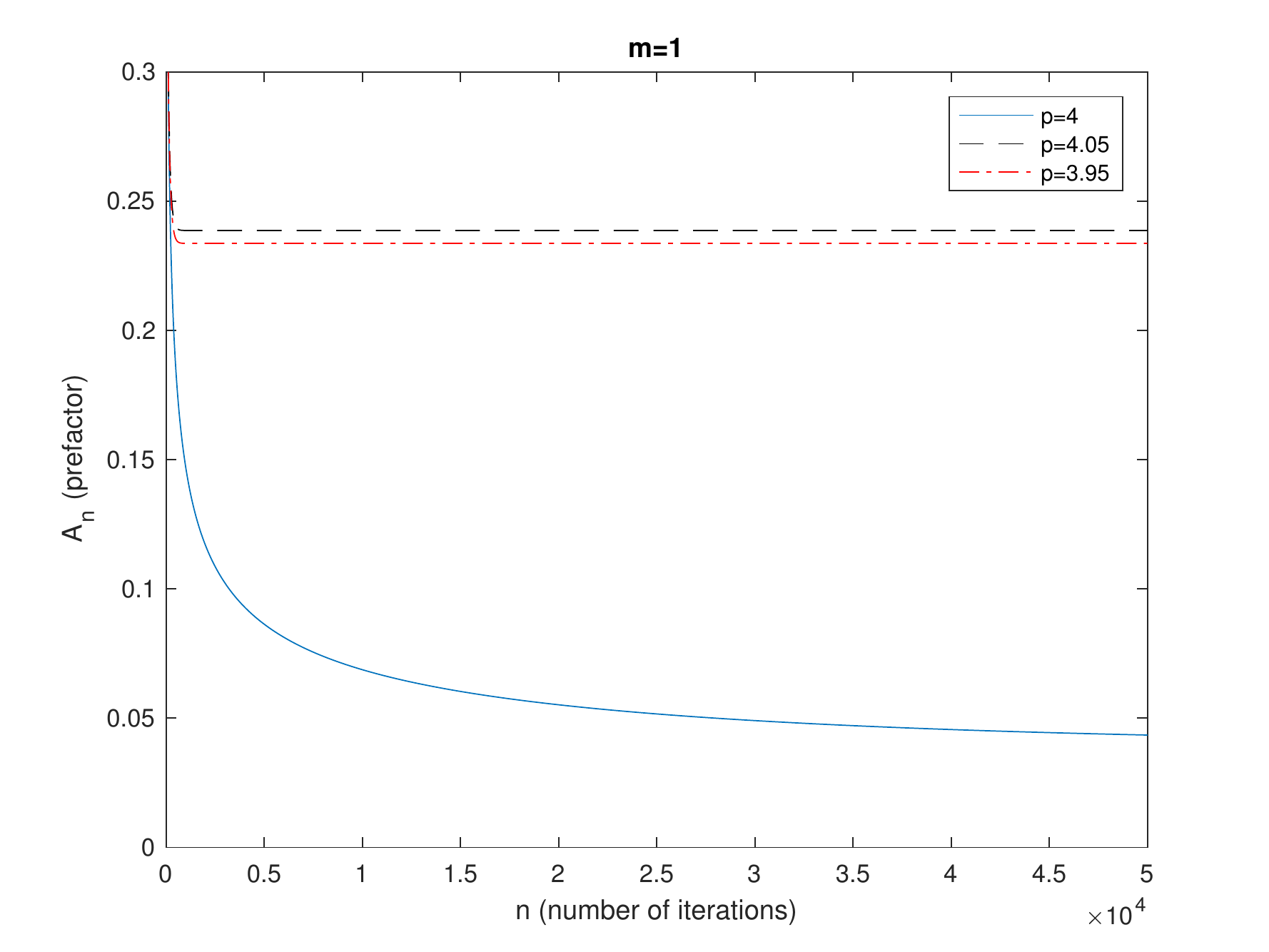}
\caption{(a) Linear diffusion ($m=0$) with constant diffusivity. $A_n\rightarrow 0$ for the marginal case $p=3$. (b) Nonlinear diffusion ($m=1$) with constant diffusivity. $A_n\rightarrow 0$ for the marginal case $p=4$.}
\label{fig:prefactors_porous}
\end{figure}

It is worth pointing out that for the nonlinear diffusion ($m=1$) with constant diffusivity, we observed exactly the same behavior for the marginal case of the critical value $p=4$. Qi and Liu \cite{QL04} showed that, for the diffusion-absorption model (\ref{eq:diff-absorp}) with a constant diffusivity, $m=1$, and $\lambda=1$, for the marginal case ($p=4$), the similarity solution is 
\begin{equation}\label{eq:sim_m=1}
u(x, t) \sim (t \log t )^{-\frac{1}{3}} \phi\left(\frac{(\log t)^{1/6} x}{t^{1/3}}\right),\quad t\rightarrow\infty,
\end{equation}
for initial data $u_0$ that satisfies 
\begin{equation}
\lim_{|x|\rightarrow\infty} \sup |x|^{k} u_0 < \infty, \quad k>1.
\end{equation}
This implies that for our compactly supported initial data if the time decay factors captured by our nRG calculation approach to the theoretical prediction, i.e. $\alpha\rightarrow 1/3$ and $\beta\rightarrow 1/3$, then we have
\begin{equation}\label{eq:prefactor_m1}
\lim_{n\rightarrow\infty}A_n \sim (\log t )^{-\frac{1}{3}} \rightarrow 0,\quad \text{as}\,\, t \rightarrow\infty.
\end{equation}
In Figure \ref{fig:decayfactors}, we show that the time decay factors captured by the nRG algorithm are consistent with the theoretical values predicted in Eq. (\ref{eq:sim_m=1}), and in Figure \ref{fig:prefactors_porous}(b) we show that the computed pre-factor $A_n$ continues to approach to zero after 50,000 nRG iterations for $p=4$, as indicated in Eq. (\ref{eq:prefactor_m1}), while $A_n$ quickly settles into a nonzero constant for $p=4.05$ and $p=3.95$ that are slightly deviated away from the critical $p$ value, $p=4$.


We comment that the numerical experiments in Figure \ref{fig:prefactors_porous} use the compactly supported initial condition (\ref{eq:IC_experiment}) with the constant diffusivity $D\equiv 1$ and the normalized absorption coefficient $\lambda =1$. The parameters for the nRG iterations are $L=2$, $\Delta x=0.1$, and $\Delta t=10^{-4}$. Similar to the previous experiments, the explicit second-order discretization is applied to the spatial derivative and the Euler's method is used for the time evolution.  

The peculiar phenomenon observed in Figure \ref{fig:prefactors_porous} motivates us to modify our Algorithm 1 in the next section, in order to capture the hidden logarithmic time decay exponent illustrated in this section. 

\section{Cubic autocatalytic chemical reaction system}\label{sec:cubic_autocatalytic}

The nRG procedure described in Algorithm 1 assumes that the asymptotic solutions decay or expand at a rate obeying the power law. However, there are differential equations (or systems of differential equations) whose solutions decay at a rate other than the power law, such as the logarithmic decay discussed in the previous Section. For these solutions, the aforementioned Algorithm 1 is not sufficient to capture the correct decay at the asymptotic region. Nevertheless, the procedure could provide sufficient information that allows us to modify the current algorithm to capture the similarity solutions of those equations.

To illustrate the modification, we consider the Cauchy problem of the chemical reaction system
\beq\label{eq:chem-reaction}
\begin{split}
u_t=u_{xx} - u^{p}v^{q},\\
v_t=dv_{xx} + u^{p}v^{q},
\end{split}
\eeq
where $p+q=3, 1\le p, q\le 2,$ and $d>0$. This system arises as a model for cubic autocatalytic chemical reactions of the type
\beq
pK_1+qK_2\longrightarrow 3K_2
\eeq
with isothermal reaction rate proportional to $u^{p}v^{q}$, where $u$ is the concentration of reactant $K_1$ and $v$ is the concentration of auto-catalyst $K_2$ \cite{LQ03}. The system is subject to the initial data $u(x, 0)=a_1(x)$ and $v(x,0)=a_2(x)$, where $a_1, a_2 \ge 0$ and $a_1, a_2 \in L^{1}(\mathbb{R})\cap L^{\infty}(\mathbb{R})$. The above system has been applied for modeling thermal-diffusive combustion problems \cite{BKX96} or mathematical biology \cite{FHMW02}.  

The similarity solutions of this system, based on different values of $p$ and $q$, are investigated in the papers by Bricmont et al. \cite{BKX96} and Li and Qi \cite{LQ03}. For $p=1$ $q=2$, Bricmont et al. show that 
\beq
\begin{split}
& t^{1/2+E_{A}}u(\sqrt{t}x, t) \longrightarrow B\psi_{A}(x),\\
& t^{1/2}v(\sqrt{t}x, t) \longrightarrow A\phi_{d}(x),
\end{split}
\eeq
as $t\rightarrow \infty$. Here the total mass $A=\int_{\mathbb{R}}\left(a_1(x) + a_2(x)\right)dx$ is conserved along time. $B$ depends continuously on $(a_1, a_2)$, and the extra decay power $E_{A}$ in time is due to the critical cubic nonlinearity of the system \cite{BKX96}. $\phi_{d}$ is the Gaussian
\beq\label{eq:gaussian}
\phi_{d}(x) = \frac{1}{\sqrt{4\pi d}}e^{-x^2/4d}.
\eeq
Li and Qi extend the above result by considering the values $1< p, q <2$ and $p+q=3$. The nontrivial initial data $a_i\ge 0$, for  $i=1, 2$ are the same as before, whereas the total mass $A$ defined as before is positive. Li and Qi show that 
\beq\label{eq:log_decay}
\begin{split}
& \sqrt{t}(\log t)^{1/(p-1)}u(\sqrt{t}x, t)\longrightarrow B\phi_1(x),\\
& \sqrt{t}v(\sqrt{t}x, t)\longrightarrow A\phi_d(x),
\end{split}
\eeq
as $t\rightarrow\infty$, where 
\beq\label{eq:B}
B=\left(\frac{4\pi d^{q/2}(p+q/d)^{1/2}}{(p-1)A^{q}} \right)^{1/(p-1)},
\eeq
and $\phi_1$ is $d=1$ in Eq. (\ref{eq:gaussian}). The peculiar phenomena of the similarity solution (\ref{eq:log_decay}) is that the $u$-component contains two decays, the regular power-law decay and a logarithmic decay. We illustrate below that the nRG algorithm described in section \ref{sec:nrg} is not sufficient to capture the second decay. However, the procedure will provide a clue for the existence of the second decay, and allows us to design a nRG procedure to capture the similarity solutions in Eq. (\ref{eq:log_decay}). We start with the regular nRG procedure stated in Section \ref{sec:nrg}, i.e. the scaling for $t$ and $x$ is the same as that in Eq. (\ref{eq:scale_tx}), which results in $u$ and $v$ being scaled by 
\beq\label{eq:sys_scaling}
\begin{split}
u_L(\tilde{x},\tilde{t}) &= L^{\alpha_1}\,u(x,t) = L^{\alpha_1}\,u(L^{\beta_1}\tilde{x},L\tilde{t}),\\
v_L(\tilde{x},\tilde{t}) &= L^{\alpha_2}\,v(x,t) = L^{\alpha_2}\,v(L^{\beta_2}\tilde{x},L\tilde{t}).\\
\end{split}
\eeq
Thus the scalings result in a system of PDEs (dropping the subscript $L$ and $\,\tilde{}$ )
\beq\label{eq:sys-scaled}
\begin{split}
u_t & = L^{-2\beta_1+1}u_{xx}-L^{-p\alpha_1-q\alpha_2+\alpha_1+1}u^{p}v^{q},\\
v_t & = L^{-2\beta_2+1} d v_{xx}+L^{-p\alpha_1-q\alpha_2+\alpha_2+1}u^{p}v^{q}.\\
\end{split}
\eeq
Similar to the Burgers equation, we choose to keep the diffusion coefficients invariant. Thus $\beta_1=\beta_2=1/2$ at all time, whereas $\alpha_1$ and $\alpha_2$ are computed by step (2) in the nRG algorithm described in Section \ref{sec:nrg}, respectively. With this choice of $\beta_1$ and $\beta_2$, the scaled PDE at the $n^{th}$ iteration is 
\beq\label{eq:power_decay_pde}
\begin{split}
(u_n) &= (u_n)_{xx} - L^{n}\left(L^{-n\bar{\alpha}_{1, n}}  \right)^{p-1}\left(L^{-n\bar{\alpha}_{2, n}}\right)^{q}(u_n)^{p}(v_n)^{q},\\
(v_n) &= d (v_n)_{xx} + L^{n}\left(L^{-n\bar{\alpha}_{1,n}}  \right)^{p}\left(L^{-n\bar{\alpha}_{2, n}}\right)^{q-1}(u_n)^{p}(v_n)^{q},
\end{split}
\eeq
where $\bar{\alpha}_{1, n}$ and $\bar{\alpha}_{2, n}$ are defined as $\bar{\alpha}_n$ in Section \ref{sec:seq}. At this stage, we assume that the power-law scaling, based on the hypothesis, is
\beq\label{eq:similarities}
\begin{split}
u(x, L^{n}) \sim \frac{A_{u}}{L^{n/2}}\phi_{u}(\frac{x}{L^{n/2}}),\\
v(x, L^{n}) \sim \frac{A_{v}}{L^{n/2}}\phi_{v}(\frac{x}{L^{n/2}}),
\end{split}
\eeq
where $A_{u}$ and $A_{v}$ are non-zero constant. The nRG iteration, based on the power-law decay assumption, in principle will show  $A_{u, n}=L^{n(\alpha_{1, n} - \bar{\alpha}_{1, n})}\sim A_{u}$, and $A_{v, n}=L^{n(\alpha_{2, n} - \bar{\alpha}_{2, n})}\sim A_{v}$ (cf. Eq. (\ref{eq:An_Bn})), as $n\rightarrow\infty$. Unfortunately (or fortunately), this is not the case. The numerical experiment, in fact, shows that $A_{u, n}\rightarrow 0$, while $A_{v, n}\ne 0$ and converges to some constant proportional to the total mass $A_{total}$, as $n\rightarrow\infty$. Based on this result, we conjecture that the $v$-component follows the power-law decay, similar to the Burgers equation, while the  $u$-component exists a hidden decay that is not captured by solely assuming the power-law decay. 

\subsection{A numerical experiment}\label{sec:no-decay}
\begin{figure}[bhtp]
\centering
(a)\includegraphics[width=2.8in]{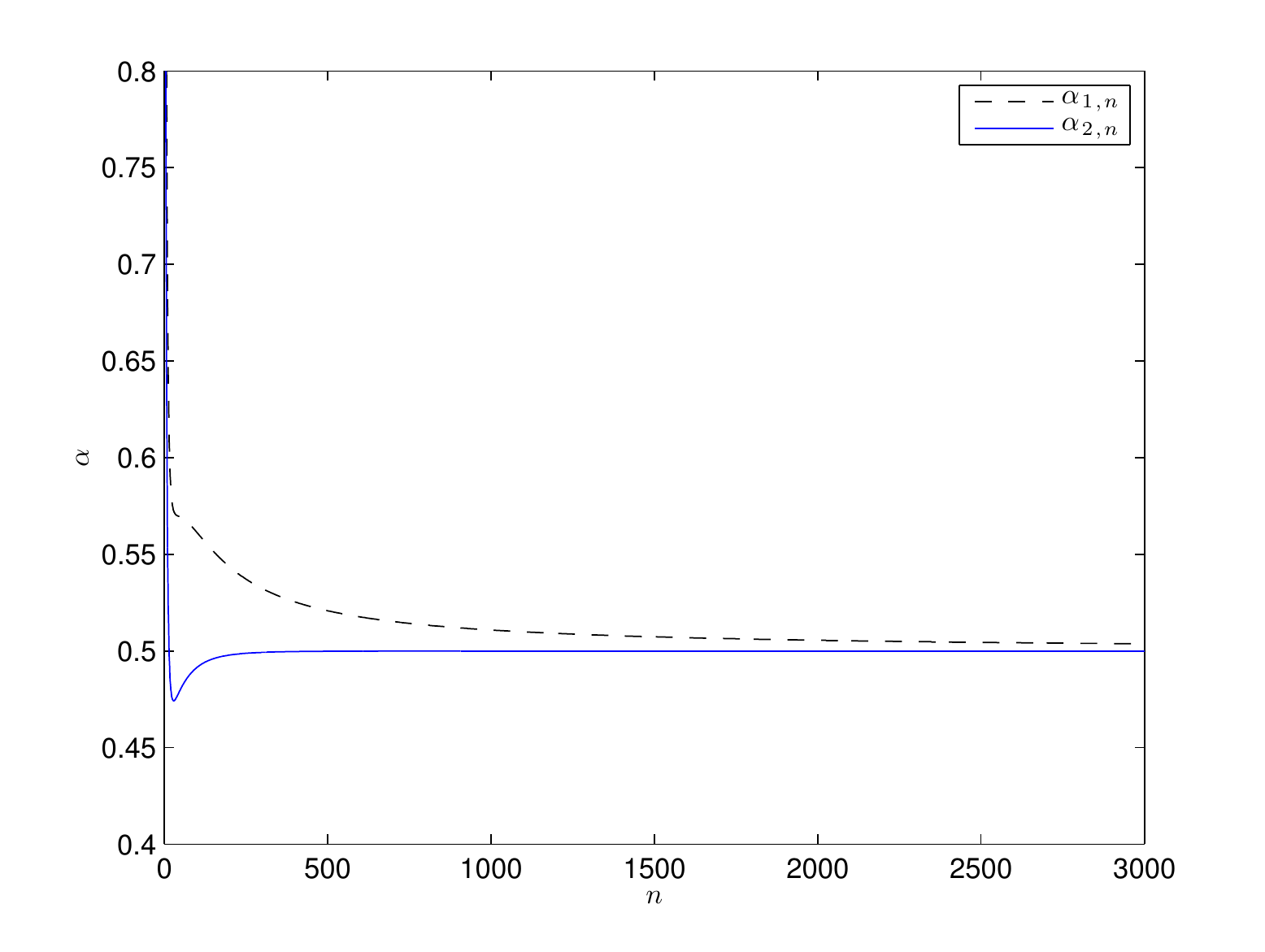} 
(b)\includegraphics[width=2.8in]{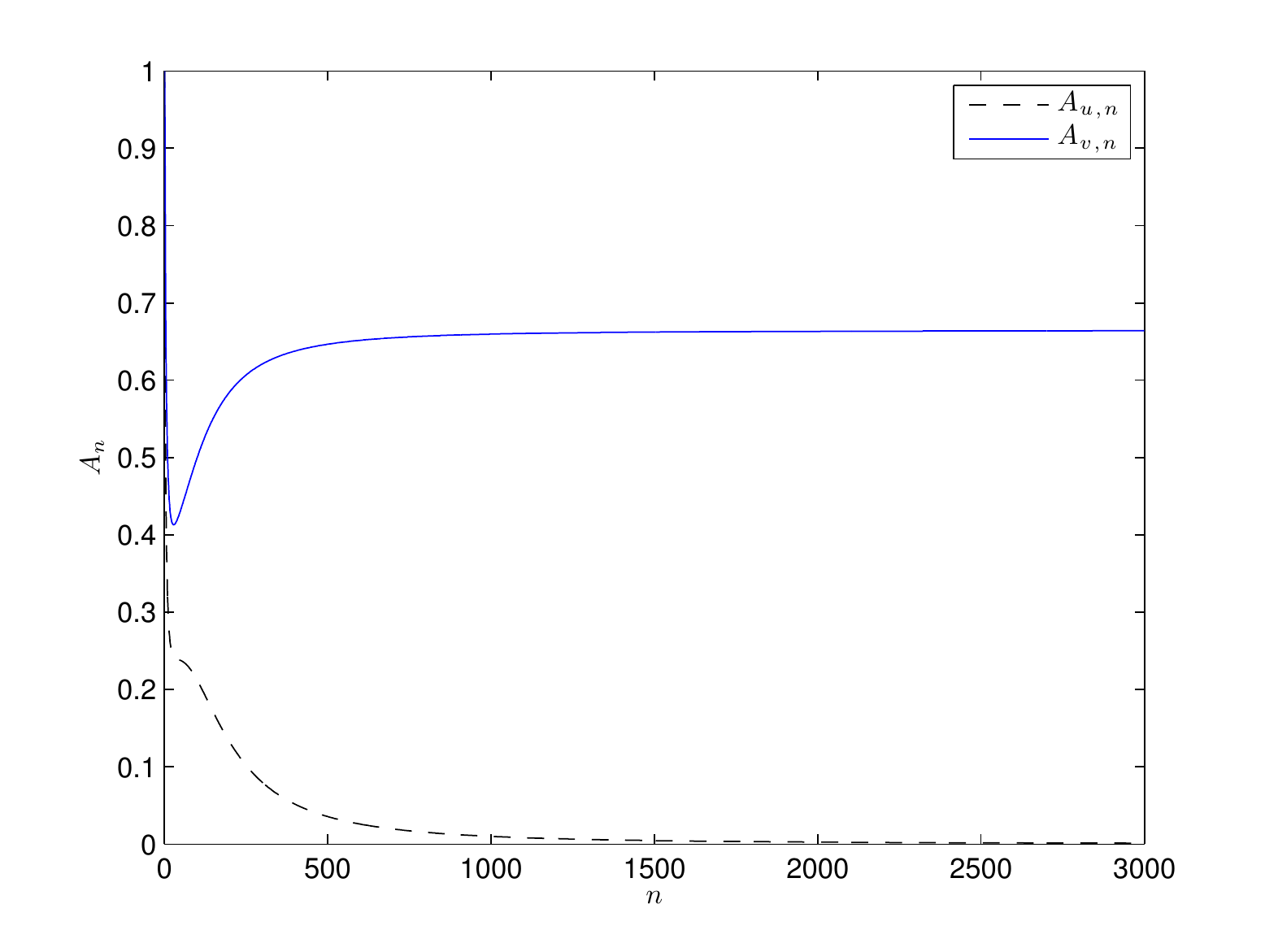} 
\caption{Computed scaling factors by the nRG procedure stated in Section \ref{sec:nrg} for the chemical reaction system. (a) $\alpha_{1, n}$ and $\alpha_{2, n}$ (b) $A_{u, n}$ and $A_{v, n}$.}
\label{fig:An_and_alpha}
\end{figure}

We conduct a nRG experiment using the power-law scaling, Eq. (\ref{eq:similarities}), for the above chemical reaction problem with the parameters, $p=q=1.5$, $d=0.75$, and $L=1.2$. The initial data are
\beq
u(x,0) = v(x,0) = \chi_{[-\ell,\ell]}(x) =\begin{cases}
1,  & -\ell \le x \le \ell, \\
0,  &\text{else}.
\end{cases}
\eeq
We choose $\ell = 0.5$ and the computational domain to be $[-10, 10]$.  For these initial data, the total conserved mass is $A=2$. Figure \ref{fig:An_and_alpha}(a) is a plot for $\alpha_{1, n}$ and $\alpha_{2, n}$ versus $n$. From the figure, we expect that $\alpha_{1, n}$ and $\alpha_{2, n}$ both converge to 1/2 as $n\rightarrow\infty$, although the figure suggests that $\alpha_{1, n}$ may converge much slower than $\alpha_{2, n}$. Since $\beta_n=1/2$ for all $n$, from Eq. (\ref{eq:An_Bn}), the convergences of $\alpha_{1, n}$ and $\alpha_{2, n}$, leads to Eq. (\ref{eq:similarities}). Moreover, Figure \ref{fig:An_and_alpha}(b) are the computed $A_{u, n}$ and $A_{v, n}$.  As expected, $A_{u, n}$ (the dashed-line) approaches to 0 as $n\rightarrow\infty$, while $A_{v, n}$ approaches to a constant $A_{v}$. Note that from Eqs. (\ref{eq:log_decay}) and (\ref{eq:similarities}), fo $t=L^n$, we have
\begin{equation}
A_{v, n} \phi_{v} \rightarrow A_{v} \phi_v = A\phi_d.
\end{equation}
Since $\phi_{v}= \sqrt{4\pi d}\, \phi_d$, $\sqrt{4\pi d}\,A_{v} = A$, or  $A_{v}=\displaystyle\frac{1}{\sqrt{4\pi d}} A$. For $d=0.75$, $A=2$, $A_{v}\approx 0.6515$. i.e. $A_{v,n }\rightarrow A_{v}\approx 0.6515$, and this is exactly what we observe in Figure \ref{fig:An_and_alpha}(b).

\begin{figure}[bhtp]
\centering
(a) \includegraphics[width=2.8in]{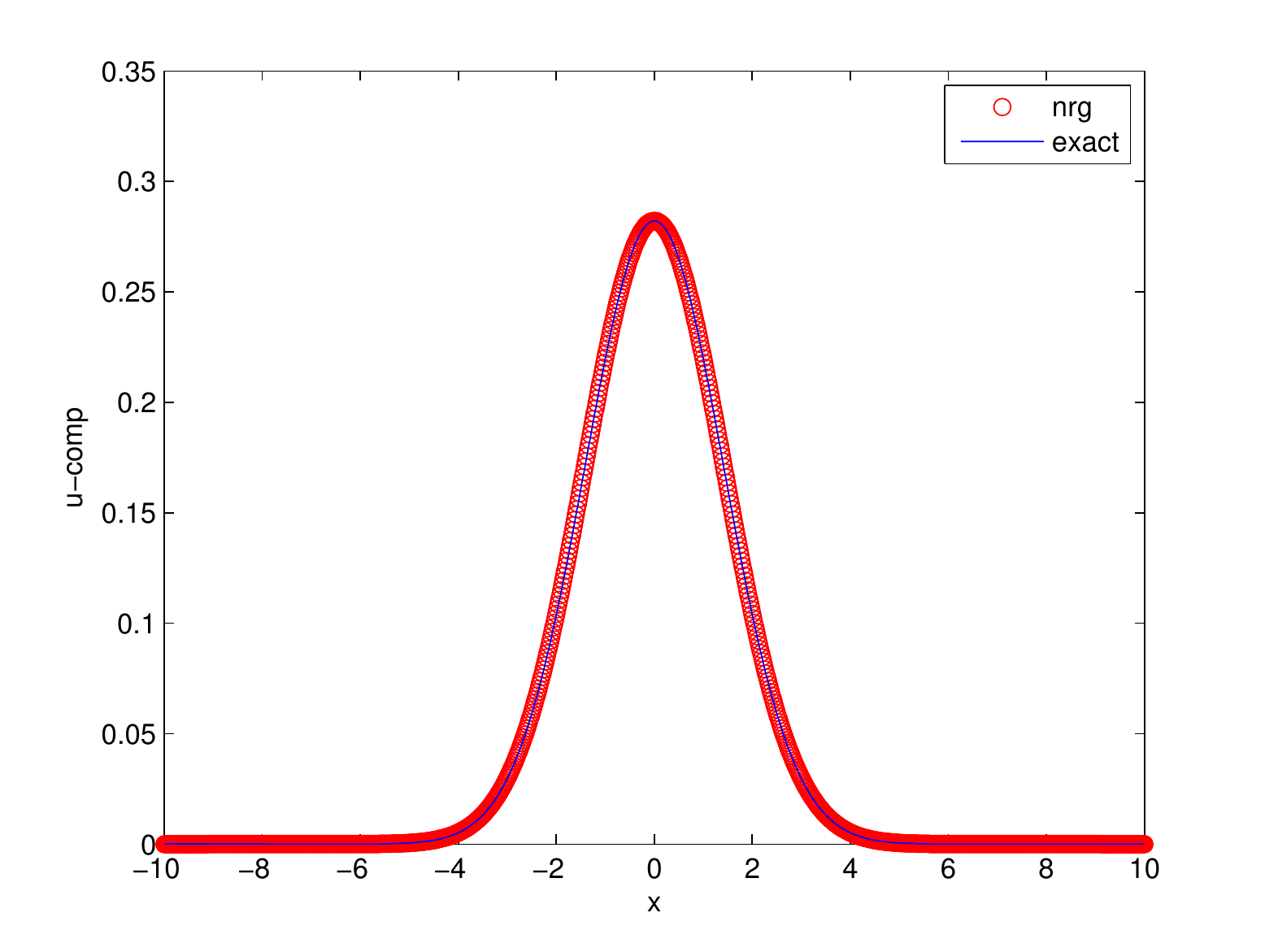} 
(b) \includegraphics[width=2.8in]{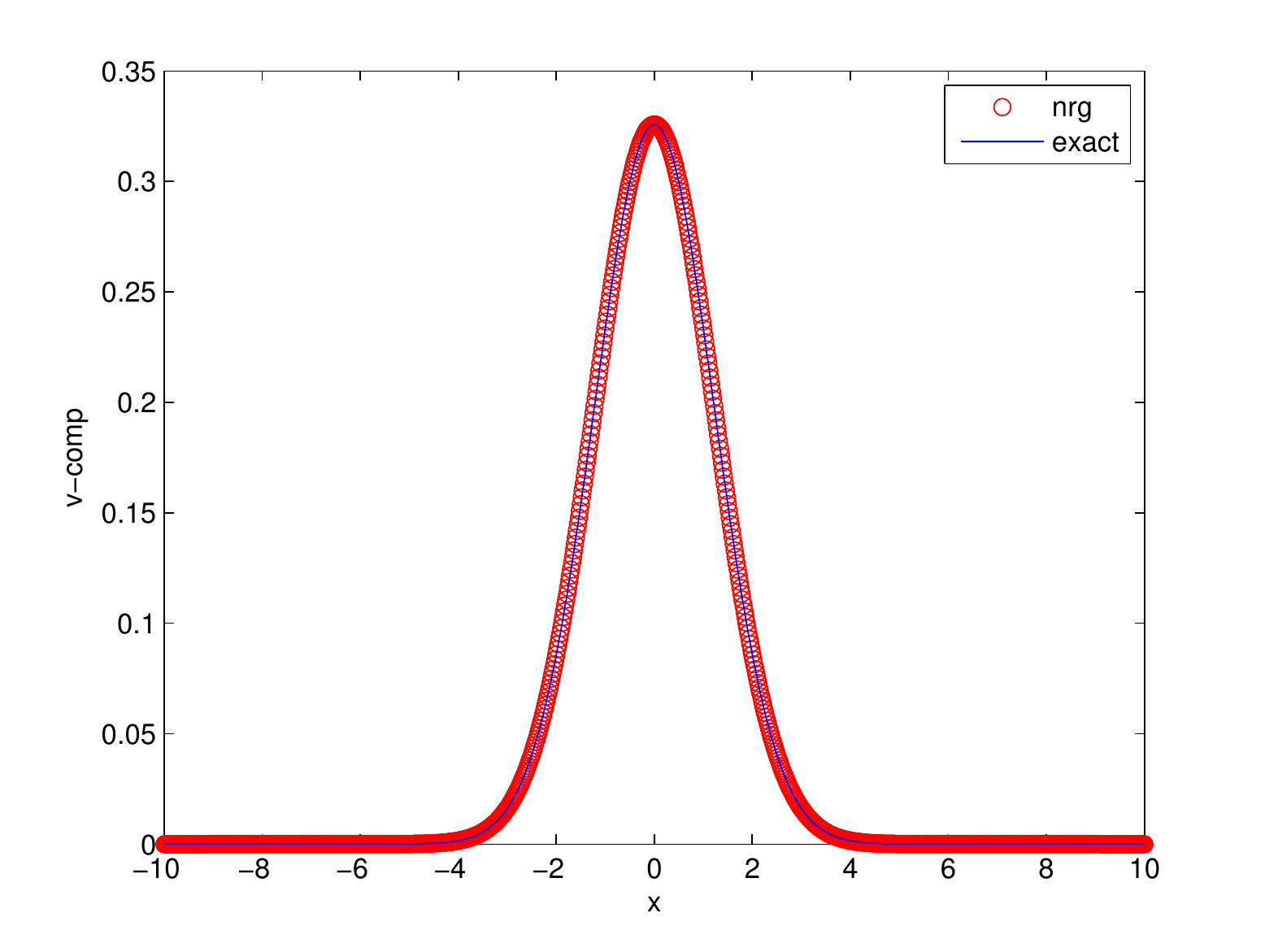} 
\caption{Comparison between the computed Gaussian similarity profile by Algorithm 1 and the predicted theoretical profile in \cite{LQ03} at $n=3000$, after adjusting the amplitudes.  (a) $u$-component, (b) $v$-component.} 
\label{fig:Gaussian_profile}
\end{figure}

Figure \ref{fig:Gaussian_profile} shows the comparison between the computed Gaussian similarity  profile and the predicted theoretical profile in \cite{LQ03} at $n=3000$, after adjusting the amplitudes. Both components correctly match the prediction, even though the hidden logarithmic decay in $u$ is not captured by the nRG algorithm.

\begin{figure}[bhtp]
\centering
\includegraphics[width=2.8in]{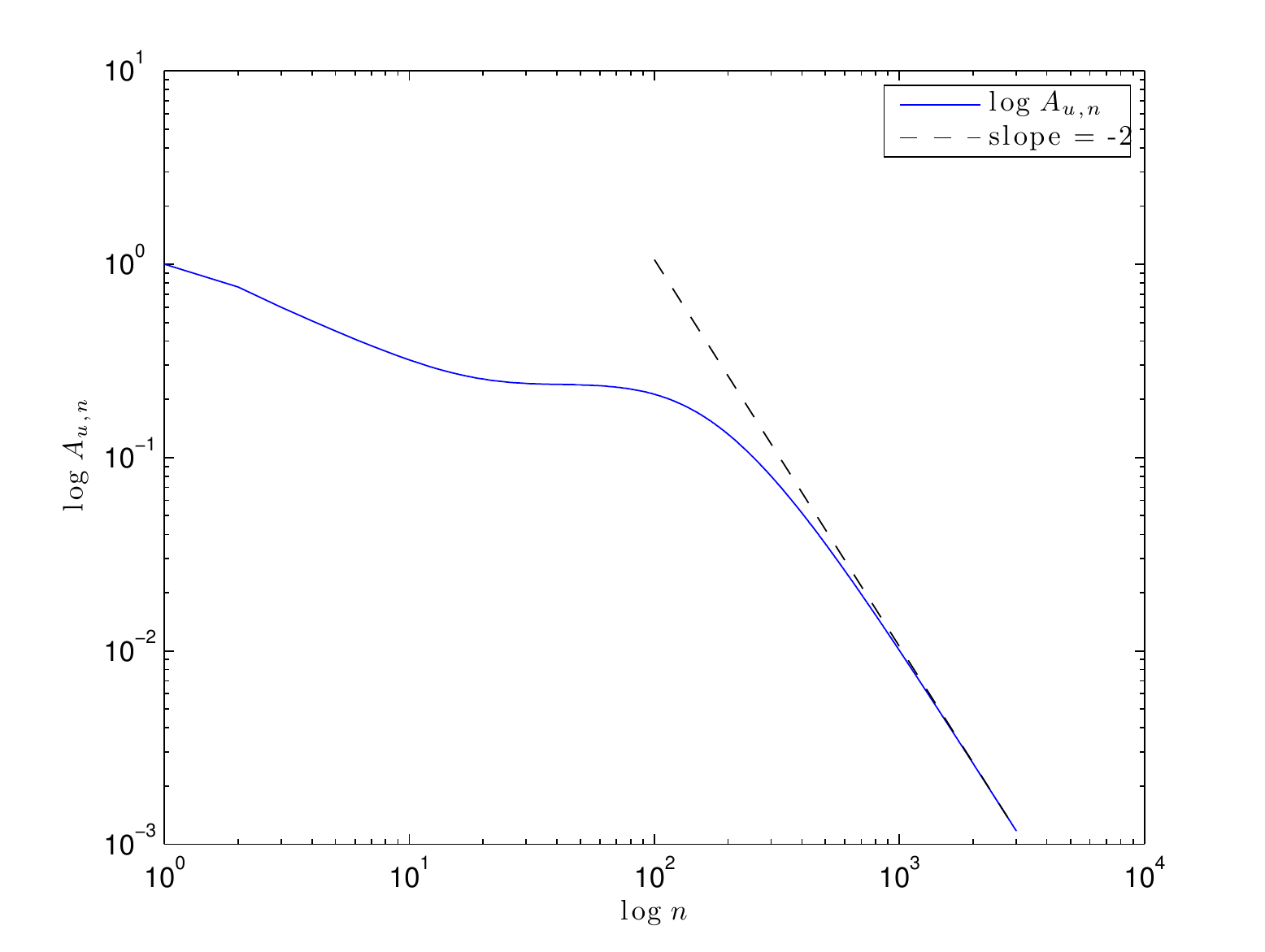} 
\caption{log$-$log plot for $A_{u, n}$ and $n$. The dashed-line is the straight line whose slope is equal to $-2$. }
\label{fig:An_no_log_decay}
\end{figure}

Now let's turn our attention to $A_{u, n}$. The fact that $A_{u, n}\rightarrow 0$ as $n\rightarrow\infty$ indicates that there was a ``hidden'' decay factor that was not captured by the current nRG procedure.  Taking a log-log plot for  $A_{u, n}$ and $n$,  Figure \ref{fig:An_no_log_decay} shows that  $\log A_{u, n} = (-2) (\log n) + \log C$, as $n\rightarrow\infty$. If we suppose that the hidden decay factor is related to $\log t$, then we could choose $C$ to be 
$C=A(\log L)^{-2}$, this results in $A_{u, n} = A(\log L^n)^{-2}$, and thus Eq. (\ref{eq:similarities}) becomes 
\beq\label{eq:log-decay-nrg}
\begin{split}
u(x, L^{n}) &\sim \frac{A}{L^{n/2}(\log L^n)^{2}}\phi_{u}(\frac{x}{L^{n/2}}),\\
v(x, L^{n}) &\sim \frac{A_{v}}{L^{n/2}}\phi_{v}(\frac{x}{L^{n/2}}).
\end{split}
\eeq
Eq. (\ref{eq:log-decay-nrg}) is evidently the similarity (asymptotic) solutions of the system of chemical-reaction equations for $t=L^{n}$ and $p=q=3/2$ (cf. Eq. (\ref{eq:log_decay})).

\section{Modified RG algorithm for logarithmic decay}\label{sec:AM2}

The above experiment suggests that the component $v$ has the decay factor $\sqrt{t}$,  but the component $u$ may have more than one decay factor. Without the asymptotic formula (\ref{eq:log_decay}), in principle, we do not know what the decay factors are. However, if we ``guess'' that one of them is also $\sqrt{t}$ based on Figures \ref{fig:An_and_alpha} and \ref{fig:An_no_log_decay}, and suppose that the other is related to $\log t$ with some unknown power $\gamma$, then
%
%
the solution in the asymptotic region gives us an idea how to compute the power $\gamma_n$ at the end of each iteration. Taking the hint from Eq. (\ref{eq:log_decay}), at times $t$ and $Lt$, the ratio of the solutions is
\beq\label{eq:modified-RG0}
 \frac{||u(x,t)||_{\infty}}{||u(x, Lt)||_{\infty}} = \frac{L^{1/2}t^{1/2}(\log Lt)^{\gamma}}{t^{1/2}(\log t)^{\gamma}}=L^{1/2}\left(\frac{(\log Lt)}{(\log t)}\right)^{\gamma}
\eeq
To modify the nRG procedure for this case, we observe that at the end of the $(n-1)^{th}$ iteration $t=L^{n-1}$, 
\beq\label{eq:modified-RG1}
L^{1/2}\left(\frac{\log L^{n}}{\log L^{n-1}}\right)^{\gamma_{n}} =  L^{1/2}\left(\frac{n}{n-1}\right)^{\gamma_{n}} =  \frac{||u_{n-1}(\cdot,1)||_{\infty}}{||u_{n-1}(\cdot,L)||_{\infty}},\quad n > 1,
\eeq
following Eq. (\ref{eq:modified-RG0}). 
%
%
Note that for the case $n=1$, $\gamma_{1}$ is computed by the power-law scaling,
\beq\label{eq:n=0-gamma}
L^{\gamma_{1}} = \frac{||u_0(\cdot,1)||_{\infty}}{||u_0(\cdot,L)||_{\infty}}.
\eeq
Eqs. (\ref{eq:modified-RG1}) and (\ref{eq:n=0-gamma}) suggest that for the $u$-component, the initial condition for the next iteration is set by
\beq\label{eq:modified-RG3}
u_{n}(x, 1)= L^{1/2}\left(\frac{n}{n-1}\right)^{\gamma_{n}}u_{n-1}(L^{1/2}x, L),\quad \text{for}\,\,n > 1,
\eeq
and
\beq\label{eq:modified-RG4}
u_{1}(x, 1)=L^{\gamma_{1}}u(L^{1/2}x, L),\quad \text{for}\,\, n=1.
\eeq
Here we have chosen $\beta_1=\beta_2=1/2$ in order to keep the diffusion coefficients unchanged. Note that from Eq. (\ref{eq:modified-RG3}), at the end of the $n^{th}$ iteration ($n > 1$), the iterative solutions $u_{n}$ and $v_{n}$ are related to the solutions of the PDE's by 
\beq\label{eq:u_v_sys}
\begin{split}
u_{n}(x,t)&=L^{\gamma_{1}+(n-1)/2}\,\displaystyle\overset{n}{\underset{k=2}{\Pi}}\left(\frac{k}{k-1}\right)^{\gamma_{k}}\,u(L^{n/2}x, L^{n}t),\\
v_{n}(x,t)&= L^{n\bar{\alpha}_{2, n}}v(L^{n/2}x, L^{n}t),
\end{split}
\eeq
where $\bar{\alpha}_{2, n} = \left(\alpha_{2,1}+\cdots+\alpha_{2,n}\right)/n$. Eq. (\ref{eq:u_v_sys}) implies that 

\beq\label{eq:u_v_sys_inv}
\begin{split}
u(x,t)&=L^{-\gamma_{1}-(n-1)/2}\,\displaystyle\overset{n}{\underset{k=2}{\Pi}}\left(\frac{k-1}{k}\right)^{\gamma_{k}}\,u_n(L^{-n/2}x, L^{-n}t),\\
v(x,t)&= L^{-n\bar{\alpha}_{2, n}}v_n(L^{-n/2}x, L^{-n)}t).
\end{split}
\eeq
Hence for the $n^{th}$ iteration ($n \ge 1$), the scaled system of PDEs for $u_n$ and $v_n$ is 
\beq\label{eq:un-sys}
\begin{split}
(u_{n})_t &= (u_{n})_{xx} - L^{n}\,L^{(-p+1)\gamma_{1}}(L^{-(n-1)/2})^{(p-1)}\,\left(\displaystyle\overset{n}{\underset{k=2}{\Pi}}\left(\frac{k-1}{k}\right)^{\gamma_{k}}\right)^{p-1}\, (L^{-n\bar{\alpha}_{2, n}})^{q}(u_n)^{p}(v_n)^{q},\\
(v_{n})_t &= d\,(v_{n})_{xx} + L^{n}\,L^{-p\gamma_{1}}(L^{-(n-1)/2})^{p}\,\left(\displaystyle\overset{n}{\underset{k=2}{\Pi}}\left(\frac{k-1}{k}\right)^{\gamma_{k}}\right)^{p}\, (L^{-n\bar{\alpha}_{2, n}})^{q-1}(u_n)^{p}(v_n)^{q},
\end{split}
\eeq
for $\beta_1=\beta_2=1/2$.  For $n=0$, the unscaled equation (\ref{eq:chem-reaction}) is solved.

%

Now, similar to the steps from Eq. (\ref{eq:rg_similarity_Ln}) to Eq. (\ref{eq:An_Bn}), we can define a variable $A_{u,n}$ for the $u$-component, so that we can monitor $A_{u, n}$ for convergence.  From the first equation in Eq. (\ref{eq:log-decay-nrg}) and the first equation in Eq. (\ref{eq:u_v_sys}), we have 
\beq\label{eq:An_sys}
\begin{split}
u(L^{n/2}x, L^{n}) & \sim L^{-n/2}(\log L^{n})^{-\gamma}A\phi(x),\\
u_n(x,1) & =L^{\gamma_{1}+(n-1)/2}\,\displaystyle\overset{n}{\underset{k=2}{\Pi}}\left(\frac{k}{k-1}\right)^{\gamma_{k}}\,u(L^{n/2}x, L^n).
\end{split}
\eeq
If we let $A_{*}=A(\log L)^{-\gamma}$ and assume that $\gamma_n\rightarrow\gamma$ as $n\rightarrow\infty$, Eq. (\ref{eq:An_sys}) implies
\beq\label{eq:An_sys2}
\begin{split}
u_n(x, 1) & \sim L^{\gamma_{1}+(n-1)/2}\,\displaystyle\overset{n}{\underset{k=2}{\Pi}}\left(\frac{k}{k-1}\right)^{\gamma_{k}}\,L^{-n/2}n^{-\gamma_{n}}A_{*}\phi(x)\\
& \sim L^{\gamma_{1}-1/2}\displaystyle\overset{n}{\underset{k=2}{\Pi}}\left(\frac{k}{k-1}\right)^{\gamma_{k}-\gamma_{n}}A_{*}\phi(x),
\end{split}
\eeq
The above equation holds, since  $\gamma_{n}\rightarrow\gamma$ for $n\rightarrow\infty$ and thus $n^{-\gamma_{n}}=\displaystyle\overset{n}{\underset{k=2}{\Pi}}\left(\frac{k}{k-1}\right)^{-\gamma_{n}}$ for $n\rightarrow\infty$. Eq. (\ref{eq:An_sys2})  is equivalent to 
\beq
L^{1/2-\gamma_1}\displaystyle\overset{n}{\underset{k=2}{\Pi}}\left(\frac{k}{k-1}\right)^{\gamma_n-\gamma_k}u_n(x, 1) \sim A_{*}\phi(x).
\eeq
If we define 
\beq
A_{u, n}= L^{1/2-\gamma_{1}}\displaystyle\overset{n}{\underset{k=2}{\Pi}}\left(\frac{k}{k-1}\right)^{\gamma_{n}-\gamma_{k}},\quad n > 1,
\eeq
we expect that $A_{u, n}\rightarrow A_{*}$ for $n$ large enough, provided $u_n(x, 1)\rightarrow \phi(x)$. 
Since $\phi=\sqrt{4\pi}\phi_1$, where $\phi_1$ is the Gaussian function in Eq. (\ref{eq:gaussian}) with $d=1$, this implies that  $A\sqrt{4\pi}=B$, where $B$ is the theoretical prediction in Eq. (\ref{eq:B}).  Therefore
\begin{equation}\label{eq:Astar}
A_{u, n}\rightarrow A_{*}=A(\log L)^{-\gamma} = \frac{B (\log L)^{-\gamma}}{\sqrt{4\pi}}.
\end{equation}
For $A_{v, n}$ we expect $A_{v,n}\rightarrow A_v = \displaystyle\frac{A}{\sqrt{4\pi}}$, where $A$ is the conserved total mass, the same as before.
We summarize the modified nRG procedure for the chemical reaction problem with the choice of parameters $\beta_1=\beta_2=1/2$ in Algorithm 2.

\begin{algorithm}[t]
\begin{algorithmic}
\label{alg2:nrg}
\For {$n=0,1, 2,\ldots,$ until convergence}
\begin{enumerate}
\item[1.]  Start with the IVP (\ref{eq:chem-reaction}) for $n=0$. Evolve $u_n$ and $v_n$ from $t=1$ to $t=L$, using the IVP (\ref{eq:un-sys})  for $n \ge 1$. 
\item[2.] Compute  $\gamma_{n}$ for the $u$-component by
\beq\label{eq:modified-RG-u-comp}
\begin{split}
&L^{\gamma_{1}} = \frac{||u_0(\cdot,1)||_{\infty}}{||u_0(\cdot,L)||_{\infty}}, \\
&L^{1/2}\left(\frac{n}{n-1}\right)^{\gamma_{n}} =  \frac{||u_{n-1}(\cdot,1)||_{\infty}}{||u_{n-1}(\cdot,L)||_{\infty}},\quad n \ge 2.
\end{split}
\eeq

Compute $\alpha_{2, n}$ for the $v$-component by 
\begin{equation*}
L^{\alpha_{2, n}} = \frac{||v_{n-1}(\cdot,1)||_{\infty}}{||v_{n-1}(\cdot,L)||_{\infty}}.
\end{equation*}

\item[3.] Compute $A_{u, n}= L^{1/2-\gamma_{1}}\displaystyle\overset{n}{\underset{k=2}{\Pi}}\left(\frac{k}{k-1}\right)^{\gamma_{n}-\gamma_{k}}$, $n > 1$;  $A_{v, n}=L^{n(\alpha_{2, n}-\bar{\alpha}_{2, n})}$, where $\bar{\alpha}_{2, n} = \left(\alpha_{2,1}+\cdots+\alpha_{2,n}\right)/n$. 

\item[4.] Set initial data for the next iteration by $f_{u, n+1} = L^{1/2}\left(\frac{n}{n-1}\right)^{\gamma_{n}}u_n\left(L^{1/2}x,L\right)$ and $f_{v, n+1}(x)=L^{\alpha_{2, n}}v_n\left(L^{1/2}x,L\right)$.
\end{enumerate}
\EndFor
\end{algorithmic}
\caption{The nRG procedure for the chemical reaction system}
\end{algorithm}

\section{A numerical experiment for the logarithmic decay}\label{sec:AM2_example}

\begin{figure}[bhtp]
\centering
(a) \includegraphics[width=2.8in]{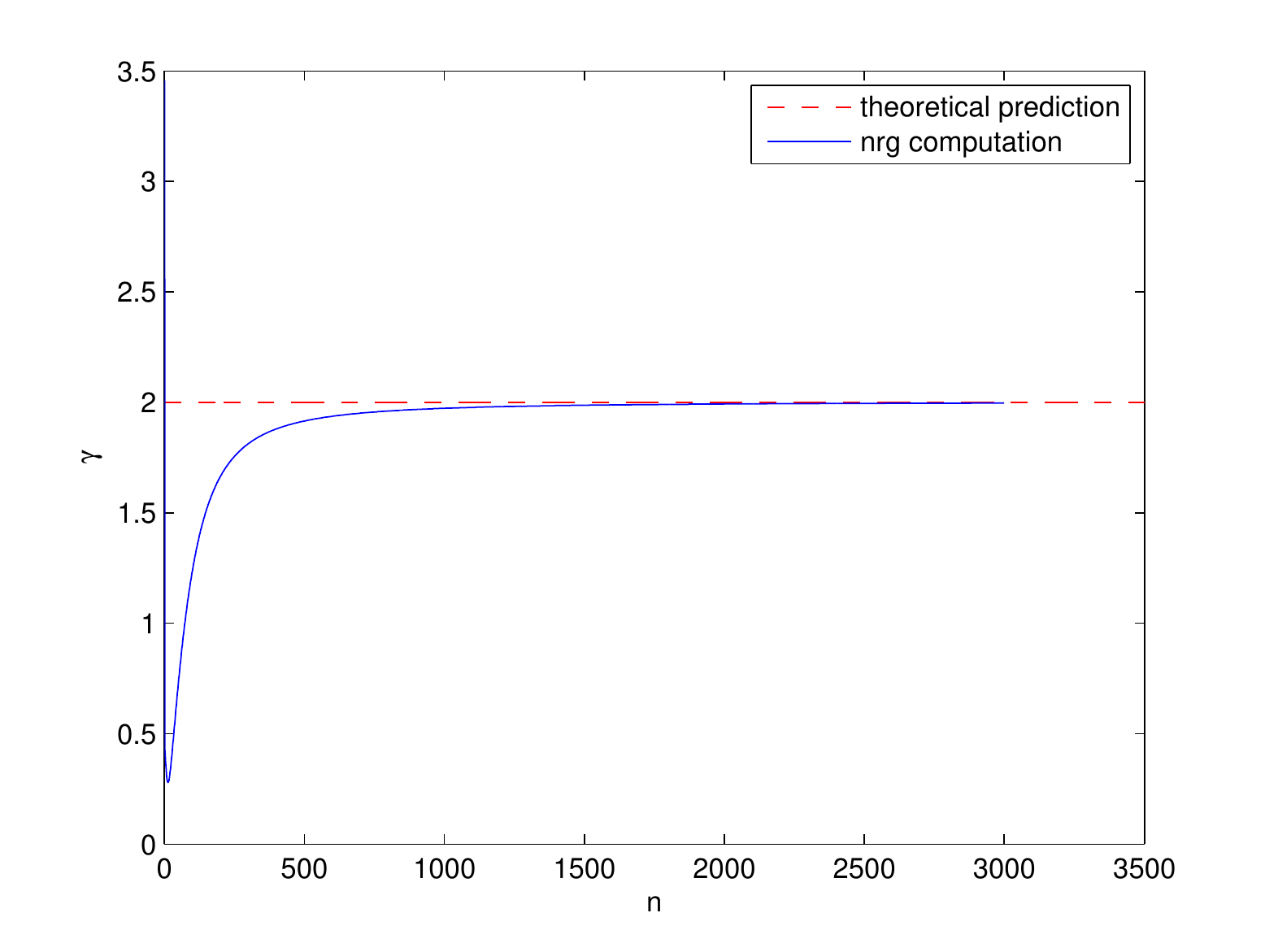} 
(b) \includegraphics[width=2.8in]{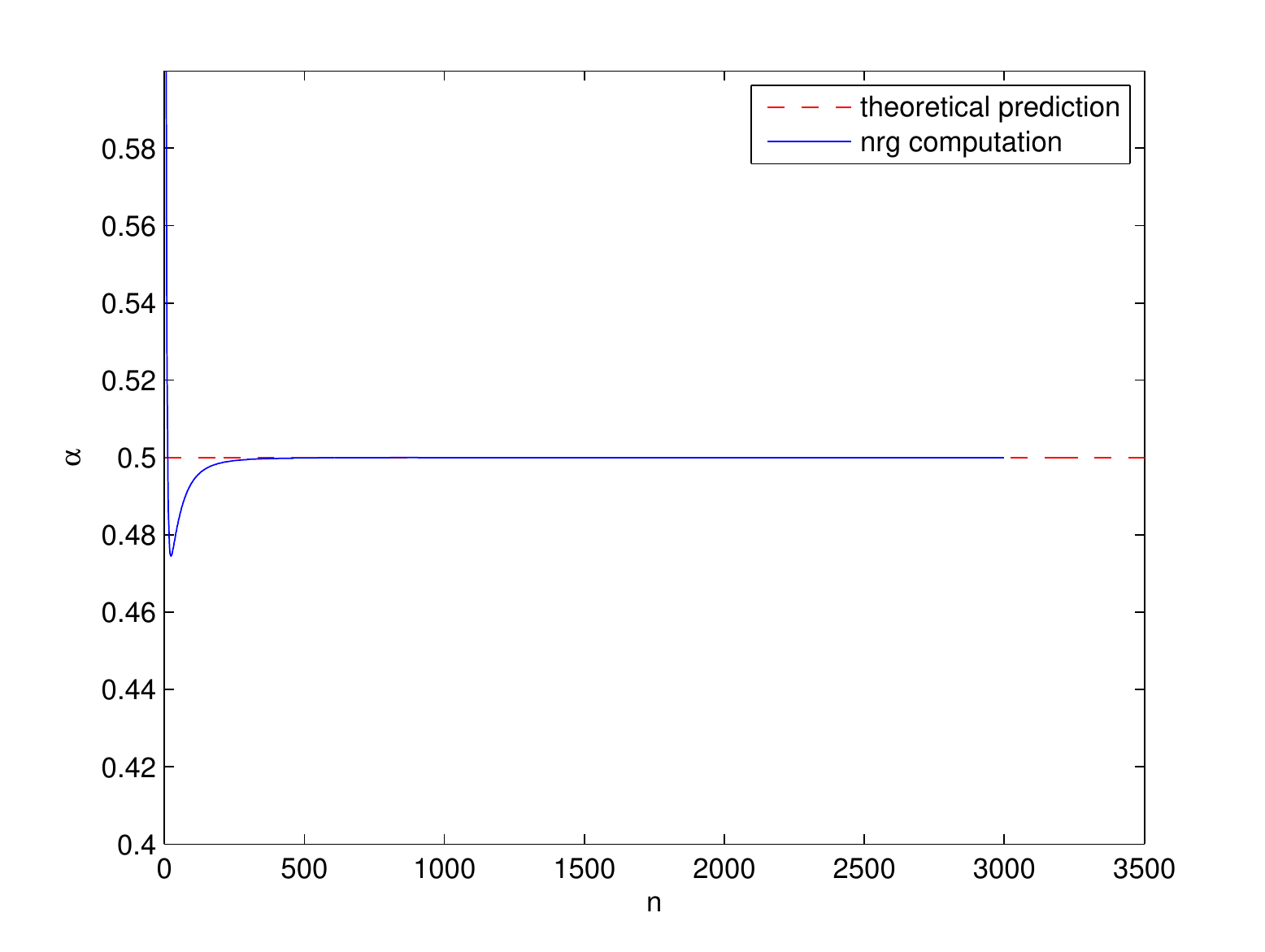} 
\caption{Comparison of the theoretical prediction and the nRG computation for the exponents: (a) the sequence of logarithmic decay exponent $\gamma_n\rightarrow\gamma =2$ for $u$, and (b) the sequence of the power law decay exponent $\alpha_n\rightarrow \alpha =0.5$ for $v$.}
\label{fig:exponents}
\end{figure}

\begin{figure}[bhtp]
\centering
(a) \includegraphics[width=2.8in]{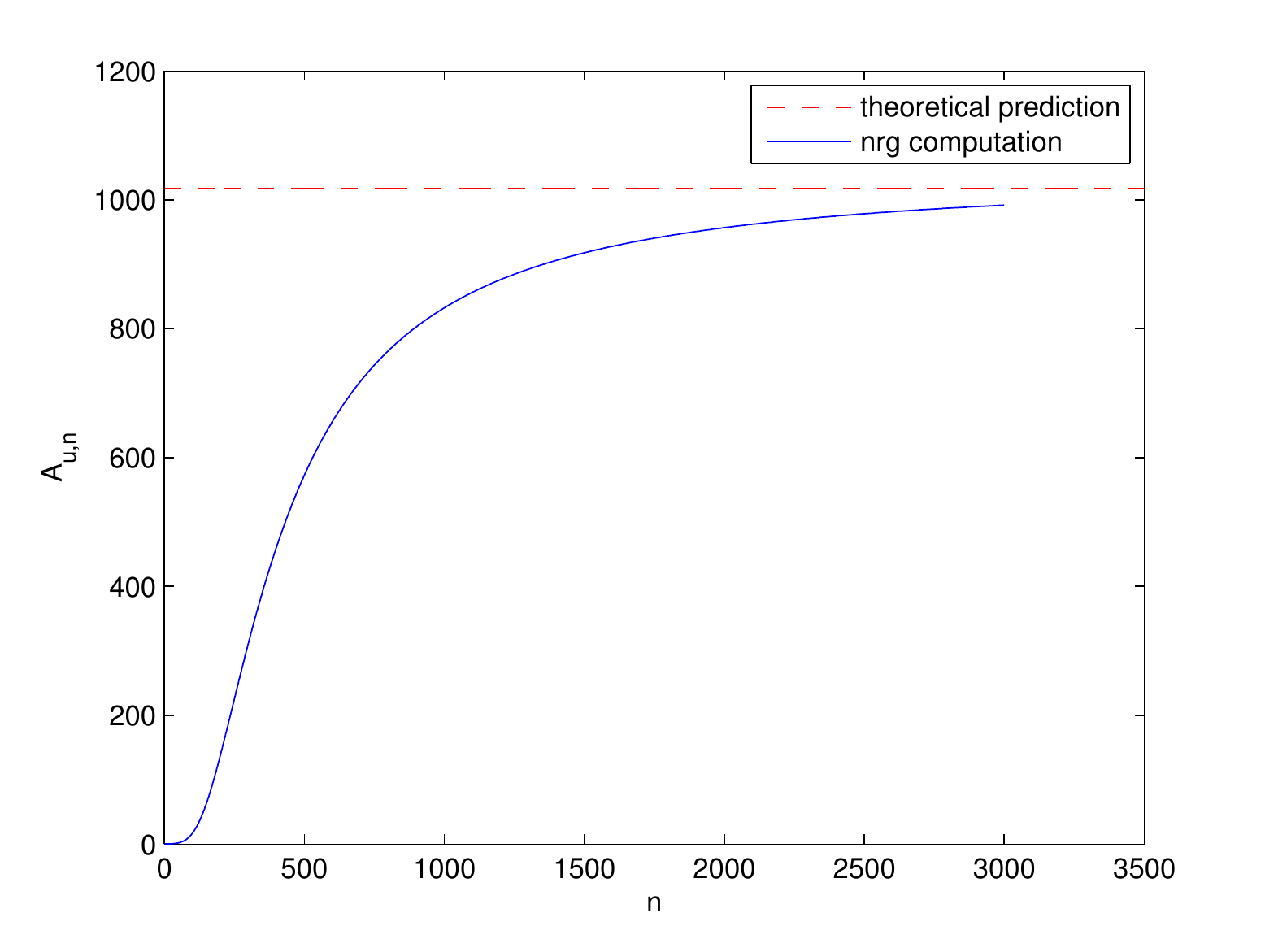} 
(b) \includegraphics[width=2.8in]{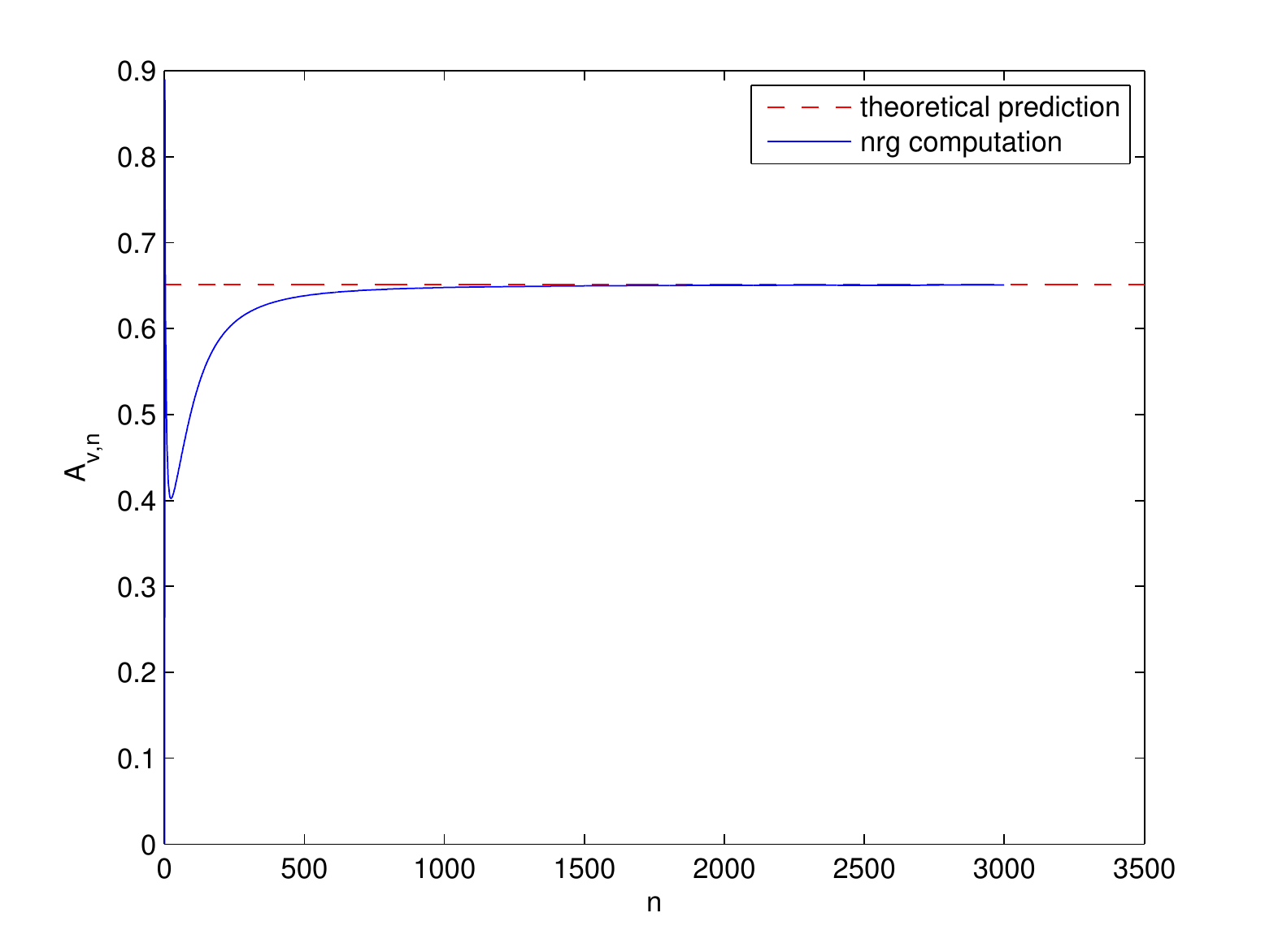} 
\caption{Comparison of the theoretical prediction and the nRG computation for the pre factors: (a)  $A_{u, n}\rightarrow A_{*} \approx 1016.89$, and (b) $A_{v, n}\rightarrow A_v\approx 0.6515$. }
\label{fig:prefactors}
\end{figure}

To illustrate that the modified RG algorithm accurately captures the hidden logarithmic decay, we apply Algorithm 2 to the diffusion-reaction equations  Eq. (\ref{eq:chem-reaction}) with $p=q=3/2$. For this set of parameters, the asymptotical behavior of the solutions follows Eqs. (\ref{eq:log_decay}) and (\ref{eq:B}). Thanks to the choice of $\beta_1=\beta_2=1/2$, the diffusivities in Eq. (\ref{eq:un-sys}) are kept to be 1 and $d$, respectively. Here we choose $d=3/4$. The system of PDEs (\ref{eq:un-sys}) is discretized by an explicit second-order method (forward Euler  for the time derivative and the second-order center difference for the spatial second derivative). The parameters used in our nRG algorithm are $L=1.25$, $x\in [-10, 10]$, $dx= 0.04$, and $dt=0.00025$. The number of iterations for nRG is 3000. The theoretical prediction for the critical exponents is $\gamma =2$ (power of the logarithmic decay for $u$) and $\alpha = 0.5$ (power of the power law decay for $v$). Figure \ref{fig:exponents} shows that Algorithm 2 accurately captures these two exponents. At the mean time, the theoretical prediction for the pre-factor $A_{*}$ is $A_{*}\approx 1016.89$, computed by Eq. (\ref{eq:Astar}), and the pre-factor $A_{v}$ is $A_{v}\approx 0.6515$, the same as our previous calculation, and we observe from Figure \ref{fig:prefactors} that both $A_{u, n}$ and $A_{v, n}$ numerically converges to their theoretical values, respectively.

\begin{figure}[bhtp]
\centering
(a) \includegraphics[width=2.8in]{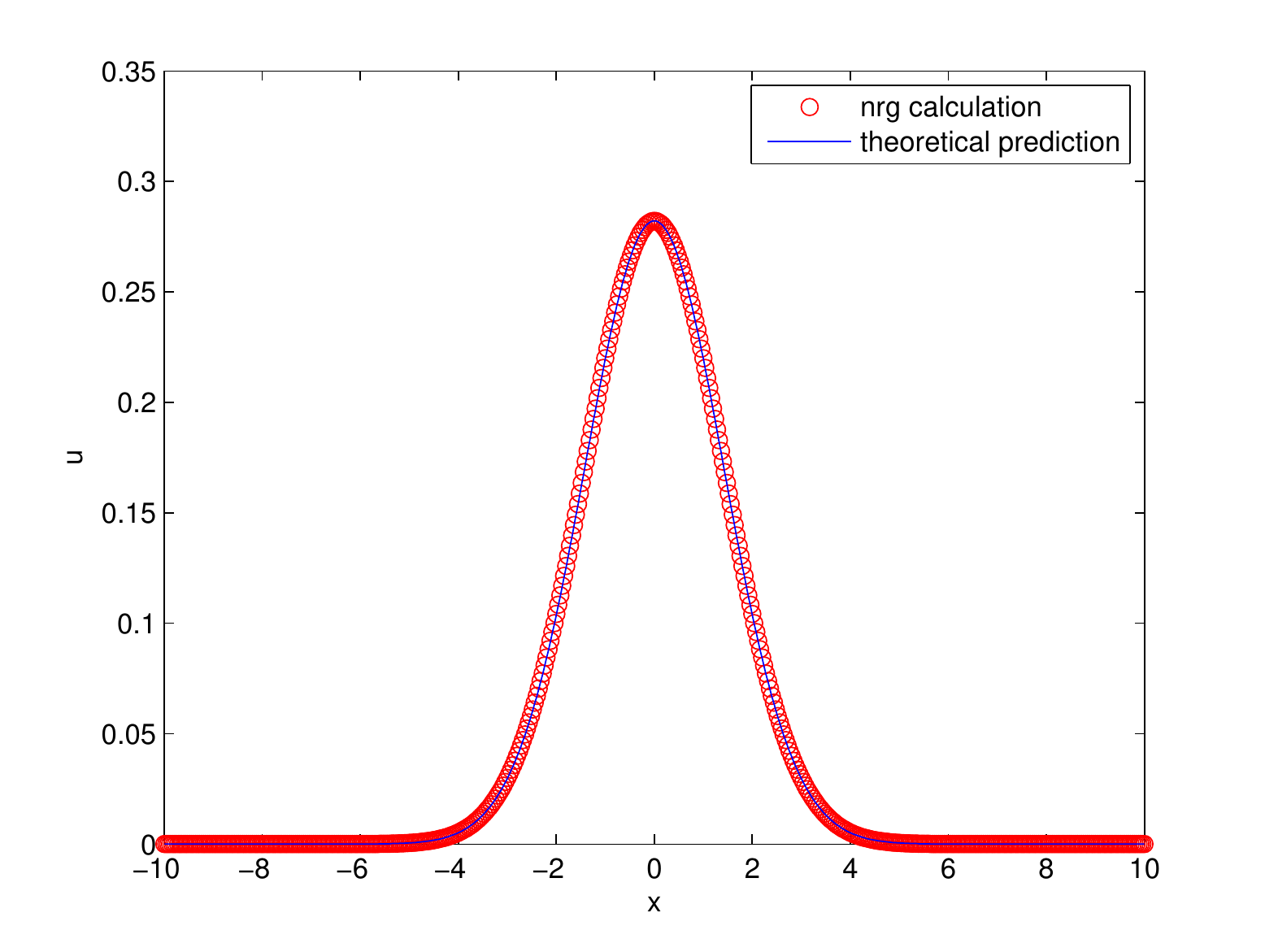} 
(b) \includegraphics[width=2.8in]{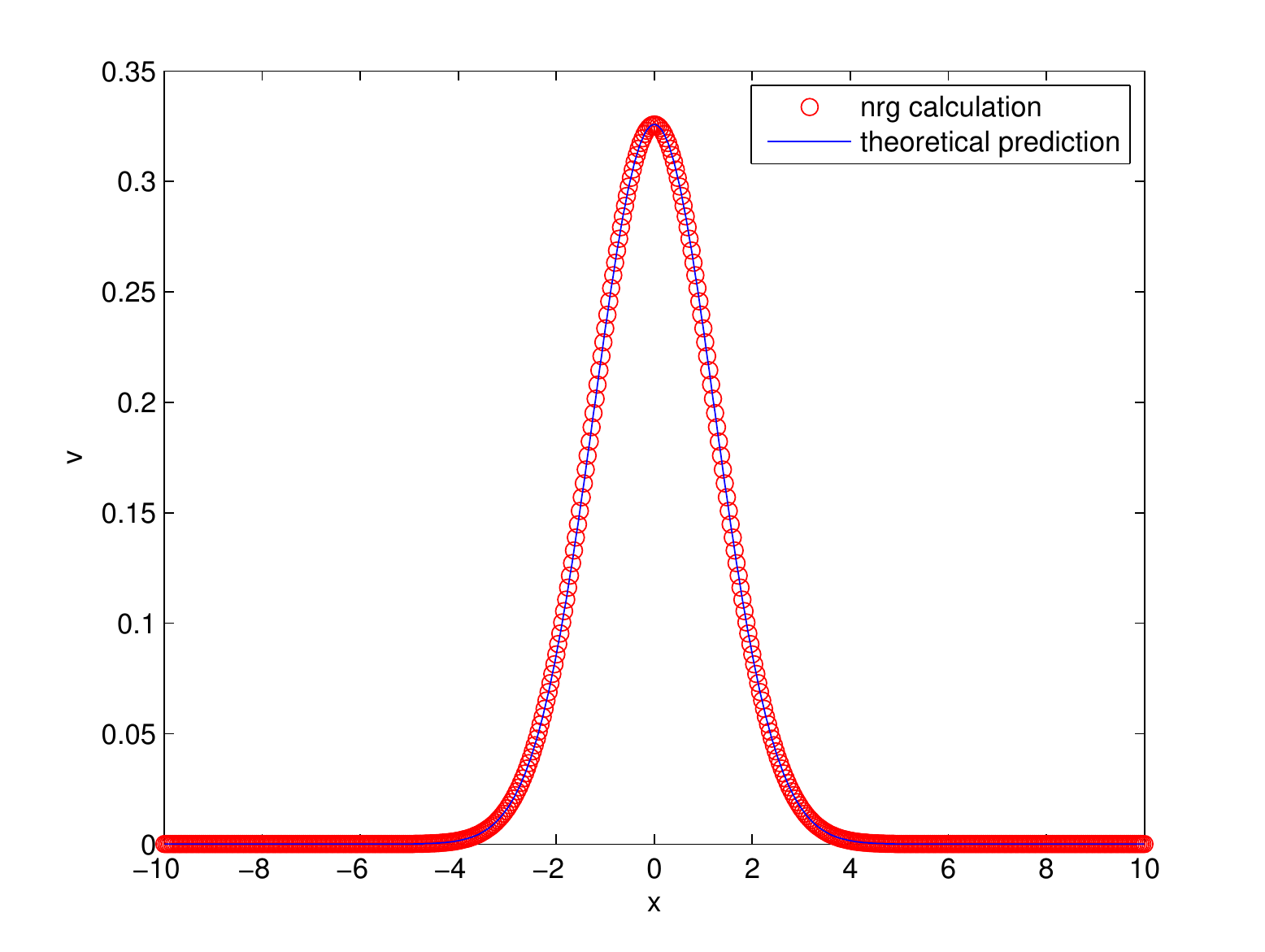} 
\caption{Comparison between the computed Gaussian similarity profile by Algorithm 2 and the predicted theoretical profile in \cite{LQ03} at $n=3000$, after adjusting the amplitudes.  (a) $u$-component, (b) $v$-component.} 
\label{fig:Gaussian_profile_mnrg}
\end{figure}

Finally, in the the previous section \ref{sec:no-decay}, our calculation suggests that the original nRG algorithm captures the power law exponents and produces the final similarity profiles that match the theoretical prediction in \cite{LQ03} without taking into account the logarithmic decay. In this experiment, we use the hint from the previous calculation to assume the exponent of the power law decay for the $u$ component. We modified the RG algorithm to include the logarithmic decay. The modified RG algorithm captures the critical exponents and render the numerically convergent pre-factors for both components. It remains to show whether the similarity profiles produced by the modified RG algorithm match the theoretical prediction. Sure enough, Figure \ref{fig:Gaussian_profile_mnrg} shows that the modified RG algorithm produces the similarity profiles that match the theoretical prediction exactly after we adjust the amplitudes by multiplying the factor $1 / \sqrt{4\pi d}$, where $d=1$ for $u$ and $d=3/4$ for $v$, respectively.

\section{Concluding Remarks}

We have presented  and systematically examined a numerical procedure, based on the RG theory for PDEs, that renders the detailed and efficient computation of asymptotically self-similar dynamics in solutions of PDEs. The effectiveness and robustness of the nRG algorithms were illustrated through several examples of quasilinear and nonlinear PDEs combining diffusive, reactive and nonlinear propagation effects. It is worth noting that the modified RG algorithm presented in Sections \ref{sec:AM2} and \ref{sec:AM2_example} for the nonlinear system of cubic autocatalytic chemical reaction equations nicely responds to the remark made by Li and Qi \cite{LQ03}:
\begin{quote}
``The appearance of  $\log t$ indicates the analysis is more involved and subtle. In particular, it is well known in the scientific computation field that a scaling of $\log t$ is hardly detectable in computation."
\end{quote}
by detecting the extra decay and capturing the power of logarithmic decay.

We refer readers to \cite{Isaia_thesis} for some preliminary results of the calculations of multidimensional problems by using the similar numerical scaling strategy described in this paper. A proper modification of the described RG algorithm can be used to compute traveling waves and is currently under our investigation. We are also investigating the applicability of an adapted version of the RG algorithm to blow-up problems. We expect to report our results in the future.

\section{Acknowledgement}
GAB thanks the Department of Mathematics at the University of Wyoming for the hospitality which made possible this collaboration.  LL is partially supported by NSF DMS 1413273.

\appendix

\section{Data of the validation example of the Barenblatt's equation}\label{app:data1}

\begin{table}[htpbt]
\caption{}
\label{tab:validation}
\centering
\begin{tabular}{c|c|c|c}
\hline
\backslashbox{$\epsilon$}{$\alpha(\epsilon)$} & linear perturbative &quadratic perturbative &  computed  \\\hline
-0.9 &0.282226347932771 &0.230753893532771 &0.177628930214765  \\\hline 
-0.8 &0.306423420384685 &0.265753826784685 &0.237625966873617  \\\hline 
-0.7 &0.330620492836600 &0.299482835236600 &0.285818489281125\\\hline 
-0.6 &0.354817565288514 &0.331940918888514 &0.336450068396061 \\\hline 
-0.5 &0.379014637740428 &0.363128077740428 &0.369213156837354\\\hline 
-0.4 &0.403211710192343 &0.393044311792343 &0.399078625788562  \\\hline 
-0.3 &0.427408782644257 &0.421689621044257 &0.426763625391393\\\hline
-0.2 &0.451605855096171 &0.449064005496171 &0.452587853032610 \\\hline
-0.1 &0.475802927548086 &0.475167465148086 &0.476940935172314\\\hline
0 &0.500000000000000 &0.500000000000000 &0.499991367558985 \\\hline
0.1 &0.524197072451914 &0.523561610051914 &0.521936778063526  \\\hline
0.2 &0.548394144903829 &0.545852295303829 &0.542938675675171  \\\hline
0.3 &0.572591217355743 &0.566872055755743 &0.563094313618048\\\hline
0.4 &0.596788289807657 &0.586620891407657 &0.582429649916552  \\\hline
0.5 &0.620985362259572 &0.605098802259572 &0.601176410464571  \\\hline
0.6 &0.645182434711486 &0.622305788311486 &0.619283036774065  \\\hline
0.7 &0.669379507163400 &0.638241849563400 &0.636638537608389 \\\hline
0.8 &0.693576579615315 &0.652906986015315 &0.653808099549284  \\\hline
0.9 &0.717773652067229 &0.666301197667229 &0.670383569387994 \\\hline
1.0 &0.741970724519143 &0.678424484519143 &0.686217411435144 \\\hline

\end{tabular}
\end{table}

\pagebreak

\end{document}